\documentclass[12pt]{article}
\usepackage{amssymb}
\usepackage{mathrsfs}
\newtheorem{theorem}{Theorem}

\newtheorem{claim}[theorem]{Claim}
\newtheorem{proposition}[theorem]{Proposition}
\newtheorem{lemma}[theorem]{Lemma}
\newtheorem{corollary}[theorem]{Corollary}
\newtheorem{remark}[theorem]{Remark}
\newtheorem{remarks}[theorem]{Remarks}

\newtheorem{definition}[theorem]{Definition}

\newcommand{\R}{\mathbb{R}}
\newcommand{\Q}{\mathbb{Q}}

\newcommand{\Sf}{\mathbb{S}}

\newcommand{\spa}{\mbox{span\,}}
\newcommand{\hess}{\mbox{\em  Hess\,}}

\newcommand{\po}{{\hspace*{-1ex}}{\bf .  }}

\newcommand{\sus}{submanifolds }

\def\Ral{{\cal R}}
\def\P{{\cal P}}

\def\ral{{\cal R}}

\def\Fes{{\cal F}}
\def\Fes{{\cal F}}

\def\<{\langle}

\def\>{\rangle}
\def\a{\alpha}
\def\va{\varphi}

\def\o{\omega}
\def\O{\Omega}

\def\bea{\begin{eqnarray*} }
\def\eea{\end{eqnarray*} }
\def\be{\begin{equation} }
\def\ee{\end{equation} }

\def\proof{\noindent{\it Proof: }}
\def\qed{\ifhmode\unskip\nobreak\fi\ifmmode\ifinner\else
\hskip5 pt \fi\fi\hbox{\hskip5 pt \vrule width4 pt
height6 pt  depth1.5 pt \hskip 1pt }}
\begin{document}
\title{The vectorial Ribaucour transformation for submanifolds of constant sectional curvature
%\footnote {1991 Mathematics Subject Classification. Primary: ??}
}
\author{D. Guimar\~aes and R. Tojeiro}
\date{}
\maketitle

%\author{S. Canevari and R. Tojeiro}
%        \thanks{E-mail: marcos@impa.br, rovenski@math.haifa.ac.il and
%tojeiro@dm.ufscar.br}
%\\
%(1) IMPA, Brazil (2) University of Haifa, (3) UFSCar, Brazil}
\date{}
\maketitle

\begin{abstract}
We obtain a reduction of the vectorial Ribaucour transformation that preserves the class of submanifolds of constant sectional  curvature of space forms, which we call the $L$-transformation. It allows to construct a family of such submanifolds starting with a given one and a vector-valued solution of a system of linear partial differential equations. We prove a decomposition theorem for the $L$-transformation, which is a far-reaching generalization of the classical permutability formula for the Ribaucour transformation of surfaces of constant curvature in Euclidean three space. As a consequence, we derive  a Bianchi-cube theorem, which  allows to produce, from $k$ initial scalar $L$-transforms of a given  submanifold of constant  curvature,  a whole $k$-dimensional cube all of whose remaining $2^k-(k+1)$ vertices are submanifolds with the same constant sectional curvature given by explicit algebraic formulae. We also obtain  further reductions, as well as  corresponding decomposition  
 and Bianchi-cube theorems, for the classes of $n$-dimensional flat Lagrangian submanifolds of  $\mathbb{C}^n$ and 
$n$-dimensional Lagrangian submanifolds with constant  curvature $c$ of the complex projective space $\mathbb C\mathbb P^n(4c)$ 
or the complex hyperbolic space $\mathbb C\mathbb H^n(4c)$ of complex dimension $n$ and constant holomorphic curvature~4c.
\end{abstract}

The study of isometric immersions of space forms into space forms is a central topic in submanifold theory, having its origins in the study of surfaces $f\colon M^2(c)\to \Q^3(\tilde c)$ of constant Gauss curvature $c$ in three-dimensional space forms of constant sectional curvature~$\tilde c$. If $c\neq 0$ (respectively, $c=0$), the Gauss and Codazzi equations of such surfaces reduce to 
the sin-Gordon or sinh-Gordon equations (respectively, wave and Laplace equations), depending on whether $c<\tilde c$ or $c>\tilde c$,
respectively. Understanding the complicated structure of these equations has led to some well-known global nonexistence results as well as to beautiful transformation theories for the construction of local examples.

The algebraic structure of the second fundamental form of an isometric immersion $f\colon M^n(c)\to \Q^{n+p}(c)$ of higher dimension and codimension was investigated by E. Cartan by means of his theory of exteriorly orthogonal quadratic forms. He also studied 
isometric immersions $f\colon M^n(c)\to \Q^{n+p}(\tilde c)$ with $c<\tilde c$  by looking at their compositions $i\circ f$ with an umbilical inclusion
$i\colon \Q^{n+p}(\tilde c)\to \Q^{n+p+1}(c)$, and proved   that for such an isometric immersion one must have $p\geq n-1$, and that for $p=n-1$ the normal bundle is flat. 

Cartan's algebraic results were later extended to the dual case $c>\tilde c$ by Moore \cite{mo}.  If $p\leq n-2$, it turns out that the only possible structure for the second fundamental form 
is that of a composition $f=h\circ i$ of an umbilical inclusion $i\colon M^n(c)\to \Q^{n+1}(\tilde c)$ with an  isometric immersion $h\colon U\subset \Q^{n+1}(\tilde c)\to \Q^{n+p}(\tilde c)$ with $i(M^n(c))\subset U$. If $p=n-1$, exactly one further possibility arises, in which case the immersion has also flat normal bundle.

Flatness of the normal bundle is thus a natural condition for isometric immersions $f\colon M^n(c)\to \Q^{n+p}(\tilde c)$, 
and implies that $M^n(c)$ can be locally endowed with principal coordinates  with respect to which the integrability conditions
give rise to nonlinear partial differential equations that reduce to the wave, Laplace, sin and  sinh-Gordon equations when $n=2$, $p=1$ and $c\neq \tilde c$ (see Propositions \ref{thm:hiebetaij} and \ref{teo:curvdif} below).

For the construction of explicit examples of submanifolds of constant sectional curvature of o space forms with higher dimension and codimension,
a useful tool is the Ribaucour transformation  developed in \cite{dt4} (see also \cite{dt5}) as an extension of the classical 
Ribaucour transformation  for surfaces in Euclidean three space. It allows to produce  a family of new submanifolds 
of constant sectional curvature of a space form starting with a given one and a solution of a linear system of partial differential equations.

Aiming at a better understanding of the iteration of Ribaucour transformations for submanifolds, with an eye towards deriving a procedure for the  construction of all submanifolds with flat normal bundle,  a vectorial version of the transformation was developed in \cite{dft},
based on its version in \cite{lm} for orthogonal systems, which shed light on its permutability properties. 

In this paper we obtain  a reduction of the vectorial Ribaucour transformation that preserves the class of submanifolds of constant sectional curvature of space forms. It depends on a linear operator $L$ of a vector space $V$, so we call it the $L$- Ribaucour transformation, or simply the $L$-transformation.  In the scalar case, that is, when $V$ has dimension one, it reduces to the scalar Ribaucour transformation for submanifolds of constant sectional curvature studied in \cite{dt4}.
%, called the ${\cal R}_C$-transformation in that paper,  which in turn reduces to the Ribaucour transformation  of  constant curvature surfaces in $\R^3$ studied by Bianchi \cite{bi}. 
The $L$-transformation allows to construct a family of new submanifolds of constant sectional curvature starting with a given one and a vector-valued solution of a linear system of PDE's. Proving the existence of $L$-transforms of a  submanifold of constant sectional curvature with given initial conditions requires  looking for \emph{invertible} solutions of a  
certain system of Sylvester-type matrix equations.

We prove a decomposition theorem for the $L$-transformation, which is a far-reaching generalization of the  permutability formula for the Ribaucour transformation of surfaces of constant curvature in Euclidean three space. It implies that an $L$-transformation given by a \emph{diagonalizable} operator $L$ is the iterate of $n=\dim V$ scalar Ribaucour transformations of the  type considered in \cite{dt4}. In particular, we show that  $k$ such Ribaucour transforms of a given submanifold of constant sectional curvature give rise to a whole $k$-dimensional cube, all of whose remaining $2^k-(k+1)$ vertices  are submanifolds with the same constant sectional curvature given by means of explicit algebraic formulae. On the other hand, an $L$-transformation given by a non-diagonalizable linear operator $L$
yields new submanifolds of constant sectional curvature that can not be obtained by an iteration of a sequence of scalar Ribaucour transformations as those in \cite{dt4}.

We develop a further reduction of the $L$-transformation that preserves  the classes  of $n$-dimensional flat Lagrangian submanifolds of  $\mathbb{C}^n=\R^{2n}$ and $n$-dimensional  submanifolds with constant  curvature $c$ of  $\Sf_\epsilon^{2n+1}(c)$ that are horizontal with respect to the Hopf fibration $\pi\colon \Sf_\epsilon^{2n+1}(c)\to \tilde M^n(4c)$. Here $\Sf_\epsilon^{2n+1}(c)$ stands for either the standard Euclidean sphere or the anti-de-Sitter space time of dimension  $2n+1$ and constant  curvature $c$, corresponding to $\epsilon=1$ or $\epsilon =-1$, respectively, and $\tilde M^n(4c)$  denotes either the complex projective space $\mathbb C\mathbb P^n(4c)$ or complex hyperbolic space $\mathbb C\mathbb H^n(4c)$ of complex dimension $n$ and constant holomorphic sectional curvature 4c, corresponding to $c>0$ or $c<0$, respectively.  Horizontal $n$-dimensional  submanifolds with constant  curvature $c$ of  $\Sf_\epsilon^{2n+1}(c)$  
%with respect to the Hopf fibration 
project down to $n$-dimensional 
Lagrangian submanifolds with constant  curvature $c$ of $\tilde M^n(4c)$.

 The new transformation also depends on an operator $P$ of a vector space $V$, so we call it the $P$-transformation.
 Proving the existence of $P$-transforms of a given submanifold satisfying certain initial conditions now requires investigating under which conditions a certain  Lyapunov-type matrix equation admits invertible solutions. The $P$-transformation  provides a process to construct a family of $n$-dimensional flat Lagrangian submanifolds of  $\mathbb{C}^n$ and $n$-dimensional horizontal  submanifolds with constant  curvature $c$ of  $\Sf_\epsilon^{2n+1}(c)$  starting with a given one and a vector-valued solution of a linear system of partial differential equations.

 We prove a decomposition result for the $P$-transformation and a corresponding Bianchi $P$-cube theorem, which shows how to generate a $k$-dimensional cube of 
$n$-dimensional flat Lagrangian submanifolds of  $\mathbb{C}^n$ (respectively,  $n$-dimensional horizontal submanifolds with constant sectional curvature $c$ of  $\Sf_\epsilon^{2n+1}(c)$) starting with $k$ such submanifolds. As before, the remaining  $2^k -(k+1)$ vertices of the cube are given by explicit algebraic formulae.   We also obtain a further reduction of the $P$-transformation, as well as  corresponding decomposition and Bianchi-cube theorems,  that preserves the class  of $n$-dimensional flat Lagrangian submanifolds of  $\R^{2n}$ that are contained in $\Sf^{2n-1}$. These are the lifts by the Hopf projection $\pi\colon \Sf^{2n-1}\to 
\mathbb C\mathbb P^{n-1}$ of $(n-1)$-dimensional flat Lagrangian submanifolds of   $\mathbb C\mathbb P^{n-1}$. 

We illustrate the procedures in the paper  by applying the $P$-transformation to the vacuum solution of the system of partial differential equations associated to $n$-dimensional flat Lagrangian submanifolds of  $\mathbb{C}^n$. First, starting with scalar $P_i$-transforms, $1\leq i\leq k$, of the associated (degenerate) submanifold (which are themselves immersed submanifolds), with $P_i\neq \pm P_j$ for $1\leq i\neq j\leq k$, we write down the formulae for the remaining  $2^k -(k+1)$ vertices of the corresponding Bianhi $P$-cube.  We also produce an example of an  $n$-dimensional flat Lagrangian submanifold of  $\mathbb{C}^n$  by applying the $P$-transformation to the vacuum solution with a non-diagonalizable $P$,  which therefore  can not be produced by an iteration of scalar $P$-transformations.

We point out that the vectorial Ribaucour transformation of flat Lagrangian submanifolds of  $\mathbb C^n$ and $\mathbb C\mathbb P^{n-1}$ was also  studied in \cite{tw} under a different approach based on the dressing action of a rational loop group on these immersions.

\section{Preliminaries}

In this section we recall some known results on systems of partial differential equations associated to isometric immersions of space forms into space forms, as well as on the Ribaucour transformation and its vectorial version. Since our study of Lagrangian submanifolds of constant curvature $c$ of $\mathbb C\mathbb H^n(4c)$ requires working on Euclidean space with a flat metric of index two, we consider pseudo-Riemannian ambient space forms.

\subsection{PDE's associated to submanifolds of constant curvature}

Given an isometric immersion $f\colon M^n\rightarrow \mathbb Q_s^{n+p}(c)$ of a Riemannian manifold into a pseudo-Riemannian space form of constant sectional curvature $c$ and index $s$, we denote by $N_1^f(x)$ the \emph{first normal space} of $f$ at $x\in M^n$, which is the subspace of the normal space $N_fM(x)$ of $f$ at $x$ spanned by the image of the second second fundamental form $\alpha_f$ at $x$. We say that $N_1^f(x)$ is \emph{nondegenerate} if $N_1^f(x)\cap N_1^f(x)^\bot=\{0\}$. We also denote by $\nu_f(x)$ the \emph{index of relative nullity} of $f$ at $x$, defined as the dimension of the kernel of $\alpha_f$ at $x$.\vspace{1ex}

The following results are well known (cf. \cite{t} and \cite{dt4}).

\begin{proposition}\label{thm:hiebetaij}
Let $f\colon M^n(c)\rightarrow \mathbb Q^{n+p}_s(c)$ be an isometric immersion with flat normal bundle and $\nu_f\equiv 0$. If $s\geq 1$, suppose further that $N_1^f(x)$ is nondegenerate everywhere. Then $p\geq n$ and there exist  locally  principal coordinates $(u_1,..., u_n)$, an orthonormal normal frame $\xi_1,...,\xi_{p}$ and smooth functions $v_1,...,v_n$ and $h_{ir}$, $1\leq i\leq n$, $1 \leq r\neq i\leq p$, with $v_1, \ldots, v_n$ positive, such that
\begin{equation}\label{alphasist}
ds^2=\sum_j v_j^2du_j^2, \ \ \ \alpha\left(\partial_i,\partial_j\right)=v_i\delta_{ij}\xi_i, \ \ \ \nabla_{\partial_i}X_j= h_{ji}X_i\ \ \ \mbox{and}\ \ \ \nabla^\bot_{\partial_i}\xi_r=h_{ir}\xi_i
\end{equation} 
where $\partial_i=\frac{\partial}{\partial u_i}$,  
$X_i = (1/v_i)\partial_i$ and 
$h_{ij} = (1/v_i
)\partial_i v_j$ for $1\leq i\neq j\leq n$. Moreover, the pair $(v, h)$,
where $v = (v_1,..., v_n)$ and $h = (h_{ir})$, satisfies the system of PDE's
\begin{equation}\label{eq:i}
\left\{
\begin{array}{lc}
i) \ \partial_j v_i = h_{ji}v_j, & ii) \ \partial_i h_{ij} +\partial_j h_{ji}+\sum_{\ell=1}^n h_{\ell i}h_{\ell j}+cv_iv_j=0,\vspace{0.1cm}\\
iii) \ \partial_j h_{ir} = h_{ij}h_{jr}, & iv) \ \epsilon_j\partial_j h_{ij} +\epsilon_i\partial_i h_{ji}+\sum_{r=1}^p \epsilon_rh_{ir}h_{jr}=0,\vspace{0.1cm}\\
\end{array}
\right.
\end{equation}
where always  $i\neq j, \ \{\ell,r\}\cap \{i,j\}=\emptyset$ and $\epsilon_r=\left<\xi_r,\xi_r\right>.$

Conversely, if $(v, h)$ is a solution of (\ref{eq:i}) on an open simply connected subset $U \subset \R^n$ such that $v_i > 0$ everywhere then there exists an immersion  $f: U\rightarrow  \Q^{n+p}_s(c)$ with flat normal bundle, $\nu_f \equiv 0$, nondegenerate first normal bundle of rank $n$ and induced metric $ds^2 = \sum_iv_i^2 du_i^2$ of constant sectional curvature $c$. 
\end{proposition}

 In the next proposition,  $\mathbb O_s^t(n\times p)$ stands for the subspace of $M_{n\times p}(\mathbb R)$ of all matrices $V$ that satisfy $V^tJV=\tilde J$, where $J_{ij}=\epsilon_i\delta_{ij}$ and $\tilde J_{ij}=\tilde \epsilon_i\delta_{ij}$, the $\epsilon_i$ being $-1$ for $s$ of the indices $1,...,p$ and 1 for the others, and the $\tilde \epsilon_i$ being $-1$ for $t$ of the indices $1,...,n$ and 1 for the others. 

\begin{proposition}\label{teo:curvdif}
Let $f\colon M^n(c)\rightarrow \mathbb Q^{n+p}_s(\widetilde c)$,  $c\neq \tilde c$, be an isometric immersion with flat normal bundle. Let $X_1,...,X_n$ be an orthonormal  frame that diagonalizes the second fundamental form of $f$ and let $\eta_i=\alpha_f(X_i, X_i)$, $1\leq i\leq n$,  be the principal normal vector fields of $f$. Assume that
$$
\theta_i=\<\eta_i, \eta_i\>+\tilde c-c\neq 0
$$
for $1\leq i\leq n$. Then $p\geq n-1$ and there exist locally  principal coordinates $u_1,...,u_n$ on $M^n(c)$ for~$f$. Moreover, if 
$\xi_1,...,\xi_{p}$ is a parallel orthonormal  frame of $N_fM$ and  $V_{ir}$, $1\leq  i\leq n$, $1\leq  r\leq p$, are defined by
$$
A_{\xi_r}X_i=v_i^{-1}V_{ir}X_i,
$$
where $v_i=|\partial_i|$ for $1\leq i\leq n$, then  the triple $(v,h,V)$, with $v=(v_1, \ldots, v_n)$, $h=(h_{ij})$ for $h_{ij} = (1/v_i
)\partial_i v_j$,  $1\leq i\neq j\leq n$, and $V=(V_{ir})$, 
satisfies the system of PDE's   
\begin{equation}\label{eq:ii}
\left\{
\begin{array}{ll}
i) \ \partial_i v_j=h_{ij}v_i & ii) \ \partial_\ell V_{ir} = h_{\ell i}V_{\ell r},\vspace{0.1cm}\\
iii) \ \partial_i h_{ij} +\partial_j h_{ji}+\sum_{\ell=1}^n h_{\ell i}h_{\ell j}+cv_iv_j=0, & iv) \ \partial_j h_{i\ell} = h_{ij}h_{j\ell}, \vspace{0.1cm}\\
\end{array}
\right.
\end{equation}
where $1\leq i\neq j\neq \ell\neq i\leq n$. Furthermore, the matrix $\hat V\in M_{n\times (p+1)}(\mathbb R)$  defined by 
\be\label{eq:hatV}
\hat V_{ir}=V_{ir},  \ 1\leq r\leq p, \ \mbox{and} \ \hat V_{i(p+1)}=\sqrt{|c-\tilde c|}v_i
\ee
belongs to $\mathbb O^t_{s+\epsilon_0}(n\times (p+1))$, where  $\epsilon_0=0$ or $1$ according to whether $\tilde c>c$ or $\tilde c <c,$ respectively, and $t$ is the number of indices for which $\theta_i<0.$

Conversely,  if $(v,h,V)$ is a solution of (\ref{eq:ii}) on an open simply connected subset $U\subset \mathbb{R}^n$ such that $v_i\neq 0$ everywhere and such that $\hat V\in M_{n\times (p+1)}(\mathbb R)$,  defined by  (\ref{eq:hatV}) for $\tilde c\neq c$, belongs to $\mathbb O^t_{s+\epsilon_0}(n\times (p+1))$, then there exists an  immersion $f:U\rightarrow \mathbb Q_s^{n+ p}(\tilde c)$, with $\theta_i<0$ for $t$ of the indices, that has $(v,h,V)$ as associated triple and whose induced metric $ds^2=\sum_iv_i^2du_i^2$ has constant sectional curvature $c$.
\end{proposition}

\subsection{The vectorial Ribaucour transformation}

If $M^n$ is an $n$-dimensional Riemannian manifold
and  $\xi$ is a pseudo-Riemannian vector bundle over $M^n$ endowed
with a compatible connection $\nabla^\xi$, we
denote by $\Gamma(\xi)$ the space of smooth sections of $\xi$.
If
$\zeta=\xi^*\otimes \eta=\mbox{Hom}(\xi,\eta)$ is the
tensor product of the vector bundles $\xi^*$ and $\eta$, where $\xi^*$
stands for the dual vector bundle of $\xi$   and $\eta$ is a
pseudo-Riemannian vector bundle over $M^n$, then
$\nabla Z\in\Gamma(T^*M\otimes \zeta)$  is given by
$$
(\nabla^\zeta_X Z)(v)
=\nabla^\eta_X Z(v)-Z(\nabla^\xi_X v)
$$
for all $X\in \mathfrak{X}(M)=\Gamma(TM)$,   $v\in\Gamma(\xi)$. For $Z\in \Gamma(\xi^*\otimes \eta)$, we define
$Z^t\in\Gamma(\eta^*\otimes \xi)$ by
$$
\<Z^t(u),v\>=\<u,Z(v)\>, \;\;\;u\in
\Gamma(\eta),  v\in \Gamma(\xi).
$$
The exterior derivative $d\omega\in\Gamma(\Lambda^2T^*M\otimes \xi)$
of $\o$ is related to its covariant derivative  by
$$
d\omega(X,Y)=\nabla \omega(X,Y) - \nabla\omega(Y,X)=(\nabla_X \omega)Y - (\nabla_Y\omega)X.
$$
The one-form $\omega$ is {\em closed\/} if
$d\omega=0$. If $Z\in \Gamma(\xi)$, then
$\nabla Z=dZ\in\Gamma(T^*M\otimes \xi)$
is the one-form given by $\nabla
Z(X)=\nabla^\xi_X Z$. In case $\xi=M\times V$ is a trivial vector
bundle over $M^n$,  with $V$ a Euclidean vector space, that is, a
vector space endowed with an inner product, then $\Gamma(T^*M\otimes \xi)$
is identified with the space of smooth one-forms with values in $V$. \vspace{1ex}
%We use the same notation for the vector space $V$ and the trivial
%vector bundle $\xi=M\times V$.

  The vectorial Ribaucour transformation for Euclidean
  submanifolds was introduced in \cite{dft}, and can be extended to submanifolds of $\R_s^{n+p}$ as follows.\vspace{1ex}
  
  Let $f\colon\,M^n\to \R_s^{n+p}$ be an isometric
immersion of a simply connected Riemannian manifold. Let $V$ be a vector space and let $\va\colon M^n\to V$ be a smooth map,
which we  regard as a section of the trivial vector bundle $M^n\times V$, also denoted by $V$. Fix an inner product $\<\cdot , \cdot\>$ on $V$, thus making
$M\times V$ into a Riemannian vector bundle, which we endow with the trivial compatible connection $\nabla=d$.
% Denote $\psi=\va^t\in \Gamma(V^*)$, also regarding $\va$ as a section of $\Gamma(\R^*\otimes V)$. Thus
%$$
%\psi(v)=\<\va, v\>
%$$
%for all $v\in \Gamma(V)$. 
%Denote also $\omega=d\va$ and $\nu=\omega^t\in \Gamma(V^*\otimes TM)$. Now 
Assume that 
$\beta\in
\Gamma(V^*\otimes N_fM)$ is such that 
\be\label{eq:alphao}\alpha(X,\omega^t(v))+(\nabla_X\beta)v=0\ee
 for all $v\in \Gamma(V)$, where $\omega=d\va\in \Gamma(TM^*\otimes V)$,
 %, with $\o=d\va$,
and 
\begin{equation}\label{eq:phicom}
[\Phi_u,\Phi_v]=0 
\end{equation}
for all $u,v\in \Gamma(V)$, where
 $\Phi=\Phi(\omega, \beta)\in \Gamma(T^*M\otimes V^*\otimes TM)$ is  given by 
  \be\label{eq:Phi}\Phi_vX:=\Phi(X)(v)
=(\nabla_X\omega^t)v-A_{\beta(v)}X .\ee 
It follows from (\ref{eq:alphao}) that $\Fes\in \Gamma(V^*\otimes f^*T\R_s^{n+p})$ given by  $\Fes=f_*\omega^t+ \beta$ satisfies
\be\label{eq:combe0}
d\Fes(X)(v)=f_*\Phi_vX
\ee
for all $X\in \mathfrak{X}(M)$ and  $v\in \Gamma(V)$. Now let $\Omega\in \Gamma(V^*\otimes V)$ be a solution of
the completely integrable first order system
\be\label{eq:o1}d\Omega=\Fes^td\Fes\ee such that
\be\label{eq:o2}\Omega+\Omega^t=\Fes^t\Fes.\ee 
We point out that the integrability conditions of (\ref{eq:o1}) 
are precisely equations (\ref{eq:phicom}) (see \cite{dft}), and that (\ref{eq:o1}) implies (\ref{eq:o2}) up to a parallel section of $V^*\otimes V$. 

\begin{definition}\po\label{ribvet} {\em If $\Omega$ is invertible everywhere, then the isometric immersion $\tilde{f}\colon\tilde{M}^n\to \R_s^{n+p}$ given
by \be\label{eq:rb2} \tilde{f}=f-\Fes\Omega^{-1}\va, \ee where $\tilde{M}^n$ is the  subset  of regular points of $\tilde f$ endowed with the induced metric, is called  the {\em vectorial Ribaucour
transform\/} of $f$ determined by  $(\va,\beta,\Omega)$ (or, more precisely, by $(\va,\beta,\Omega, \<\cdot , \cdot\>)$; see Remark \ref{innerproduct} below).
We write $\tilde{f}={\cal R}_{\va,\beta,\Omega}(f)$.}
\end{definition}
  
 \begin{remark}\po\label{innerproduct} {\em Let $\<\cdot , \cdot\>^\sim$ be another inner product on $V$, related to $\<\cdot , \cdot\>$ by
 %\be\label{eq:innerp}
 $
 \<\cdot , \cdot\>^\sim = \<\cdot, B\cdot\>
 $
 for some invertible  $B\in V^*\otimes V$, which we regard as a parallel section of the trivial vector bundle $V^*\otimes V$,
 %we endow the trivial vector bundle $M\times V$ with the compatible connection given by
 \be\label{eq:cov}  \nabla_X Bv=B\nabla_X v\ee
 %is a compatible connection on  $M\times V$ with respect to $\<\cdot, \cdot\>^\sim$.
 %for all $X\in \mathfrak{X}(M)$ and $v\in \Gamma(V)$.
 for all  $X\in \mathfrak{X}(M)$ and $v\in \Gamma(V)$.  The transpose $\hat \nu$ of $\omega=d\va$ with respect to  $\<\cdot , \cdot\>^\sim$ is related to its transpose $\nu=\omega^t$ with respect to  $\<\cdot , \cdot\>$ by 
 $\hat \nu=\nu B$
 and, in view of  (\ref{eq:cov}), also the covariant derivatives $\nabla \hat \nu$ and $\nabla \nu$ are related by
 $$(\nabla_X \hat \nu)v=(\nabla_X\nu)Bv$$
 for all $X\in \mathfrak{X}(M)$ and $v\in \Gamma(V)$.  Thus, if we define $\hat \beta\in \Gamma(V^*\otimes N_fM)$ by $\hat \beta= \beta B$, then
 $$\hat \Fes=f_*\hat \nu+ \hat \beta=\Fes B$$ and 
 $$\hat\Phi_vX=(\hat\nabla_X\hat\nu)v-A_{\hat\beta(v)}X=\Phi_{Bv} X.$$
 Hence (\ref{eq:alphao}) is satisfied by $(\hat \nu, \hat \beta)$, and 
 $$[\hat\Phi_u, \hat\Phi_v]=[\Phi_{Bu}, \Phi_{Bv}]=0$$ for all $u, v\in \Gamma(V)$. Note also that the transposes $\hat \beta^t$ and $\beta^t$ with respect to $\<\cdot, \cdot\>^\sim$ and $\<\cdot, \cdot\>$, respectively, coincide, for
 $$\<\hat \beta^t \xi, v\>^\sim=\<\xi, \hat \beta v\>=\<\xi,  \beta Bv\>=\<\beta^t \xi,  Bv\>=\<\beta^t \xi,  v\>^\sim.$$
 Similarly, $\hat \Fes^t=\Fes^t$,  and hence   $\hat \Omega\in \Gamma(V^*\otimes V)$, given by $\hat \Omega=\Omega B$, satisfies 
 $\hat \Omega^t=\Omega^tB$, 
  $$d\hat\Omega=\hat\Fes^td\hat\Fes\;\;\mbox{and}\;\;\hat\Omega+\hat\Omega^t=\hat\Fes^t\hat\Fes.$$
 Moreover,
 $$f-\hat\Fes\hat\Omega^{-1}\va=f-(\Fes B)(B^{-1}\Omega^{-1})\va=f-\Fes\Omega^{-1}\va$$
 and thus $(\va, \beta, \Omega, \<\cdot, \cdot\>)$ and $(\hat \va=\va, \hat \beta, \hat \Omega, \<\cdot, \cdot\>^\sim)$ give rise to the same vectorial Ribaucour transform. 
%  Therefore, in the definition of the vectorial Ribaucour transformation we identify quadruples $(\va, \beta, \Omega, \<\cdot, \cdot\>)$ and $(\hat \va=\va, \hat \beta, \hat \Omega, \<\cdot, \cdot\>^\sim)$ that are related as above.
}
 \end{remark}

If $\dim V=1$,  after identifying $V^*\otimes N_fM$ with
 $N_fM$ then $\va$ and $\beta$ become elements of $C^\infty(M)$ and
$\Gamma(N_fM)$, respectively, and equation
(\ref{eq:alphao}) reduces to
\be\label{eq:reduced}\alpha(X,\nabla\va)+\nabla^\perp_X\beta=0.\ee
Moreover, $\Omega$ reduces to the function  $\Omega=(1/2)\<\Fes,\Fes\>$, where $\Fes=f_*\nabla\va+\beta$ is the \emph{Combescure transform} of $f$ determined by $(\varphi, \beta)$, 
and  (\ref{eq:rb2}) becomes the parameterization of a {\em
scalar\/}  Ribaucour transform of $f$ obtained in Theorem $17$ of
\cite{dt5}. 
%Geometrically, $f$ and $\tilde f$  envelope  a common $n$-sphere congruence in $\R^N$. 
In this case, since $\Omega$ is determined by $\va$
and $\beta$ we write $\tilde{f}={\cal R}_{\va,\beta}(f)$  instead
of $\tilde{f}={\cal R}_{\va,\beta,\Omega}(f)$. The tensor $\Phi=\Phi_{\varphi, \beta}$ given by (\ref{eq:Phi}) reduces in this case to
$$\Phi= \hess \varphi-A_\beta$$
and is a Codazzi tensor on $M^n$ that is related to ${\cal
F}$ by ${\cal F}_*X=f_*\Phi(X)$ for all $X\in \mathfrak{X}(M)$. It is called the \emph{Codazzi tensor associated to the Ribaucour transform $\tilde f$}.

Given two pairs $(\varphi,\beta)$ and $(\varphi',\beta')$ satisfying (\ref{eq:reduced}) such that  $[\Phi_{\varphi',\beta'},\Phi_{\varphi,\beta}]=0$, taking linear combinations $(\tilde \varphi, \tilde \beta)=c(\varphi,\beta)+c'(\varphi',\beta')$, $c,c'\in \R$, yields a \emph{one}-parameter family of scalar Ribaucour transforms 
${\cal R}_{\tilde\va,\tilde\beta}(f)$ of $f$, because pairs $(\varphi,\beta)$ and $(\varphi',\beta')$ related by $(\varphi,\beta)=\lambda(\varphi',\beta')$ for some $\lambda\neq 0$ give rise to the same Ribaucour transform. It is called the 
\textit{associated family} determined by ${\cal R}_{\va,\beta}(f)$ and ${\cal R}_{\va',\beta'}(f)$.\vspace{1ex}

 The vectorial Ribaucour transformation can be extended to isometric immersions $f\colon M^n\rightarrow \mathbb Q_s^{n+p}(\tilde c)$
 of a Riemannian manifold $M^n$ into a pseudo-Riemannian manifold of constant sectional curvature $\tilde c$ as follows. 
 Here, and in the sequel, $i$ stands for the umbilical inclusion
 $i:\mathbb Q^{n+p}_s(\tilde c)\rightarrow \mathbb R^{n+p+1}_{s+\epsilon_0}$, where  $\epsilon_0=0$ or $\epsilon_0=1$ depending on whether $\tilde c>0$ or $\tilde c<0$, respectively.

 \begin{definition}\label{ribvetgeral} {\em An isometric immersion 
  $\tilde f:\tilde M^n\rightarrow \mathbb Q_s^{n+p}(\tilde c)$ is said to be the vectorial Ribaucour transform of $ f: M^n\rightarrow \mathbb Q_s^{n+p}(\tilde c)$ with data $(\varphi,\beta,\Omega)$ if $\tilde F=i\circ \tilde f: \tilde M^n\rightarrow \mathbb R^{n+p+1}_{s+\epsilon_0}$ is a vectorial Ribaucour transform of  $F=i\circ f: M^n\rightarrow \mathbb R^{n+p+1}_{s+\epsilon_0}$ with data $(\varphi,\tilde \beta,\Omega)$, where 
 %\begin{equation}\label{eq:betatil}
 $$
\tilde \beta=i_*\beta +\tilde cF\va^t.
$$
%\end{equation} 
}
%\emph{Here, and in the sequel, $i$ stands for the umbilical inclusion
% $i:\mathbb Q^{n+p}_s(\tilde c)\rightarrow \mathbb R^{n+p+1}_{s+\epsilon_0}$, where  $\epsilon_0=0$ or $\epsilon_0=1$ depending on whether $\tilde c>0$ or $\tilde c<0$, respectively}. \vspace{1ex}
\end{definition}

In this case, the corresponding vectorial Combescure transform of $F$ is
$${\cal G}={F}_*\omega^t+i_*\beta+\tilde cF\va^t$$
and the  tensor $\Phi=\Phi(\omega, \beta)$ is now   given by 
  \be\label{eq:Phi2}\Phi_vX:=\Phi(X)(v)
=(\nabla_X\omega^t)v-A_{\beta(v)}X +\tilde c\va^t(v)X.\ee 
Equations (\ref{eq:o1}) and (\ref{eq:o2}) become
\be\label{eq:o1b}d\Omega={\cal G}^td{\cal G}\ee 
and
\be\label{eq:o2b}\Omega+\Omega^t={\cal G}^t{\cal G},\ee
respectively, and (\ref{eq:o1b}) is equivalent to
\be\label{eq:o1c}d\Omega=\omega \Phi.\ee

The following basic properties of the vectorial 
transformation were derived in \cite{dft} when  $\tilde c=0$.
The extension to the case $\tilde c\neq 0$ is straightforward.

\begin{proposition}
\po\label{prop:basic} $(i)$  The bundle map
$$\P=I-{\cal G}\Omega^{-1}{\cal G}^t\in\Gamma((F^*T\R_{s+\epsilon_0}^{n+p+1})^*\otimes \tilde{F}^*T\R_{s+\epsilon_0}^{n+p+1})$$
 is a vector bundle
isometry satisfying 
$$\P F=\tilde F=F-{\cal G}\Omega^{-1}\va\,\,\,\,\mbox{and}\,\,\,\,\tilde{f}_*=\P f_*D$$
 where
$$D=I-\Phi_{\Omega^{-1}\va}\in \Gamma(T^*M\otimes TM).$$ In
particular, the
 metrics 
 %$\<\;,\;\>$  and $\<\;,\;\>^\sim$
  induced by $f$ and $\tilde{f}$ are related by
$\<\;\;,\;\>^\sim=D^*\<\;\;,\;\>.$\vspace{1ex}\\
$(ii)$ The normal
connections and second fundamental forms of $f$ and $\tilde{f}$
are related by
\be\label{eq:ncon}\tilde{\nabla}_X^\perp{\P}\xi=\P\nabla_X^\perp\xi\ee
and \be\label{eq:sffs}
\tilde{A}_{\P\xi}=D^{-1}(A_\xi+\Phi_{\Omega^{-1}\beta^t\xi}), \ee
or equivalently,
\be\label{eq:sffs2}\tilde{\alpha}(X,Y)=\P(\a(X,DY)+\beta(\Omega^{-1})^t\Phi(X)^tDY).
\ee
$(iii)$  The Levi-Civita connections of the
 metrics $\<\;,\;\>$  and
$\<\;,\;\>^\sim$ are related by \be\label{eq:lcivita}
D\tilde{\nabla}_XY=\nabla_X
DY-(\Phi(X)\Omega^{-1}\o-(\Phi(X)\Omega^{-1}\o)^t)(DY). \ee
\end{proposition}

\subsection{The decomposition theorem}

The following is a straightforward extension for nonflat ambient space forms of 
a decomposition property of  the vectorial Ribaucour transformation  proved
in \cite{dft} for Euclidean submanifolds.

\begin{theorem}\label{c:viiinf}
Let $\ral_{\varphi,\beta,\Omega}(f)\colon\tilde M^n\rightarrow \mathbb Q_s^{n+p}(\tilde c)$ be a vectorial Ribaucour transform of an isometric immersion $f\colon M^n\rightarrow \mathbb Q_s^{n+p}(\tilde c)$ determined by $(\varphi,\beta,\Omega)$ as in Definition \ref{ribvetgeral}. For a direct sum decomposition $V=V_1\oplus V_2$, which we assume to be orthogonal after a suitable change of the inner product on $V$ (see Remark \ref{innerproduct}), let $\pi_j\colon V\to V_j$, $1\leq j\leq 2$, denote the orthogonal projection and define 
\be\label{vabo2nf}\va_j=\pi_{V_j}\circ
\va,\,\,\,\,\beta_j=\beta|_{V_j}\,\,\,\mbox{and}\,\,\,\Omega_{ij}=\pi_{V_i}\circ
\Omega|_{V_j}\in \Gamma(V_j^*\otimes V_i),\,\,\,1\leq i,j\leq
2.\ee 
Then the triple $(\varphi_j,\beta_j,\Omega_{jj})$  satisfies the conditions of Definition \ref{ribvetgeral} for $1\leq j\leq 2$, except possibly for the invertibility of  $\Omega_{jj}$. Suppose that  the latter condition is also satisfied    and define $f_j=\ral_{\varphi_j,\beta_j,\Omega_{jj}}(f)$,
$F_j=i\circ f_j$,  
${\cal G}_j={F_j}_*\omega_j^t+i_*\beta_j+\tilde cF_j\varphi_j^t$, where $\omega_j=d\varphi_j$, and 
\be\label{eq:bar2nf}
\bar{\va}_i=\va_i-\Omega_{ij}\Omega_{jj}^{-1}\va_j,\;\;\;
\bar{\beta}_i=
\P_j(\beta_i-\beta_j(\Omega^{-1}_{jj})^t\Omega_{ij}^t)\;\;\;\mbox{and}\;\;\;
\bar{\Omega}_{ii}=
\Omega_{ii}-\Omega_{ij}\Omega_{jj}^{-1}\Omega_{ji}, \ee $1\leq
i\neq j\leq 2$, where ${\cal P}_j=I-{\cal G}_j\Omega_{jj}^{-1}{\cal G}_j^t$. Then the triple  $(\bar \varphi_i,\bar \beta_i,\bar \Omega_{ii})$ satisfies the conditions of Definition \ref{ribvetgeral} with respect to $f_j$, $1\leq j\neq i\leq 2$, and  
\begin{equation}\label{eq:permnf2}
\ral_{\varphi,\beta,\Omega}(f)=\ral_{\bar \varphi_i, {\bar \beta_i},\bar \Omega_{ii}}(\ral_{\varphi_j,\beta_j,\Omega_{jj}}(f)).
\end{equation}
\end{theorem}

 Theorem \ref{c:viiinf} implies that   any
vectorial Ribaucour transformation whose associated data
$(\va,\beta,\Omega)$ are defined on a vector space $V$ can be
regarded as the iteration of $k=\mbox{dim} V$ scalar Ribaucour transformations.
We  discuss below a  consequence of this fact.
\vspace{1ex}

The isometric immersions $f_i\colon M_i^n\rightarrow \mathbb Q_s^{n+p}(\tilde c), \ 1\leq i\leq 4$,  form a \emph{Bianchi quadrilateral\/} if for each of them  the preceding and subsequent ones (thought of as points on an oriented circle) are Ribaucour transforms of it and the associated Codazzi tensors  commute\vspace{1ex}.

A \textit{Bianchi cube} is a $(k+1)$-tuple $({\cal{C}}_0, . . . ,{\cal{C}}_k)$, where each ${\cal{C}}_r$, $1\leq r\leq k$, is a family of isometric immersions $f_{\alpha_{r}}\colon M_{\alpha_r}^n\rightarrow \mathbb Q_s^{n+p}(\tilde c)$ indexed in the set of multi-indices
$$
\Lambda_r=\{\alpha_r=\{i_1,...,i_r\}\subset \{1,...,k\}\,:\, \alpha_r\,\mbox{with $r$ distinct elements}\}
$$
satisfying the following conditions for all $1\leq s\leq k-1$:
\begin{itemize}
	\item[(i)] Each $f_{\alpha_{s+1}}\in {\cal C}_{s+1}$ with  $\alpha_{s+1}=\alpha_{s}\cup \{i_j\}$  is a Ribaucour transform of $f_{\alpha_{s}}\in {\cal C}_s$.
		\item[(ii)]   $\{f_{\alpha_{s-1}},f_{\alpha_{s-1}\cup \{i_l\}},f_{\alpha_{s-1}\cup \{i_j\}},f_{\alpha_{s+1}}\}$ is a Bianchi quadrilateral if  $\alpha_{s+1}=\alpha_{s-1}\cup \{i_l,i_j\}$.
\end{itemize}

 Theorem 2 in \cite{dft} can be easily extended to nonflat ambient space forms as follows.

%\begin{theorem}\label{teo:cuboor}
%Let $f\colon M^n\rightarrow \mathbb Q_s^{n+p}(\tilde c)$ be an isometric immersion and let $f_1,...,f_k$ be Ribaucour transforms of $f$.  Then, for a generic choice of  isometric immersions $f_{ij}\colon M_{ij}^n\rightarrow \mathbb Q_s^{n+p}(\tilde c)$,  $1\leq i\neq j\leq k$, such that $\{f,f_i,f_j,f_{ij}\}$ is a Bianchi quadrilateral, there exists a  Bianchi k-cube $({\cal C}_0,...,{\cal C}_k)$ such that ${\cal C}_0=\{f\}, \ {\cal C}_1=\{f_1,...,f_k\}$ and ${\cal C}_2=\{f_{ij}\}_{1\leq i\neq j\leq k}$,  which is unique if none of $f_1,...,f_k$ belongs  to the associated family determined by any two of the others. 
%\end{theorem}

\begin{theorem}\label{teo:cuboor}
Let $f\colon M^n\rightarrow \mathbb Q_s^{n+p}(\tilde c)$ be an isometric immersion and let $f_1,...,f_k$ be Ribaucour transforms of $f$ none of which belongs  to the associated family determined by any two of the others. Then, for a generic choice of  isometric immersions $f_{ij}\colon M_{ij}^n\rightarrow \mathbb Q_s^{n+p}(\tilde c)$,  $1\leq i\neq j\leq k$, such that $\{f,f_i,f_j,f_{ij}\}$ is a Bianchi quadrilateral, there exists a unique Bianchi k-cube $({\cal C}_0,...,{\cal C}_k)$ such that ${\cal C}_0=\{f\}, \ {\cal C}_1=\{f_1,...,f_k\}$ and ${\cal C}_2=\{f_{ij}\}_{1\leq i\neq j\leq k}$.  
\end{theorem}

\section{The $L$-transformation}

  In this section we obtain a reduction of  the vectorial Ribaucour
  transformation that preserves the class of  submanifolds with constant  curvature.
  First we relate the curvature tensors of submanifolds that are associated by a vectorial Ribaucour transformation.\vspace{1ex}

  Let $f\colon\,M^n\to \Q^{n+p}_s(\tilde c)$ be an isometric immersion and let  $\tilde{f}=\ral_{\varphi,\beta,\Omega}(f)\colon\,\tilde
M^n\to \Q^{n+p}_s(\tilde c)$ be a  vectorial Ribaucour transform of $f$.
 Set $S=\Omega^{-1}\beta^t\in
\Gamma(N_fM^*\otimes V)$ and define $U\in \Gamma(T^*M\otimes T^*M\otimes T^*M\otimes TM)$ by
$$U(X,Y)=\Phi(X) S(A(Y)+\frac{1}{2}\Phi(Y) S)^t.$$
Given $T\in \Gamma(T^*M\otimes T^*M\otimes T^*M\otimes TM)$, define
$\hat{T}\in \Gamma({\cal A}_2(T^*M)\otimes {\cal A}(TM))$ by
$$\hat{T}(X,Y)=(T(X,Y)-T(X,Y)^t)- (T(Y,X)-T(Y,X)^t).$$

\begin{proposition}\po\label{prop:curvs} The curvature tensors of $M^n$
and $\tilde{M}^n$
 are related by
%\be\label{eq:curvs}
$$
D\tilde{R}(X,Y)=(R(X,Y)-\tilde c(X\wedge Y)+\hat{U}(X,Y)+\tilde c (DX\wedge DY))
D.
$$
%\ee
\end{proposition}
\proof By (\ref{eq:sffs}) and   (\ref{eq:sffs2}) we have
\be\label{eq:datil}\begin{array}{l} D\tilde{A}_{\tilde{\alpha}(Y,Z)}X=A(X)(\alpha(Y,DZ)+S^t\Phi(Y)^tDZ)\vspace{1ex}\\\hspace*{13ex}
+\Phi(X)S(\alpha(Y,DZ)+S^t\Phi(Y)^tDZ).
\end{array}\ee
From (\ref{eq:datil}) and the Gauss equations for $\tilde{f}$ and $f$  we obtain
$$\begin{array}{l}D\tilde{R}(X,Y)Z=
D\tilde{A}_{\tilde{\alpha}(Y,Z)}X-D\tilde{A}_{\tilde{\alpha}(X,Z)}Y+D(\bar R(X,Y)Z)^t\vspace{1ex}\\
\hspace*{14.8ex}=
R(X,Y)DZ-\tilde c (X\wedge Y)DZ+A(X)S^t\Phi(Y)^tDZ\vspace{1ex}\\
\hspace*{17ex}+\Phi(X)SA(Y)^tDZ+\Phi(X)SS^t\Phi(Y)^tDZ-A(Y)S^t\Phi(X)^tDZ\vspace{1ex}\\
\hspace*{17ex}-
\Phi(Y)SA(X)^tDZ-\Phi(Y)SS^t\Phi(X)^tDZ+\tilde c(DX\wedge DY)DZ\vspace{1ex}\\
\hspace*{15ex}=R(X,Y)DZ-\tilde c(X\wedge Y)DZ+\hat{U}(X,Y)DZ+\tilde c(DX\wedge DY)DZ. \vspace{2ex}\qed
\end{array}
$$
Now define $\rho\in  \Gamma(V^*\otimes V)$,  $\Psi\in \Gamma(T^*M\otimes V^*\otimes TM))$ and $Q\in \Gamma(T^*M\otimes T^*M\otimes T^*M\otimes TM)$ by
\be\label{eq:rho}\rho=\beta^t\beta-(c-\tilde c)\va\va^t,\ee  
\be\label{eq:Psi}\Psi(Y)=A(Y)\beta+(c-\tilde c)Y\va^t+\frac{1}{2}\Phi(Y)\Omega^{-1}\rho\ee
and
\be\label{eq:Q} 
Q(X,Y)=\Phi(X)\Omega^{-1} \Psi(Y)^t.
\ee

\begin{corollary}\po\label{cor:csc} If $M^n$ has constant sectional curvature
$c$, then the same holds for $\tilde{M}^n$ if and only if $\hat{Q}$ vanishes identically.
\end{corollary}
\proof First observe that $\tilde{M}^n$ has  constant sectional
curvature $c$ if and only if
$$D\tilde{R}(X,Y)Z=cD(\<Y,Z\>^\sim X -
\<X,Z\>^\sim
Y)=c(DX\wedge DY)DZ.
$$
By Proposition \ref{prop:curvs} and the fact that $M^n$ has
constant sectional curvature $c$ we have
\begin{eqnarray}\label{eq:1}
D \tilde{R}(X,Y)-c(DX\wedge DY) D&=&(R(X,Y)-\tilde c(X\wedge Y)+\hat{U}(X,Y)\nonumber\\
&&+\tilde c(DX\wedge DY)-c(DX\wedge DY)) D \\
&=&((c-\tilde c)(X\wedge Y)-(c-\tilde c)(DX\wedge DY)
\nonumber
+\hat{U}(X,Y)) D.
\end{eqnarray}
Using that $D=I-\Phi_{\Omega^{-1}\va}$  we obtain
\be\label{eq:2}\begin{array}{l}
X\wedge Y-DX\wedge DY=XY^t-YX^t-DX(DY)^t+DY(DX)^t\vspace{1ex}\nonumber\\
\hspace*{20.5ex}=X(\Phi(Y)\Omega^{-1}\va)^t+(\Phi(X)\Omega^{-1}\va)Y^t-
(\Phi(X)\Omega^{-1}\va)(\Phi(Y)\Omega^{-1}\va)^t\vspace{1ex}\nonumber\\
\hspace*{21.5ex}-Y(\Phi(X)\Omega^{-1}\va)^t-(\Phi(Y)\Omega^{-1}\va)X^t
+(\Phi(Y)\Omega^{-1}\va)(\Phi(X)\Omega^{-1}\va)^t\vspace{1ex}\nonumber\\
\hspace*{20.5ex}=\hat{H}(X,Y)
\end{array}
\ee where
$H(X,Y)=\Phi(X)\Omega^{-1}\va(Y-\frac{1}{2}\Phi(Y)\Omega^{-1}\va)^t.$
It follows from (\ref{eq:1}) and (\ref{eq:2}) that
$$D\tilde{R}(X,Y)-c(DX\wedge DY) D=(\hat{U}(X,Y)+(c-\tilde c)\hat{H}(X,Y)) D,$$
and the proof is completed by checking that
$U(X,Y)+(c-\tilde c)H(X,Y)=Q(X,Y)$.\qed 
\vspace{2ex}

Looking for  triples $(\varphi, \beta, \Omega)$  for which $\hat Q=0$ leads to the following definition.

\begin{definition}\po\label{df:ribconst} {\em Let  $\tilde f={\cal
R}_{\va,\beta,\Omega}(f)\colon\,\tilde{M}^n\to \Q^{n+p}_s(\tilde c)$ be a vectorial Ribaucour transform of $f\colon\,M^n(c)\to \Q^{n+p}_s(\tilde c)$  determined by
$(\va,\beta,\Omega)$, with $\va\in \Gamma(V)$, $\beta\in \Gamma(V^*\otimes N_fM)$ and $\Omega\in \Gamma(V^*\otimes V)$. If there exists  $L\in V^*\otimes V$ such that  
\begin{equation}\label{eq:cconst}
\Phi(Y)L+A(Y)\beta+(c-\widetilde c)Y\varphi^t=0
\end{equation} 
where    $\Phi=\Phi(\nu, \beta)$ is given  by (\ref{eq:Phi2}), and 
\be\label{eq:cond3}
\Omega L + L^t \Omega^t=\rho
\ee
with $\rho$  given by (\ref{eq:rho}), then $\tilde f$ is said to be an {\em $L$-vectorial Ribaucour transform} of $f$, or simply an {\em $L$-transform} of $f$.
We write $\tilde f={\cal R}_{\varphi,\beta,\Omega,L}(f)$. }
\end{definition} 
\begin{remarks}\po\label{re:ltransform}{\em $(i)$ If $\<\cdot , \cdot\>^\sim$ is another inner product on $V$, with 
 %\be\label{eq:innerp}
 $\<\cdot, \cdot\>^\sim = \<\cdot, B\cdot\>$ 
 for some invertible  $B\in V^*\otimes V$, then conditions (\ref{eq:cconst}) and (\ref{eq:cond3}) are satisfied by $(\va, \beta, \Omega, L,  \<\cdot, \cdot\>)$  if and only if they are  satisfied by $(\hat\va, \hat\beta, \hat\Omega, \hat L,  \<\cdot, \cdot\>^\sim)$, where 
  $\hat\va=\va$, $\hat\beta=\beta B$, $\hat\Omega=\Omega B$ and $\hat L=B^{-1}LB$. Thus, in the definition of the $L$-transformation one may replace $(\va, \beta, \Omega, L,  \<\cdot, \cdot\>)$ by any quintuple  $(\hat\va, \hat\beta, \hat\Omega, \hat L,  \<\cdot, \cdot\>^\sim)$  related to $(\va, \beta, \Omega, L,  \<\cdot, \cdot\>)$ in this way.
 Notice  that the transposes of $\hat L$ and $L$ with respect to $\<\cdot, \cdot\>^\sim$ and $\<\cdot, \cdot\>$ coincide, for
 $$\<\hat L^t v, w\>^\sim=\<v, \hat L w\>^\sim=\<v, B\hat Lw\>=\<v, LB w\>=\<L^t v, Bw\>=\<L^t v, w\>^\sim.$$ 
$(ii)$ In the scalar case, the $L$-transformation reduces to  the transformation for submanifolds of constant sectional curvature given by Theorem $13$ of \cite{dt4}, called the ${\cal R}_C$-transformation in that paper, which in turn reduces to the Ribaucour transformation for constant curvature  surfaces in $\R^3$ studied by  Bianchi \cite{bi}.
}
 \end{remarks}

\begin{theorem}\po\label{thm:csc} If $f\colon\,M^n(c)\to \Q^{n+p}_s(\tilde c)$ is an isometric immersion
 and  $\tilde f\colon\tilde{M}^n\to \Q^{n+p}_s(\tilde c)$ is an
$L$-transform of $f$, then $\tilde{M}^n$  also has constant curvature $c$.
\end{theorem}
\proof We get from (\ref{eq:cconst}) that 
$\Psi$  defined by (\ref{eq:Psi}) satisfies 
$\Psi(Y)=\Phi(Y)\Omega^{-1}(\rho/2-\Omega L)$, and hence $Q$ given by (\ref{eq:Q}) reduces to 
$Q(X,Y)=\Phi(X)\Omega^{-1}(\rho/2-L^t\Omega^t)(\Omega^{-1})^t\Phi(Y)^t.$
Using (\ref{eq:cond3}), we obtain  $\hat{Q}(X,Y)=0$,  and the statement
follows from Corollary \ref{cor:csc}.\vspace{1ex}\qed

To prove the existence of $L$-transforms we first write (\ref{eq:alphao}) and (\ref{eq:cconst}) in the local coordinates given by Propositions \ref{thm:hiebetaij} and \ref{teo:curvdif}.

\begin{lemma}\po\label{prop:fundconst}
Let $f\colon M^n(c)\rightarrow \mathbb Q^{n+p}_s(c)$
be an isometric immersion  satisfying 
 the assumptions of  Proposition \ref{thm:hiebetaij}, 
 and  let $u_1, \ldots, u_n$,    $\xi_1, \ldots, \xi_p$ and $(v,h)$ be, respectively,  principal coordinates on an open subset $U\subset M^n(c)$,  the orthonormal frame of $N_fU$ and the pair associated to $f$ given by that result. Given $\varphi\colon U\to V$ and 
 $\beta\in \Gamma(V^*\otimes N_fU)$, where $V$ is a Euclidean vector space, define $\gamma_1,...,\gamma_n, \beta^1,...,\beta^p\colon U\to V$ by  
 \be\label{eq:betagamma}\gamma_i=v_i^{-1}\omega(\partial_i)\,\,\,\mbox{and}\,\,\,\,\beta^r=\beta^t(\xi_r)\ee
  for $1\leq i\leq n$ and $1\leq r\leq p$, where $\omega=d\varphi$. Then $(\omega, \beta)$ satisfies (\ref{eq:alphao}) and 
  $\Phi(\omega, \beta)$ satisfies (\ref{eq:cconst}) if and only if $(\varphi,\gamma_1,...,\gamma_n,\beta^1,...,\beta^p)$ is a solution of   the linear system of PDE's
 $$
\ral=\left\{
\begin{array}{ll}
i) \ \partial_i\varphi=v_i\gamma_i, \,\,\,\, ii) \ \partial_i\gamma_j= h_{ji}\gamma_i, \,\,\,i\neq j\\
iii) \ \partial_i \gamma_i =- \sum_{j\neq i}h_{ji} \gamma_j-((L^t)^{-1}-I)\beta^i -cv_i\varphi& \vspace{0.1cm}\\
iv)\ \partial_i\beta^r=h_{ir}\beta^i,\,\,i\neq r,\,\,\,\,\,v)\ \epsilon_i\partial_i\beta^i=-\gamma_i-\sum_{r\neq i}\epsilon_rh_{ir}\beta^r.\\
\end{array}
\right.
$$ 
\end{lemma}
\proof Equation $(i)$ of $\ral$ is equivalent to $\omega=d\varphi$. Now we have
\be\label{eq:omegabeta}
\omega^t(v)=\sum_{i=1}^nv_i^{-1}\<\gamma_i, v\>\partial_i\,\,\,\,\mbox{and}\,\,\,\,\beta(v)=\sum_{r=1}^p\epsilon_r\<\beta^r,v\>\xi_r,
\ee
hence using  
(\ref{alphasist}) we obtain
$\alpha_f(\partial_i, \omega^t(v))=\<\gamma_i, v\>\xi_i$
and
$$(\nabla_{\partial_i}\beta)(v)=\sum_{r=1, r\neq i}^p\epsilon_r\<\partial_i\beta^r-h_{ir}\beta^i,v\>\xi_r+\<\epsilon_i\partial_i\beta^i+\sum_{r=1, r\neq i}^p \epsilon_r h_{ir}\beta^r, v\> \xi_i.$$
Thus (\ref{eq:alphao})  is equivalent to $(iv)$ and $(v)$. Moreover,  using the third formula in (\ref{alphasist}),  it follows that $\Phi=\Phi(\omega, \beta)$ satisfies
%\begin{equation}\label{eq:phipi}
$$
\Phi_{v} \left(\partial_i\right)= \sum_{j\neq i}\<\partial_i\gamma_j-h_{ji}\gamma_i,v\>X_j +\<\partial_i\gamma_i+\sum_{j\neq i}\gamma_jh_{ji}-\beta^i+cv_i\varphi,v\>X_i,
$$
%\end{equation}
hence  (\ref{eq:cconst}) reduces to $(ii)$ and $(iii)$.\qed

\begin{remark}\po\emph{Lemma \ref{prop:fundconst} also holds for an isometric immersion $f\colon M^n(c)\rightarrow \mathbb Q^{n+p}_s(\tilde c)$, with $c\neq \tilde c$,  satisfying the assumptions of Proposition \ref{teo:curvdif}, if system $\ral$ is replaced by
$$
\Ral^*=\left\{
\begin{array}{ll}
i) \ \partial_i\varphi=v_i\gamma_i, \,\,\,\, ii) \ \partial_i\gamma_j= h_{ji}\gamma_i, \,\,i\neq j\\
iii) \ \partial_i\gamma_i =- \sum_{j\neq i}h_{ji} \gamma_j-((L^t)^{-1}-I)\sum_r\epsilon_rV_{ir}\beta^r&-\left((L^t)^{-1}(c-\tilde c)+\tilde c\right)v_i\varphi,\\
 iv) \ \partial_i\beta^r=-V_{ir}\gamma_i,\,\,\,i\neq r.\\
\end{array}
\right.
$$}
\end{remark}

\begin{definition}\po\label{df:admissible} {\em Let $V,W_1,W_2$ be  Euclidean vector spaces.
% endowed with inner products. 
Given   $A\in V^*\otimes V$, $c, \tilde c\in \R$, $\psi\in V^*$, $\nu\in V^*\otimes  W_1$ and $\beta\in V^*\otimes W_2$, we say that $(\psi,\nu,\beta)$ is \emph{A-admissible} if the system
% of equations 
\begin{equation}\label{eq:matrixsystemb}
\left\{ \begin{array}{l} X+X^t=\nu^t\nu+\beta^t\beta+\tilde c\psi^t\psi\vspace{1ex}\\
XA+A^tX^t=\beta^t\beta-(c-\tilde c)\psi^t\psi\end{array}\right.
\end{equation}
admits a unique solution $X\in V^*\otimes V$, and if such solution  is invertible.}
\end{definition}
\begin{remark}\po {\em  Notice that if  the inner product $\<\cdot, \cdot\>$ is replaced by  $\<\cdot, \cdot\>^\sim$, given by $\<\cdot, \cdot\>^\sim=\<\cdot, B\cdot\>$ for some invertible $B\in V^*\otimes V$, then $(\psi,\nu,\beta)$ is $A$-admissible if and only if
$(\hat\psi,\hat\nu,\hat\beta)$, given by $\hat\psi=\psi B$, $\hat\nu=\nu B$ and $\hat\beta=\beta B$, is $\hat A$-admissible, where $\hat A=B^{-1}AB$. Namely, 
system (\ref{eq:matrixsystemb}) has a unique solution $X$ if and only if the corresponding system for $(\hat\psi,\hat\nu,\hat\beta, \hat A)$ has a unique solution  $\hat X$, and these are related by $\hat X=XB$.}
\end{remark}

\begin{theorem}\po\label{lem:princonst}
Let $f\colon M^n(c)\rightarrow \mathbb Q^{n+p}(\tilde c)$ be an isometric immersion  
as in either of Propositions \ref{thm:hiebetaij} or  \ref{teo:curvdif}, depending on whether $c=\tilde c$ or $c\neq \tilde c$, respectively. Given $x_0\in M^n(c)$,  if $L$ is an invertible endomorphism of a Euclidean vector space  $V$  and 
$(\psi_0,\nu_0,\beta_0)$ is an $L$-admissible  triple, with $\nu_0\in V^*\otimes  T_{x_0}M$ and $\beta_0\in V^*\otimes  N_fM(x_0)$, then there exist an open neighborhood $W$ of $x_0$ and a unique $L$-transform $\tilde f=\ral_{\varphi,\beta,\Omega, L}(f|_W)$ of $f|_W$ such that $\varphi(x_0)=\psi^t_0$, $\omega(x_0)=\nu^t_0$ and $\beta(x_0)=\beta_0$.
\end{theorem}
\proof  Let $\varphi^0,\gamma_1^0,...,\gamma_n^0,\beta^1_0,...,\beta_0^p \in V$ be given by $\varphi^0=\psi_0^t$, $\gamma_i^0=v_i^{-1}\nu_0(\partial_i)$ and $\beta^r_0=\beta_0^t(\xi_r)$. It is immediate to verify, using the fact that $(v, h)$ (resp.,  $(v, h, V)$)
is a solution of (\ref{eq:i}) (resp., (\ref{eq:ii})),
that the compatibility conditions of  system ${\cal R}$ (resp., ${\cal R}^*$)  are satisfied. Thus there exists  a solution $(\varphi,\gamma_1,...,\gamma_n,\beta^1,...,\beta^p)$  of  $\ral$ (resp.,  $\ral^*$)  such that $\varphi(x_0)=\varphi^0$, $\gamma_i(x_0)=\gamma_i^0$ and $\beta^r(x_0)=\beta^r_0$ for $1\leq i\leq n$ and $1\leq r\leq p.$ It follows from  Lemma \ref{prop:fundconst} that $\omega=d\varphi$ and $\beta$ given by (\ref{eq:omegabeta}) 
satisfy (\ref{eq:alphao}), (\ref{eq:cconst}) and the initial conditions 
$$\varphi(x_0)=\psi^t_0,\,\,\beta(x_0)=\beta_0\,\,\,\mbox{and}\,\,\,\omega(x_0)=\nu^t_0.$$
 Since the triple $(\psi_0,\nu_0,\beta_0)$ is $L$-admissible, the  system  of  equations
$$
\left\{ \begin{array}{l} X+X^t=\nu_0^t\nu_0+\beta_0^t\beta_0+\tilde c\psi_0^t\psi_0\vspace{1ex}\\
XL+L^tX^t=\beta_0^t\beta_0-(c-\tilde c)\psi_0^t\psi_0
\end{array}\right.
$$
has a unique solution $X=\Omega_0$, which 
is invertible. Let $\Omega\in \Gamma(V^*\otimes V)$ be the unique solution of (\ref{eq:o1b}) with $\Omega(x_0)=\Omega_0$ on an open simply connected neighborhood $W$ of $x_0$ where $\Omega$ remains invertible. 
Using (\ref{eq:alphao}) and (\ref{eq:o1c}) we obtain 
$$\begin{array}{l}d(\Omega L+L^t\Omega^t-\rho)(X)=\o\Phi(X)L+L^t\Phi(X)^t\o^t
-d\beta^t(X)\beta-\beta^td\beta(X)\vspace{1ex}\\\hspace*{25ex}+(c-\tilde c)\o(X)\va^t+(c-\tilde c)\va\o(X)^t
\vspace{1ex}\\\hspace*{23.4ex}=\o(\Phi(X)L+A(X)\beta+(c-\tilde c)X\va^t)\vspace{1ex}\\\hspace*{25.5ex}
+(\o(\Phi(X)L+A(X)\beta+(c-\tilde c)X\va^t))^t,
\end{array}
$$
and the right-hand-side vanishes by (\ref{eq:cconst}). Thus  (\ref{eq:cond3}) holds  on $W$ since it holds at $x_0$.

It follows that $(\varphi,\beta,\Omega)$ is the unique triple that satisfies the conditions of Definitions~\ref{ribvet} and \ref{df:ribconst} on $W$, as well as the given initial conditions at $x_0$. Therefore $\tilde f=\ral_{\varphi,\beta, \Omega}(f|_{W})$ is an $L$-transform of $f|_{W}$, and it is the unique one defined on an open neighborhood of $x_0$  which 
satisfies the given initial conditions at $x_0$.\vspace{1ex}\qed

It remains to investigate  when   $(\psi,\nu,\beta)$ is an $L$-admissible triple for a given endomorphism $L$, and  in particular for which  endomorphisms $L$ there exist  $L$-admissible triples. 
First, in the next result  we give conditions on $L$ and $(\psi,\nu,\beta)$  under which system (\ref{eq:matrixsystemb}) has exactly one solution.

An operator $A\in V^*\otimes V$ is said to be \emph{nonderogatory} if its minimal and characteristic polynomials coincide, or equivalently, if for each eigenvalue $\alpha$ of $A$ 
the eigenspace $\ker (A^c-\alpha I)$ of the complexified endomorphism $A^c\in (V^c)^*\otimes V^c$, $V^c=V\otimes \mathbb{C}$, has dimension one. If $a_{11}, \ldots, a_{1n_1}, \ldots, a_{p1}, \ldots, a_{pn_p}$ is the Jordan basis of $A^c$, with $a_{i1}, \ldots, a_{in_i}$ corresponding to the eigenvalue $\alpha_i$, we call  $a_{i1}, \ldots, a_{in_i}$ the \emph{generalized eigenvectors} associated to $\alpha_i$. Thus, for $1\leq i\leq p$ and $1\leq k_i\leq n_i$ we have
$(A^c-\alpha_i I)a_{ik_i}=a_{i,k_i-1},\,\,\,\,a_{i0}=0.$

\begin{proposition}\po\label{teo:matrixsystem2}
Let $A\in V^*\otimes V$ be nonderogatory, let $\alpha_1, \ldots, \alpha_p$ (resp., $\gamma_1, \ldots, \gamma_q$) be the real (resp., complex) eigenvalues of $A$,  let  $a_{i1},a_{i2},...,a_{in_i}$ 
(resp.,  $w_{j1},w_{j2},...,w_{jm_j}$) be the generalized eigenvectors of $A$ associated with $\alpha_i$, $1\leq i\leq p$  
(resp., $\gamma_j$, $1\leq j\leq q$ ). Then system (\ref{eq:matrixsystemb}) has a solution $X\in V^*\otimes V$ if and only if
\begin{equation}\label{eq:autcom0a}
\begin{array}{l}
\sum_{r=0}^{k_i-1}\langle a_{i,k_i-r},B a_{i,r+1}\rangle=0, \ i=1,...,p \ \mbox{and} \ k_i=1,2,...,n_i
\end{array}
\end{equation}
and
\begin{equation}\label{eq:autcom1a}
\begin{array}{l}
\sum_{r=0}^{\ell_j-1}\langle w_{j,\ell_j-r},B^c\bar{w}_{j,r+1}\rangle=0, \ \ j=1,...,q \ \mbox{and} \ \ell_j=1,2,...,m_j,
\end{array}
\end{equation}
where 
$$B=(I-A^t)\beta^t\beta-A^t\nu^t\nu-((c-\tilde c)I+\tilde cA^t)\psi^t\psi,$$
$B^c$ is its complexification and $\langle \;, \;\rangle\colon V^c\times V^c \to \mathbb{C}$ is the canonical hermitian inner product on $V^c=V\otimes \mathbb{C}$. Moreover, if (\ref{eq:autcom0a}) and (\ref{eq:autcom1a}) are satisfied then the solution $X$ 
is unique.
\end{proposition}
\proof Let $X_s=\frac{1}{2}(X+X^t)$ and $X_a=\frac{1}{2}(X-X^t)$ be the symmetric and anti-symmetric parts of $X\in V^*\otimes V$, respectively. Then  $X$ is a solution of  (\ref{eq:matrixsystemb}) if and only if 
\be\label{eq:xs}
2X_s=\nu^t\nu+\beta^t\beta+\tilde c\psi^t\psi 
\ee
and
\be\label{eq:xa} 
X_aA-A^tX_a=B_s=(B+B^t)/2.
\ee
Let $\Psi\colon {\cal A}(V)\to {\cal S}(V)$ be the linear map from the space of anti-symmetric endomorphisms of $V$ into the 
space of symmetric endomorphisms of $V$ given by
$$\Psi(Y)=YA-A^tY.$$
A necessary and sufficient condition for $\Psi$ to be  injective is
that $A$ be a nonderogatory endomorphism (see \cite{TZ}).  Therefore, a  solution $X$ of  (\ref{eq:matrixsystemb}) exists if and only if $B_s$  belongs to the range of $\Psi$, and in this case $X$ is uniquely determined by $X=X_s+X_a$, with $X_s$ given by (\ref{eq:xs}) and $X_a$ the unique element of ${\cal A}(V)$ such that $\Psi(X_a)=B_s$.

Now, $B_s$  belongs to the range of $\Psi$ if and only if  its complexification $B_s^c$ belongs to the range of $\Psi$ regarded as a linear map from ${\cal A}(V^c)$ to ${\cal S}(V^c)$, where $V^c$ is endowed with the canonical hermitian inner product and  ${\cal S}(V^c)$, ${\cal A}(V^c)$ denote the spaces of self-adjoint and anti-self-adjoint endomorphisms of $V^c$, respectively. This is the case  if and only if $B_s$, and hence $B$,  belongs to the orthogonal complement of the kernel of the adjoint $\Psi^*$ of $\Psi$ with respect to the inner product on $V^c\otimes V^c$ given by
$$\<X, Y\>=\mbox{trace}Y^*X$$
where $Y^*$ stands for the adjoint of $Y$. It is easily checked that
$$\Psi^*Y=YA^t-AY$$
and that a basis of $\ker \Psi^*$ is  $\{Y_{ik_i}, Z_{j\ell_j}, \bar Z_{j\ell_j}, \,\,1\leq k_i\leq  n_i,\,\,1\leq \ell_j\leq m_j\}$, where
$$Y_{ik_i}=a_{ik_i}a_{i1}^t+a_{i,k_i-1}a_{i2}^t+\ldots, a_{i1}a^t_{ik_i}, \,\,\,1\leq i\leq  p\,\,\\,\mbox{and}\,\,\,1\leq k_i\leq  n_i,$$
and
$$Z_{j\ell_j}=w_{j\ell_j}w_{j1}^t+w_{j,\ell_j-1}w_{j2}^t+\ldots, w_{j1}w^t_{j\ell_j}, \,\,\,1\leq j\leq  q\,\,\ \,\mbox{and}\,\,\, 1\leq \ell_j \leq  m_j.$$
(cf. \cite{K}). Therefore $B^c_s\in \Psi({\cal A}(V^c))$ if and only if
$\<B^c, Y_{ik_i}\>=0=\<B^c, Z_{j\ell_j}\>$
for all $1\leq i\leq p$, $1\leq k_i\leq n_i$, $1\leq j\leq q$,  $1\leq \ell_j\leq m_j$, which are equivalent to (\ref{eq:autcom0a}) and (\ref{eq:autcom1a}).\qed

\begin{remark}\po {\em If (\ref{eq:autcom0a}) and (\ref{eq:autcom1a}) are satisfied then the unique solution $X$ of (\ref{eq:matrixsystemb}) 
is $X=X_s+X_a$, where $X_s$ is given by (\ref{eq:xs}) and $X_a$ is the unique solution of (\ref{eq:xa}), an explicit expression of which can be found in \cite{M}.}
\end{remark}

Next, for $A$ and $(\psi,\nu,\beta)$ as in the preceding proposition, we give further sufficient conditions for the unique solution 
$X$ of (\ref{eq:matrixsystemb}) to be invertible, that is, for $(\psi,\nu,\beta)$ to be $A$-admissible, under the assumption that $(c,\tilde c)$ belongs to the subset
$$D(c, \tilde c)=\R^2\setminus \{(c, \tilde c)\,:\, \tilde c< 0\,\,\mbox{and}\,\,c\in [\tilde c, 0]\}.$$

\begin{proposition}\label{lem:invertsol}\po Let $A\in V^*\otimes V$ be  nonderogatory and let $(\psi,\nu,\beta)$ satisfy conditions (\ref{eq:autcom0a}) and (\ref{eq:autcom1a}).  Set 
\be\label{eq:s}S= \left\{ \begin{array}{l} \ker \nu \cap \ker \beta,\,\,\,\mbox{if}\,\,\,\,(c, \tilde c)= (0,0);\\
\ker \psi\cap \ker \nu \cap \ker \beta, \,\,\,\mbox{if}\,\,\,\,(c, \tilde c)\neq (0,0)
\end{array}\right.
\ee
and $S^c=S\otimes \mathbb{C}$. If $(c, \tilde c)\in D(c, \tilde c)$, then
 the unique solution $X\in V^*\otimes V$ of (\ref{eq:matrixsystemb}) satisfies 
$$\ker X=\ker X^t\subset S\,\,\,\,\mbox{and}\,\,\,\,A(\ker X)\subset \ker X.$$
 In particular,  if  
$$E_{\alpha_i}\cap S=\{0\}= E_{\gamma_j}\cap S^c$$
for all $1\leq i\leq p$ and $1\leq j\leq q$ then $X$ is invertible, and hence
 $(\psi,\nu,\beta)$ is $A$-admissible.
\end{proposition}
\proof
Subtracting the second equation from the first one in (\ref{eq:matrixsystemb}), we obtain
$$
X(I-A)+(I-A^t)X^t=\nu^t\nu+c\psi^t\psi.
$$ 
Therefore, for all $u\in \ker X^t$ we have
$$
0=|\beta u|^2+|\nu u|^2+\tilde c \psi^2(u)=|\beta u|^2-(c-\tilde c)\psi^2(u)=|\nu u|^2+c\psi^2(u).
$$
For $(c, \tilde c)\in D(c, \tilde c)$, this implies  that  
$\ker X^t\subset S$. We obtain  from (\ref{eq:matrixsystemb})  that $\ker X=\ker X^t$ and  $A(\ker X)\subset \ker X$.\vspace{1ex}
\qed

The following proposition provides conditions on a nonderogatory endomorphism $A$ for an A-admissible triple $(\psi,\nu,\beta)$ to exist.

\begin{proposition}\label{p:ladm} \po Let $A\in V^*\otimes V$ be nonderogatory. If $(c, \tilde c)\in D(c, \tilde c)\setminus \{(0,0)\}$, then 
 there  exists an A-admissible triple $(\psi,\nu,\beta)$  if and only if all  real eigenvalues of $A$ with odd algebraic multiplicity belong to the subset $Z(c,\tilde c)$ defined as either $[\kappa,1]$, $[0,1]$, $[0,\kappa]$, $(-\infty,\kappa]\cup [0,\infty)$, $(-\infty,1]\cup [\kappa,\infty)$,  $(-\infty,1]$ or $[0, \infty)$, depending on whether $c\geq\tilde c>0$, $0\leq c<\tilde c$, $c<0<\tilde c$,  $c<\tilde c<0$, $\tilde c<0<c$, $\tilde c=c<c$ or 
 $c<0=\tilde c$, respectively, with $\kappa=(\tilde{c}-c)/\tilde c$.
 If  $(c, \tilde c)=(0,0)$ or  $(c, \tilde c)\not\in D(c, \tilde c)$,  such a triple always exists.
\end{proposition}
\proof Let  $\alpha_1, \ldots, \alpha_p$ (resp., $\gamma_1, \ldots, \gamma_q$) be the real (resp., complex)  eigenvalues of the nonderogatory endomorphism $A$,   and  $a_{i1},a_{i2},...,a_{in_i}$ 
(resp.,  $w_{j1},w_{j2},...,w_{jm_j}$)  the generalized eigenvectors of $A$ associated with $\alpha_i$, $1\leq i\leq p$ (resp., $\gamma_j$, $1\leq j\leq q$ ). 
We denote by $R_{w_{jl}}$ and $I_{w_{jl}}$  the real and imaginary parts of $w_{jl}$, respectively.
We will show that, if $(c, \tilde c)\in D(c, \tilde c)\setminus \{(0,0)\}$, then one can choose $y_{ik_i}, y_{j\ell_j}^R, y_{j\ell_j}^I\in \R$, $Y_{ik_i}, Y_{j\ell_j}^R, Y_{j\ell_j}^I\in W_1$
and $\xi_{ik_i}, \xi_{j\ell_j}^R, \xi_{j\ell_j}^I\in W_2$, $1\leq k_i\leq n_i$, $1\leq \ell_j\leq m_j$, so that  the triple $(\psi, \nu, \beta)$ defined by
$\psi(a_{ik_i})=y_{ik_i}$, $\psi(R_{w_{j\ell_j}})=y_{j\ell_j}^R$,  $\psi(I_{w_{j\ell_j}})=y_{j\ell_j}^I$,  $\nu(a_{ik_i})=Y_{ik_i}$,
$\nu(R_{w_{j\ell_j}})=Y_{j\ell_j}^R$, $\nu(I_{w_{j\ell_j}})=Y_{j\ell_j}^I$, $\beta(a_{ik_i})=\xi_{ik_i}$,  $\beta(R_{w_{j\ell_j}})=\xi_{j\ell_j}^R$ and 
$\beta(I_{w_{j\ell_j}})=\xi_{j\ell_j}^I$   is $A$-admissible if and only if $A$ satisfies the conditions in the statement. If  
$(c, \tilde c)=(0,0)$ or  $(c, \tilde c)\not\in D(c, \tilde c)$, we will show that such a choice is always possible.

\begin{claim}\label{aff3}\po  One can choose $y_{j\ell_j}^R, y_{j\ell_j}^I\in \R$,  $Y_{j\ell_j}^R, Y_{j\ell_j}^I\in W_1$ and $\xi_{j\ell_j}^R, \xi_{j\ell_j}^I\in W_2$, $1\leq j\leq q$,  $1\leq \ell_j\leq m_j,$ so that 
(\ref{eq:autcom1a}) is satisfied  and 
either $(y_{j1}^R,Y_{j1}^R, \xi_{j1}^R)$ or $(y_{j1}^I,Y_{j1}^I, \xi_{j1}^I)$ is nontrivial. 
%$(y_{j1}^R,Y_{j1}^R, \xi_{j1}^R)\neq (0,0,0)$ or $(y_{j1}^I,Y_{j1}^I, \xi_{j1}^I)\neq (0,0,0)$.
\end{claim}
\proof For $\ell_j=1$, denoting by $R_{\gamma_j}$ and $I_{\gamma_j}$ the real and imaginary parts of $\gamma_j$, respectively,  equation (\ref{eq:autcom1a}) becomes
\be\label{eq:18a}\begin{array}{l}
\hspace*{5ex}(1-R_{\gamma_j})\left(|\xi^R_{j1}|^2-|\xi^I_{j1}|^2\right)+2I_{\gamma_j}\left<\xi^R_{j1},\xi^I_{j1}\right>-R_{\gamma_j}\left(|Y^R_{j1}|^2-|Y^I_{j1}|^2\right)\vspace{0.1cm}\\
\hspace*{10ex}+2I_{\gamma_j}\left<Y^R_{j1},Y^I_{j1}\right>+\left((c-\tilde c)+R_{\gamma_j}\tilde c\right)\left((y_{j1}^I)^2-(y_{j1}^R)^2)\right)
+2\tilde cI_{\gamma_j}y_{j1}^Ry_{j1}^I=0\\
\end{array}
\ee
and 
\be\label{eq:18b}\begin{array}{l}
\hspace*{5ex} -I_{\gamma_j}\left(|\xi^R_{j1}|^2-|\xi^I_{j1}|^2\right)   +2(1-R_{\gamma_j})\left<\xi^R_{j1},\xi^I_{j1}\right>   -I_{\gamma_j}\left(|Y^R_{j1}|^2-|Y^I_{j1}|^2\right)\vspace{0.1cm}\\
\hspace*{10ex} - 2R_{\gamma_j}\left<Y^R_{j1},Y^I_{j1}\right>+\tilde cI_{\gamma_j}\left((y_{j1}^I)^2-(y_{j1}^R)^2\right)
+2((c-\tilde c)+\tilde c R_{\gamma_j})y_{j1}^Ry_{j1}^I=0.
\end{array}
\ee
It is easy to see that, for each $1\leq j\leq q$,  one can choose a solution $(\xi_{j1}^R, \xi^I_{j1}, Y^R_{j1}, Y^I_{j1}, y^R_{j1}, y_{j1}^I)$ of both of the preceding equations so that  $(y_{j1}^R,Y_{j1}^R, \xi_{j1}^R)\neq (0,0,0)$ or $(y_{j1}^I,Y_{j1}^I, \xi_{j1}^I)\neq (0,0,0)$.
Once $\xi_{j\ell_j}^R, \xi^I_{j\ell_j}, Y^R_{j\ell_j}, Y^I_{j\ell_j}, y^R_{j\ell_j}, y_{j\ell_j}^I$ have been chosen for
 $1\leq \ell_j\leq k_j-1\leq m_j-1$, then (\ref{eq:autcom1a}) for $\ell_j=k_j$ becomes a pair of linear equations on 
$(\xi_{jk_j}^R, \xi^I_{jk_j}, Y^R_{jk_j}, Y^I_{jk_j}, y^R_{jk_j}, y_{jk_j}^I)$, which are easily checked to admit a common solution.\qed
%$1\leq j\leq q$ and $1\leq \ell_j\leq m_j-1$, equation (\ref{eq:autcom1a}) for $\ell_j=k_j$ becomes a pair of linear equations on 
%$(\xi_{jk_j}^R, \xi^I_{jk_j}, Y^R_{jk_j}, Y^I_{jk_j}, y^R_{jk_j}, y_{jk_j}^I)$, which are easily checked to admit a common solution. \qed

\begin{claim}\label{aff1}\po If  $(c, \tilde c)\in D(c,\tilde c)\setminus  \{(0,0)\}$ and $\alpha_i\in Z(c,\tilde c)$ then   
one can choose $y_{ik_i}, Y_{ik_i}, \xi_{ik_i},$  $1\leq i\leq p, \  1\leq k_i\leq n_i$, such that  (\ref{eq:autcom0a}) is satisfied and $(y_{i1}, Y_{i1}, \xi_{i1})\neq (0,0,0)$. The same statement is true if  $(c,\tilde c)=(0,0)$.
\end{claim}
\proof  Given $1\leq i\leq p$,  equation  (\ref{eq:autcom0a}) for $k_i=1$ becomes
\begin{equation}\label{eq:quadr}
(1-\alpha_i)| \xi_{i1}|^2-\alpha_i|Y_{i1}|^2-\left((c-\tilde c)+\alpha_i\tilde c\right)(y_{i1})^2=0.
\end{equation} 
One can choose  $(y_{i1}, Y_{i1}, \xi_{i1})\neq (0,0,0)$ satisfying 
(\ref{eq:quadr}) if the coefficients $(1-\alpha_i)$, $-\alpha_i$ and $-((c-\tilde c)+\alpha_i\tilde c)$ are neither all positive nor all negative, which is easily seen to be the case if  $(c, \tilde c)\in D(c,\tilde c)\setminus \{(0,0)\}$ and $\alpha_i\in Z(c,\tilde c)$.
Once $(y_{i1}, Y_{i1}, \xi_{i1})\neq (0,0,0)$ satisfying 
(\ref{eq:quadr}) has been chosen, it is easy to check that $y_{ik_i},$ $Y_{ik_i},$ $\xi_{ik_i}$,  
$2\leq k_i\leq n_i$, can be chosen so that the remaining equations in  (\ref{eq:autcom0a}) are satisfied.\qed

\begin{claim}\label{aff2}\po Suppose that  $(c, \tilde c)\in D(c,\tilde c)\setminus  \{(0,0)\}$ and $\alpha_i\notin Z(c,\tilde c)$ for some $i\in \{1, \ldots, p\}$. If $\alpha_i$ has even algebraic multiplicity then $y_{il_i}=Y_{il_i}=\xi_{il_i}=0$ for $l_i\leq n_i/2$, and one can choose $y_{il_i}, Y_{il_i}, \xi_{il_i},$  $n_i/2+1\leq l_i\leq n_i$, so that $a_{i1}\not\in \ker X$, where $X$ is the unique solution of (\ref{eq:matrixsystemb}).
\end{claim}
\proof Since  $\alpha_i\notin Z(c,\tilde c)$, the numbers  $(1-\alpha_i), \ -\alpha_i, \ -((c-\tilde c)+\alpha_i\tilde c)$ are either all positive or all negative. Since  $n_i$ is even, it is easy to check that a solution $(y_{il_i},Y_{il_i},\xi_{il_i})$  of (\ref{eq:autcom0a}) must satisfy $y_{il_i}=Y_{il_i}=\xi_{il_i}=0$ for $l_i\leq n_i/2$, and that   $y_{il_i},$ $Y_{il_i}$ and $\xi_{il_i}$ may be arbitrarily chosen  for $l_i>n_i/2$. We choose a solution such that $(y_{il_i},Y_{il_i},\xi_{il_i})\neq (0,0,0)$ for  $l_i=n_i/2+1.$
Using that $A$ is nonderogatory, we obtain that the unique solution $X$ of (\ref{eq:matrixsystemb}) satisfies $\left<X a_{i1},a_{il_i}\right>=\left<X a_{i1},a_{jk_j}\right>=\left<X a_{i1},{R_{w_{j\ell_j}}}\right>=\left<X a_{i1},I_{w_{j\ell_j}}\right>=0$ 
 for $1\leq l_i\leq n_i-1$, $1\leq k_j\leq n_j$ and $1\leq \ell_j\leq m_j$, hence 
 \be\label{eq:aini} 
 X(a_{i1})\in \spa\{a_{in_i}\}.
 \ee 
Moreover,
\begin{equation}\label{eq:linda}
\begin{array}{lll}
0&=&\left<(X_a A- A^tX_a-B_s) a_{il_i},a_{ik_i}\right>=\left<X_a a_{i,l_i-1},a_{ik_i}\right>-\left<X_a a_{il_i},a_{i,k_i-1}\right>\vspace{0.1cm}\\
&&-\left<B_s a_{il_i},a_{ik_i}\right>, \ 1\leq l_i,k_i\leq n_i.\\
\end{array}
\end{equation}
Using that $B_s(a_{il_i})=0$ for $l_i\leq n_i/2$, we obtain from the preceding  relations that
$$
\begin{array}{l}
\left<X a_{i1},a_{in_i}\right>=\left<X_a a_{i1},a_{in_i}\right>=\left<X_aa_{i,\frac{n_i}{2}},a_{i,\frac{n_i}{2}+1}\right>=\left<B_s a_{i,\frac{n_i}{2}+1},a_{i,\frac{n_i}{2}+1}\right>/2,\\
\end{array}
$$
 Since  $(y_{i,\frac{n_i}{2}+1}, Y_{i,\frac{n_i}{2}+1}, \xi_{i,\frac{n_i}{2}+1})\neq (0,0,0)$ and $\alpha_i\notin Z(c,\tilde c)$, we have 
$$
\left<{B_s}a_{i,\frac{n_i}{2}+1},a_{i,\frac{n_i}{2}+1}\right>=(1-\alpha_i)|\xi_{i,\frac{n_i}{2}+1}|^2-\alpha_i|Y_{i,\frac{n_i}{2}+1}|^2-((c-\tilde c)-\alpha_i\tilde c)y^2_{i,\frac{n_i}{2}+1}\neq 0.
$$
Thus $a_{i1}\notin ker(X)$. \vspace{1ex}\qed

It follows from Proposition \ref{lem:invertsol} and Claims \ref{aff3} -\ref{aff2} that, when $(c,\tilde c)\in D(c,\tilde c)\setminus \{(0,0)\}$, the real eigenvalues of $A$ with odd algebraic multiplicity belonging to $Z(c,\tilde c)$ is  a sufficient  condition  for  an $A$-admissible triple to exist. To complete the proof of Proposition \ref{p:ladm} when $(c,\tilde c)\in D(c,\tilde c)\setminus \{(0,0)\}$, it remains to show  that  this condition is also necessary. 
In fact, if  $\alpha_i\not\in Z(c,\tilde c)$ is such that $n_i$ is odd, it is easy to check that a solution $(y_{il_i},Y_{il_i},\xi_{il_i})$  of (\ref{eq:autcom0a}) must satisfy $y_{il_i}=Y_{il_i}=\xi_{il_i}=0$ for $l_i\leq (n_i+1)/2$, and that   $y_{il_i},$ $Y_{il_i}$ and $\xi_{il_i}$ may be arbitrarily chosen  for $l_i> (n_i+1)/2$.  
If $X$ is the unique solution of (\ref{eq:matrixsystemb}), using that
$B_s(a_{il_i})=0$ for $l_i\leq (n_i+1)/2$ we obtain from  (\ref{eq:linda}) that
$$
\left<X_a a_{i1},a_{in_i}\right>=\left<X_aa_{i,(n_i+1)/2},a_{i,(n_i+1)/2}\right>=0.
$$
In view of (\ref{eq:aini}), this implies that $a_{i1}\in \ker X$.\vspace{1ex}

Now assume that  $(c,\tilde c)\notin D(c,\tilde c)$ and suppose  that $\dim W_2\leq \dim W_1$, the argument for the case
$\dim W_1 < \dim W_2$ being similar. Choose $\beta\in V^*\otimes  W_2$ so that  $R_{w_{j1}}\notin \ker \beta$, $I_{w_{j1}}\notin \ker \beta$, and so that $a_{i1}\in \ker \beta$ only if $\alpha_i=\kappa=(\tilde c- c)/\tilde c$. Let $\nu \in V^*\otimes  W_1$ be given by $\nu=\lambda P\circ \beta$, 
where $P\in W^*_2\otimes  W_1$ satisfies $P^tP=I$ and $\lambda \in \mathbb R/\{0\}$. Under these conditions, we have $\ker\nu=\ker \beta$ and $\nu^t\nu=\lambda^2\beta^t\beta$.

Again, we will show that one can choose $y_{ik_i}, y_{j\ell_j}^R, y_{j\ell_j}^I\in \R$ so that  the triple $(\psi, \nu,\beta)$,  with $\psi$ defined by
$\psi(a_{ik_i})=y_{ik_i}$, $\psi(R_{w_{j\ell_j}})=y_{j\ell_j}^R$ and $\psi(I_{w_{j\ell_j}})=y_{j\ell_j}^I$,  is $A$-admissible. 

Denoting $\beta(a_{i1})=\xi_{i1}$,  equation (\ref{eq:autcom0a}) for $k_i=1$ becomes 
$$
(1-\alpha_i(1+\lambda^2))|\xi_{i1}|^2-((c-\tilde c)+\alpha_i\tilde c)(y_{i1})^2=0,
$$
which  always admits a nontrivial solution $(y_{i1}, \xi_{i1})\neq (0,0)$  if the coefficients $(1-\alpha_i(1+\lambda^2))$ and $-((c-\tilde c)+\alpha_i\tilde c)$   are both nonzero and have opposite signs,  which we can assume to be the case by an appropriate choice of $\lambda$. It also admits a nontrivial solution if at least one such coefficient is zero. Notice that when
 $-((c-\tilde c)+\alpha_i\tilde c)=0$, that is, for  $\alpha_i=\kappa$, this is due to the fact that we have chosen  $a_{i1}\in \ker \beta$ when $\alpha_i=\kappa$.

Denoting $\beta(a_{ik_i})=\xi_{ik_i}$, $1\leq i\leq p$ and $2\leq k_i\leq n_i$, it is easy to check that $y_{ik_i}$ can be defined  so that the remaining equations in (\ref{eq:autcom0a}) are satisfied. For (\ref{eq:autcom1a}), it is sufficient to replace $\nu^t\nu=\lambda^2\beta^t\beta$ in (\ref{eq:18a}) and (\ref{eq:18b}),  denoting $\beta(R_{w_{j\ell_j}})=\xi_{j\ell_j}^R$ and 
$\beta(I_{w_{j\ell_j}})=\xi_{j\ell_j}^I$ for $1\leq j\leq q$ and $1\leq \ell_j\leq m_j$, and proceed as in the proof of that  claim. In this way, one can define $\psi$ in the Jordan basis of $A$ so that (\ref{eq:autcom0a}) and (\ref{eq:autcom1a}) are satisfied and 
\begin{equation}\label{invxcr2}
E_{\alpha_i}\cap S=\{0\}=E_{\gamma_j}\cap S^c
\end{equation}
for all $1\leq j\leq q$ and $1\leq \ell_j\leq m_j$, where $S=\ker \beta \cap \ker \psi$ and $S^c=S\otimes \mathbb C$.

Now, for a solution $X$ of (\ref{eq:matrixsystemb}),  for all $u\in \ker X^t$ we have
$$
0=(1+\lambda^2)|\beta u|^2+\tilde c \psi^2(u)=|\beta u|^2-(c-\tilde c)\psi^2(u)=\lambda^2|\beta u|^2+c\psi^2(u).
$$
Choosing $\lambda$ so that  $\lambda^2+1\neq \frac{1}{\kappa}$, the above system only admits the trivial solution $|\beta u|=\psi(u)=0$, and hence
%\begin{equation}\label{invxcr}
$\ker X^t\subset S$.
%\end{equation}
Then, we obtain from (\ref{eq:matrixsystemb}) that $\ker X= \ker X^t$ and $A(\ker X)\subset \ker X$. It follows from (\ref{invxcr2}) that 
$\ker X=\{0\}$, and hence $(\psi, \nu, \beta)$ is $A$-admissible.
\vspace{1ex}\qed

In case $A\in V^*\otimes V$ is a symmetric endomorphism with $k=\dim V$ distinct eigenvalues (see part $(ii)$ of Remarks \ref{re:ldecomp} below), 
one has the following more precise statement. 

\begin{proposition}\label{p:Asym}\po Let $A\in V^*\otimes V$ be symmetric with $k=\dim V$ distinct eigenvalues $\alpha_1, \ldots, \alpha_k$. Let
$a_1, \ldots, a_k$ be eigenvectors of $A$ associated to $\alpha_1, \ldots, \alpha_k$, respectively. Then system (\ref{eq:matrixsystemb}) has a solution if and only if
\be\label{eq:symcond} (1-\alpha_i)|\beta(a_i)|^2-\alpha_i|\nu(a_i)|^2- ((c-\tilde c) +\alpha_ic)\psi(a_i)^2=0,\,\,\,\,1\leq i\leq k.\ee
If these equations are satisfied, then  the solution of (\ref{eq:matrixsystemb}) is unique and given by
\be\label{eq:explicitsol}\<Xa_i, a_i\>= \frac{1}{2}G_{ii}\;\;\;\mbox{and}\;\;\;\<Xa_i, a_j\>=\frac{1}{\alpha_i-\alpha_j}(\rho_{ij}-\alpha_jG_{ij}),\,\,\,i\neq j,\ee
where
$$\rho_{ij}=\<\beta(a_i), \beta(a_j)\>-(c-\tilde c)\psi(a_i)\psi(a_j)$$
and
$$G_{ij}=\<\nu(a_i), \nu(a_j)\>+\<\beta(a_i), \beta(a_j)\>+\tilde c\psi(a_i)\psi(a_j).$$
Moreover, if 
$(c, \tilde c)\in D(c, \tilde c)$ then $X$ is invertible, that is, $(\psi, \nu, \beta)$ is $A$-admissible, if and only if 
$(\psi(a_i),  \nu(a_i),\beta(a_i))\neq (0,0,0)$ for all $1\leq i\leq k$.
Furthermore, there exists an $A$-admissible  triple $(\psi, \nu, \beta)$ 
if and only if 
$\alpha_i\in Z(c,\tilde c)$ for all $1\leq i\leq k$. If $(c, \tilde c)=(0,0)$ or $(c, \tilde c)\not\in D(c, \tilde c)$,  an  A-admissible triple $(\psi,\nu,\beta)$ always exists.
\end{proposition}
\proof For a symmetric endomorphism $A\in V^*\otimes V$ with $k$ distinct eigenvalues $\alpha_1, \ldots, \alpha_k$, it is easily checked that  equations (\ref{eq:autcom0a})  reduce to (\ref{eq:symcond}). Thus the first assertion follows from 
Proposition \ref{teo:matrixsystem2}. Equations (\ref{eq:explicitsol}) follow directly from system (\ref{eq:matrixsystemb}), hence they
necessarily provide the unique solution of (\ref{eq:matrixsystemb}). 

If $(c, \tilde c)\in D(c, \tilde c)$, that $X$ is invertible if $(\psi(a_i), \nu(a_i), \beta(a_i))\neq (0,0,0)$ for all $1\leq i\leq k$ follows from
Proposition \ref{lem:invertsol}. For the converse, if $(\psi(a_i), \nu(a_i), \beta(a_i))= (0,0,0)$ for some $1\leq i\leq k$, it follows from the first equation in 
(\ref{eq:matrixsystemb}) that $\<X_sa_i, a_i\>=0$, hence $\<X a_i, a_i\>=\<X_aa_i, a_i\>=0$. On the other hand, the second 
equation in  (\ref{eq:matrixsystemb}) gives
$$0=\<(XA+A^tX^t)-\beta^t\beta-(c-\tilde c)\psi^t\psi)a_i, a_j\>=(\alpha_i-\alpha_j)\<Xa_i, a_j\>,\,\,\,i\neq j,
$$
hence $Xa_i=0$. The two last assertions follow from  Proposition \ref{p:ladm}.\vspace{1ex}\qed

Before concluding this section, we  compute  the pair $(\widetilde v,\widetilde h)$ (resp., the triple $(\widetilde v,\widetilde h, \widetilde V))$ associated to the $L$-transform of an isometric immersion $f\colon M^n(c)\rightarrow \mathbb Q_s^{n+p}(c)$ (resp.,  $f\colon M^n(c)\rightarrow \mathbb Q_s^{n+p}(\tilde c)$,  $c\neq \tilde c$).

\begin{proposition}\label{lemaproposto}\po
Let $f\colon M^n(c)\rightarrow \mathbb Q_s^{n+p}(c)$ be an isometric immersion 
and  $\xi_1, \ldots, \xi_p$ an orthonormal frame of $N_fM$  as in   Proposition \ref{thm:hiebetaij}. Let  $(v,h)$ be the pair associated to $f$. If  $\tilde f=\ral_{\varphi,\beta,L}(f)$  is an $L-$transform of $f$, then $\widetilde \xi_1,...,\widetilde \xi_p$ defined by  
\begin{equation}\label{eq:xitiu}
\widetilde \xi_r={\cal P}(\xi_r-\sum_{\ell=1}^p\epsilon_\ell\<\beta^r,L^{-1}\Omega^{-1}\beta^{\ell}\>\xi_{\ell})
\end{equation}
is an orthonormal  frame of $N_{\tilde f}M$ satisfying the conditions of Proposition \ref{thm:hiebetaij},
  and the pair $(\widetilde v,\widetilde h)$  associated to $\widetilde f$  is given by
\begin{equation}\label{eq:hbetatrans}
\widetilde v_j= v_j+\left<\beta^j,L^{-1}\Omega^{-1}\varphi\right> \ \mbox{and} \ \widetilde{h}_{ir}= h_{ir}+\left<\beta^r,L^{-1}\Omega^{-1}\gamma_i\right>,
\end{equation}
where $\gamma_1,...,\gamma_n$,  $\beta^1,...,\beta^p$ are given by  (\ref{eq:betagamma}).  In particular, $(\widetilde v,\widetilde h)$ is a new solution of (\ref{eq:i}).
\end{proposition}
\proof By Proposition \ref{prop:basic}, the
 metrics $\<\;,\;\>$  and
$\<\;,\;\>^\sim$ induced by $f$ and $\tilde{f}$, respectively, are related by 
$\<\;\;,\;\>^\sim=D^*\<\;\;,\;\>$, where $D=I-\Phi_{\Omega^{-1}\va}$ for $\Phi=\Phi(\omega, \beta)$.  It follows from  (\ref{eq:cconst}) that  
\be\label{eq:phiv}\Phi_{Lv}\partial_i=-A_{\beta(v)}\partial_i=-\<\beta^i,v\>X_i,\ee
where $X_i=v_i^{-1}\partial_i$. Thus
\begin{eqnarray}\label{eq:dembase}
D\partial_i=\widetilde v_iX_i, 
\end{eqnarray}
and the first of formulas (\ref{eq:hbetatrans}) follows. 
Using (\ref{alphasist}), (\ref{eq:sffs2}), (\ref{eq:omegabeta}), (\ref{eq:phiv}) and (\ref{eq:dembase}), we have for $1\leq i\leq n$ that 
$$
\begin{array}{lll}
  \tilde\alpha\left( \partial_i,{\partial_i}\right)&=& \P \left(\alpha\left(\partial_i,D\partial_i\right)+\beta(\Omega^{-1})^t\Phi\left(\partial_i\right)^tD\partial_i\right)\vspace{0.1cm}\\
&=& \P \left( \widetilde v_i\xi_i + \widetilde v_i\beta(\Omega^{-1})^t\Phi\left(\partial_i\right)^tX_i\right)\vspace{0.1cm}\\
&=&\tilde{v}_i\P\left(\xi_i-\sum_{\ell=1}^p\epsilon_\ell\left<\beta^i,L^{-1}\Omega^{-1}\beta^{\ell}\right>\xi_{\ell}\right)\vspace{0.1cm}\\
\end{array}
$$
 hence $\tilde{\xi}_1, \ldots, \tilde{\xi}_n$ are given by  (\ref{eq:xitiu}). On the other hand, from (\ref{alphasist}),   (\ref{eq:o1c}), (\ref{eq:ncon}) and equations $iv)$ and $(v)$ of system  ${\cal R}$ in Lemma \ref{prop:fundconst} we obtain  
$$
\widetilde \nabla^\bot_{\partial_i}\widetilde \xi_r=(h_{ir}+\left<\beta^r,L^{-1}\Omega^{-1}\gamma_i\right>)\P(\xi_i - \sum_{\ell=1}^p\epsilon_\ell\left<\beta^i,L^{-1}\Omega^{-1}\beta^\ell\right>\xi_\ell),\\
$$
which shows that $\widetilde \xi_1,...,\widetilde \xi_p$ 
is an orthonormal  frame of $N_{\tilde f}M$ satisfying the conditions of Proposition \ref{thm:hiebetaij}  and that $\tilde h_{ir}$ is given by (\ref{eq:hbetatrans}) for  $1\leq i\leq n$ and $1\leq r\leq p.$ \vspace{1ex}
\qed

Arguing in a similar way one obtains the following result in the case $c\neq \tilde c$.

\begin{proposition}\label{lemaproposto2}\po
Let $f\colon M^n(c)\rightarrow \mathbb Q_s^{n+p}(\tilde c)$, $c\neq \tilde c$,  be an isometric immersion of a simply connected Riemannian manifold satisfying the assumptions of Proposition \ref{teo:curvdif},  let $\xi_1, \ldots, \xi_p$ be a parallel orthonormal frame of $N_fM$ and let $(v,h,V)$ be  the   triple associated to $f$ and $\xi_1, \ldots, \xi_p$. If  $\tilde f=\ral_{\varphi,\beta,L}(f)$  is an $L-$transform of $f$, then $\widetilde \xi_1={\cal P}\xi_1,\ldots,\widetilde \xi_p={\cal P}\xi_p$ 
is a parallel orthonormal  frame of $N_{\tilde f}M$,  and the  triple $(\tilde v,\tilde h,\tilde V)$ associated to $\widetilde f$ 
and $\widetilde \xi_1,\ldots,\widetilde \xi_p$  is given by 
%\begin{equation}\label{eq:hbetatransntc}
$$
\tilde v_j= v_j-\left<B_j,\Omega^{-1}\varphi\right>, \ \tilde{h}_{ij}= 
h_{ij}-\left<B_j,\Omega^{-1}\gamma_i\right> 
\ \mbox{and} \ \tilde V_{ir}=V_{ir}+\epsilon_r\left<B_i,\Omega^{-1}\beta^r\right>,
$$
%\end{equation}
with $B_i=-(L^t)^{-1}\left(\sum_{r=1}^p \beta^{r}V_{ir}-(c-\tilde c)v_i\varphi\right)$.
In particular, $(\tilde v,\tilde h,\tilde V)$   is a solution of~(\ref{eq:ii}).
\end{proposition}

In the last result of this section  we study the inverse of an $L$-transformation. By Proposition 13 in \cite{dft}, if  $\tilde f=\ral_{\varphi,\beta,\Omega}(f)\colon\tilde M^n\rightarrow \mathbb Q_s^{n+p}(\tilde c)$ is the vectorial Ribaucour transform of $f\colon M^n\rightarrow \mathbb Q_s^{n+p}(\tilde c)$ determined by $(\varphi,\beta,\Omega)$ as in Definition \ref{ribvetgeral}, then  $$
\tilde \varphi=\Omega^{-1}\varphi, \ \tilde \beta={\cal P}\beta(\Omega^{-1})^t \ \mbox{and} \ \tilde \Omega=\Omega^{-1}
$$ 
satisfy the conditions of Definition \ref{ribvetgeral} for $\tilde f$ and  $f=\ral_{\tilde \varphi,\tilde \beta,\tilde \Omega}(\tilde f)$.  The next result states what else can be said if $\tilde f$ is an $L$-transform of $f$.

\begin{proposition}\label{prop:tltrans} \po
If $\tilde f=\ral_{\varphi,\beta,\Omega,L}(f)$ is an $L$-transform, then 
$f=\ral_{\tilde \varphi,\tilde \beta,\tilde \Omega}(\tilde f)$
 is an $L^t$-transform of $\tilde f$.
\end{proposition}
\proof 
It follows from (\ref{eq:cond3}) that $L(\Omega^t)^{-1}+\Omega^{-1}L^t=\Omega^{-1}\rho(\Omega^t)^{-1}$. 
On the other hand, by Proposition 3 in \cite{dft} we have 
$D\tilde \Phi_v=-\Phi_{\Omega^{-1}v}$
for all $v\in V$. Hence, using (\ref{eq:sffs}) and (\ref{eq:cconst}) we obtain 
$$
\begin{array}{l}
D\tilde \Phi(X)L^t +DA^{\tilde f}(X)\tilde \beta + (c-\tilde c)  DX\tilde\varphi^t\vspace{0.1cm}\\
\hspace{10ex}=-\Phi(X)\Omega^{-1}L^t+A(X)\beta(\Omega^{-1})^t+\Phi(X)\Omega^{-1}\beta^t\beta(\Omega^{-1})^t+(c-\tilde c)DX\varphi^t(\Omega^{-1})^t\vspace{0.1cm}\\
\hspace{10ex}=\Phi(X)L(\Omega^{-1})^t-\Phi(X)\Omega^{-1}\rho(\Omega^{-1})^t+A(X)\beta(\Omega^{-1})^t+\Phi(X)\Omega^{-1}\beta^t\beta(\Omega^{-1})^t\vspace{0.1cm}\\
\hspace{13ex}+(c-\tilde c)X\varphi^t(\Omega^{-1})^t-(c-\tilde c)\Phi(X)\Omega^{-1}\varphi\varphi^t(\Omega^{-1})^t=0
\end{array}
$$
Moreover,
$$
\tilde \Omega L^t+L\tilde \Omega^t= \Omega^{-1}\beta^t{\cal P}^t{\cal P} \beta(\Omega^t)^{-1}-(c-\tilde c)\Omega^{-1}\varphi\varphi^t(\Omega^{-1})^t=\tilde \beta^t\tilde \beta- (c-\tilde c)\tilde \varphi \tilde \varphi^t.\qed
$$

\section{A decomposition theorem for the
$L$-transformation}

To prove a  decomposition theorem for the
$L$-transformation we need the following.

\begin{lemma}\po\label{le:barphi} Under the assumptions of
Theorem \ref{c:viiinf}, the tensor
$$\bar{\Phi}_i(X)=\nabla_X\bar{\o}_i^t-A^j(X)\bar{\beta}_i+\tilde cX\bar\varphi_i^t$$ 
where $\bar{\o}_i=d\bar \varphi_i$, satisfies
\be\label{eq:barphi}D_j\bar{\Phi}_i(X)=\Phi_i(X)-\Phi_j(X)\Omega_{jj}^{-1}\Omega_{ji},\,\,\,1\leq i\neq j\leq 2,\ee
where $D_j=I-\Phi^j_{\Omega_{jj}^{-1}\va_j}.$
\end{lemma}
\proof 
%We first compute $\bar{\o}_i=d\bar{\va}_i$.
 We have
\begin{eqnarray}\label{eq:baromega2}\bar{\o}_i(X)&=&\o_i(X)-d{\Omega}_{ij}(X)\Omega_{jj}^{-1}\va_j+
{\Omega}_{ij}\Omega_{jj}^{-1}d\Omega_{jj}(X)\Omega_{jj}^{-1}\va_j
-{\Omega}_{ij}\Omega_{jj}^{-1}\o_j(X)\nonumber\\
&=&\o_i(X)-\o_i(\Phi^j(X)\Omega_{jj}^{-1}\va_j)
+{\Omega}_{ij}\Omega_{jj}^{-1}\o_j(\Phi^j(X)\Omega_{jj}^{-1}\va_j)
-{\Omega}_{ij}\Omega_{jj}^{-1}\o_j(X)\nonumber\\
&=&\o_i(D_jX)-{\Omega}_{ij}\Omega_{jj}^{-1}\o_j(D_jX).
\end{eqnarray}
Denoting by $\<\;,\;\>_j$  the metric induced by $f_j$, we obtain
\begin{eqnarray*}\<\bar{\o}_i^t(v_i),X\>_j\!\!\!
&=&\!\!\!\<v_i,\o_i(D_jX)-{\Omega}_{ij}\Omega_{jj}^{-1}\o_j(D_jX)\>\\
\!\!\!&=&\!\!\!\<D_j\o_i^t(v_i)-D_j\o_j^t(\Omega_{jj}^{-1})^t{\Omega}_{ij}^t(v_i),X\>\\
\!\!\!&=&\!\!\!\<\o_i^t(v_i)-\o_j^t(\Omega_{jj}^{-1})^t{\Omega}_{ij}^t(v_i)
,D_j^{-1}X\>_j.
\end{eqnarray*}
 Since $D_j^{-1}$ is symmetric with respect to $\<\;,\;\>_j$, it follows that
$D_j\bar{\o}_i^t=\o_i^t-\o_j^t(\Omega_{jj}^{-1})^t{\Omega}_{ij}^t$. 
Using (\ref{eq:lcivita}) we get
$$\begin{array}{l}
D_j(\nabla_X\bar{\o}_i^t)(v_i)=D_j{\nabla}^j_X\bar{\o}_i^t(v_i)-
D_j\bar{\o}_i^t(\nabla_X^{V_i}v_i)\vspace{1ex}\\
\hspace*{14.5ex}=\nabla_X\o_i^t(v_i)-\nabla_X\o_j^t(\Omega_{jj}^{-1})^t{\Omega}_{ij}^t(v_i)-
\Phi_j(X)\Omega_{jj}^{-1}\o_j\o_i^t(v_i)\vspace{1ex}\\\hspace*{16ex}
+\Phi_j(X)\Omega_{jj}^{-1}\o_j\o_j^t(\Omega_{jj}^{-1})^t{\Omega}_{ij}^t(v_i)
+\o_j^t(\Omega_{jj}^{-1})^t\Phi_j(X)^t\o_i^t(v_i)\vspace{1ex}\\\hspace*{16ex}-
\o_j^t(\Omega_{jj}^{-1})^t\Phi_j(X)^t\o_j^t(\Omega_{jj}^{-1})^t{\Omega}_{ij}^t(v_i)
+\o_j^t(\Omega_{jj}^{-1})^t{\Omega}_{ij}^t(\nabla_X^{V_i}v_i)-\o_i^t(\nabla_X^{V_i}v_i).
\end{array}
$$
Now,
$$\begin{array}{l}-\nabla^j_X\o_j^t(\Omega_{jj}^{-1})^t{\Omega}_{ij}^t(v_i)
+\o_j^t(\Omega_{jj}^{-1})^t{\Omega}_{ij}^t(\nabla^{V_i}_Xv_i)=
-(\nabla_X\o_j^t)(\Omega_{jj}^{-1})^t{\Omega}_{ij}^t(v_i)\vspace{1ex}\\
\hspace*{20ex}+\o_j^t(\Omega_{jj}^{-1})^t\Phi_j(X)^t\o_j^t(\Omega_{jj}^{-1})^t{\Omega}_{ij}^t(v_i)
- \o_j^t(\Omega_{jj}^{-1})^t\Phi_j(X)^t\o_i^t(v_i).
\end{array}
$$
Therefore
$$
D_j(\nabla_X\bar{\o}_i^t)=\nabla_X \o_i^t-(\nabla_X
\o_j^t)(\Omega_{jj}^{-1})^t{\Omega}_{ij}^t-
\Phi_j(X)\Omega_{jj}^{-1}\o_j\o_i^t+\Phi_j(X)\Omega_{jj}^{-1}\o_j\o_j^t(\Omega_{jj}^{-1})^t{\Omega}_{ij}^t.
$$
%$$\begin{array}{l}
%D_j(\nabla_X\bar{\o}_i^t)=\nabla_X \o_i^t-(\nabla_X
%\o_j^t)(\Omega_{jj}^{-1})^t{\Omega}_{ij}^t-
%\Phi_j(X)\Omega_{jj}^{-1}\o_j\o_i^t\vspace{1ex}\\
%\hspace*{15ex}+\Phi_j(X)\Omega_{jj}^{-1}\o_j\o_j^t(\Omega_{jj}^{-1})^t{\Omega}_{ij}^t.
%\end{array}
%$$
 On the other hand, from (\ref{eq:sffs}) we have
 $$D_j{A}^j(X)\bar{\beta}_i=A(X)\beta_i-
A(X)\beta_j(\Omega_{jj}^{-1})^t{\Omega}_{ij}^t+
\Phi_j(X)\Omega_{jj}^{-1}\beta_j^t\beta_i
-\Phi_j(X)\Omega_{jj}^{-1}\beta_j^t\beta_j(\Omega_{jj}^{-1})^t{\Omega}_{ij}^t.$$
%$$\begin{array}{l}D_j{A}^j(X)\bar{\beta}_i=A(X)\beta_i-
%A(X)\beta_j(\Omega_{jj}^{-1})^t{\Omega}_{ij}^t+
%\Phi_j(X)\Omega_{jj}^{-1}\beta_j^t\beta_i\vspace{1ex}\\
%\hspace*{15ex}
%-\Phi_j(X)\Omega_{jj}^{-1}\beta_j^t\beta_j(\Omega_{jj}^{-1})^t{\Omega}_{ij}^t.\end{array}$$
Also,
$$
\begin{array}{lll}
D_j(X)\bar \varphi_i^t
%&=&D_j(X)(\varphi_i-\Omega_{ij}\Omega_{jj}^{-1}\varphi_j)^t\vspace{0.1cm}\\
&=&(X-\Phi_j(X)\Omega_{jj}^{-1}\varphi_j)(\varphi_i^t-\varphi_j^t(\Omega_{jj}^{-1})^t\Omega_{ij}^t)\vspace{0.1cm}\\
&=&X\varphi_i^t-X\varphi_j^t(\Omega_{jj}^{-1})^t{\Omega_{ij}}^t-\Phi_j(X)\Omega_{jj}^{-1}\varphi_j\varphi_i^t
%\vspace{0.1cm}\\&&
+\Phi_j(X)\Omega_{jj}^{-1}\varphi_j\varphi_j^t(\Omega_{jj}^{-1})^t\Omega_{ij}^t.
\end{array}
$$
Thus 
\be\label{eq:j}\begin{array}{l}{D_j}\bar{\Phi}_i(X)=
\Phi_i(X)-\Phi_j(X)(\Omega_{jj}^{-1})^t{\Omega}_{ij}^t
-\Phi_j(X)\Omega_{jj}^{-1}\o_j\o_i^t
+\Phi_j(X)\Omega_{jj}^{-1}\o_j\o_j^t(\Omega_{jj}^{-1})^t{\Omega}_{ij}^t\vspace{1ex}\\\hspace*{13ex}
-\Phi_j(X)\Omega_{jj}^{-1}\beta_j^t\beta_i
+\Phi_j(X)\Omega_{jj}^{-1}\beta_j^t\beta_j(\Omega_{jj}^{-1})^t{\Omega}_{ij}^t
-\tilde c\,\Phi_j(X)\Omega_{jj}^{-1}\varphi_j\varphi_i^t\vspace{1ex}\\\hspace*{13ex}+\tilde c\,\Phi_j(X)\Omega_{jj}^{-1}\varphi_j\varphi_j^t(\Omega_{jj}^{-1})^t\Omega_{ij}^t.\\
\end{array}
\ee 
% \be\label{eq:j}\begin{array}{l}{D_j}\bar{\Phi}_i(X)=
%\Phi_i(X)-\Phi_j(X)(\Omega_{jj}^{-1})^t{\Omega}_{ij}^t
%-\Phi_j(X)\Omega_{jj}^{-1}\o_j\o_i^t
%\vspace{1ex}\\\hspace*{13ex}
%+\Phi_j(X)\Omega_{jj}^{-1}\o_j\o_j^t(\Omega_{jj}^{-1})^t{\Omega}_{ij}^t
%-\Phi_j(X)\Omega_{jj}^{-1}\beta_j^t\beta_i(v_i)\vspace{1ex}\\\hspace*{13ex}
%+\Phi_j(X)\Omega_{jj}^{-1}\beta_j^t\beta_j(\Omega_{jj}^{-1})^t{\Omega}_{ij}^t
%\vspace{1ex}\\\hspace*{13ex}
%-\tilde c\Phi_j(X)\Omega_{jj}^{-1}\varphi_j\varphi_i^t+\tilde c\Phi_j(X)\Omega_{jj}^{-1}\varphi_j\varphi_j^t(\Omega_{jj}^{-1})^t\Omega_{ij}^t.\\
%\end{array}
%\ee 
Using that ${\cal G}_j^t{\cal G}_j=\o_j\o_j^t+\beta_j^t\beta_j+\tilde c \varphi_j\varphi_j^t$ and
${\cal G}_j^t{\cal G}_i=\o_j\o_i^t+\beta_j^t\beta_i+\tilde c\varphi_j\varphi_i^t$ we obtain
$$\begin{array}{l}(\Omega_{jj}^{-1})^t{\Omega}_{ij}^t
+\Omega_{jj}^{-1}\o_j\o_i^t-\Omega_{jj}^{-1}\o_j\o_j^t(\Omega_{jj}^{-1})^t{\Omega}_{ij}^t+
\Omega_{jj}^{-1}\beta_j^t\beta_i
-\Omega_{jj}^{-1}\beta_j^t\beta_j(\Omega_{jj}^{-1})^t{\Omega}_{ij}^t\vspace{1ex}\\\hspace*{20ex}
+\tilde c\,\Omega_{jj}^{-1}\varphi_j\varphi_i^t-\tilde c\,\Omega_{jj}^{-1}\varphi_j\varphi_j^t(\Omega_{jj}^{-1})^t\Omega_{ij}^t\vspace{1ex}\\\hspace*{20ex}
=(\Omega_{jj}^{-1})^t{\Omega}_{ij}^t
+\Omega_{jj}^{-1}{\cal G}_j^t{\cal G}_i-\Omega_{jj}^{-1}{\cal G}_j^t{\cal G}_j(\Omega_{jj}^{-1})^t{\Omega}_{ij}^t\vspace{1ex}\\\hspace*{20ex}
=(\Omega_{jj}^{-1})^t{\Omega}_{ij}^t+\Omega_{jj}^{-1}{\Omega}_{ij}^t
+\Omega_{jj}^{-1}{\Omega}_{ji}-(\Omega_{jj}^{-1})^t{\Omega}_{ij}^t
-\Omega_{jj}^{-1}{\Omega}_{ij}^t\vspace{1ex}\\\hspace*{20ex}=\Omega_{jj}^{-1}{\Omega}_{ji}.
\end{array}$$
Substituting  this into (\ref{eq:j}) yields (\ref{eq:barphi}).\qed

\begin{theorem}\po\label{thm:permconst} Let $f\colon\,M^n(c)\to \Q^{n+p}(\tilde c)$, $(c, \tilde c)\in D(c, \tilde c)$, be an isometric immersion
 and let ${\cal R}_{\va,\beta,\Omega,L}(f)\colon\, \tilde{M^n}\to \Q^{n+p}(\tilde c)$ be
 an $L$-transform of $f$. Assume that $L$ decomposes as $L=L_1\oplus L_2$ with respect to an orthogonal decomposition $V=V_1\oplus V_2$.
%$\tilde c\geq 0$ or that $\tilde c<0$ and $c\in (-\infty,\tilde c)\cup (0,\infty)$.
If 
%${\cal R}_{\va,\beta,\Omega,L}(f)\colon\, \tilde{M^n}\to \Q^{n+p}(\tilde c)$ be
 %an $L$-transform of $f$,  and let  
 $(\varphi_j,  \beta_j,  \Omega_{ij})$ and $(\bar{\va}_i,\bar{\beta}_i,\bar{\Omega}_{ii})$ are given as in (\ref{vabo2nf})  and   
%${\cal R}_{\va,\beta,\Omega}(\va_i,\beta_i,\Omega_{ii}):=(\bar{\va}_i,\bar{\beta}_i,\bar{\Omega}_{ii})$ 
 (\ref{eq:bar2nf}), respectively, then $(\varphi_j,\beta_j,\Omega_{jj},L_j)$ defines an $L_j$-transform of $f$ for $1\leq j\leq 2$,
$(\bar{\va}_i,\bar{\beta}_i,\bar{\Omega}_{ii},L_i)$ an $L_i$-transform of $f_j$ for $1\leq i\neq j\leq 2$, and
%\be\label{eq:perm}
$$ {\cal R}_{\va,\beta,\Omega,L}(f) ={\cal
R}_{\bar{\va}_i,\bar{\beta}_i,\bar{\Omega}_{ii},L_i} ({\cal
R}_{\va_j,\beta_j,\Omega_{jj},L_j}(f)). 
$$
%\ee
The same conclusion  holds if $(c, \tilde c)\not\in D(c, \tilde c)$
%$\tilde c<0$ and $c\in [\tilde c,0]$,
 as long as  $\Omega_{jj}$ is invertible for $1\leq j\leq 2$.
\end{theorem}
\proof Since ${\cal R}_{\va,\beta,\Omega,L}(f)$ is
an $L$-transform of $f$, equations (\ref{eq:cconst}) and (\ref{eq:cond3}) hold. 
%$$\Phi(X)L+A(X)\beta+(c-\tilde c)X\va^t=0$$
%and $$\Omega L +L^t\Omega^t-\rho=0,$$ with $\rho$ given by
%(\ref{eq:rho}). 
By the assumption that $L=L_1\oplus L_2$, these equations are
equivalent to
\be\label{eq:dec1}\Theta_{j}(X):=\Phi_j(X)L_j+A(X)\beta_j+(c-\tilde c)X\va_j^t=0\ee
%$$\Omega_{ii}L_i+L_i^t\Omega_{ii}^t-(\beta_i^t\beta_i-(c-\tilde c)\va_i\va_i^t)=0,$$
and
\be\label{eq:dec2}
\Omega_{ij}L_j+L_i^t\Omega_{ji}^t-\beta_i^t\beta_j+(c-\tilde c)\va_i\va_j^t=0,
\ee
$1\leq
i,j\leq 2.$ Therefore, in order to prove that $(\varphi_j,\beta_j,\Omega_{jj},L_j)$ defines an $L_j$-transform of $f$ for $1\leq j\leq 2$ it remains to show that $\Omega_{jj}$ is invertible for $1\leq j\leq 2$. Assume otherwise, say,  that $\ker \Omega_{11}\neq \{0\}$.  Let $\alpha_1, \ldots, \alpha_p$ (resp., $\gamma_1, \ldots, \gamma_q$) be the real (resp., complex) eigenvalues of $L$,  and let  
%$A\in M_m(\mathbb{R})$ be a nonderogatory matrix. For each $i=1,...,p$ (resp., $j=1,...,q$),  let
  $a_{i1},a_{i2},...,a_{in_i}$ 
(resp.,  $w_{j1},w_{j2},...,w_{jm_j}$) be the generalized eigenvectors of $L$ associated to $\alpha_i$, $1\leq i\leq p$ 
%Jordan block $J_{n_i}(\alpha_i)$ (resp., $J_{m_j}(\alpha_j)$) of $A$ correspondent to the real eigenvalue (resp., complex) $\alpha_i$ 
(resp., $\gamma_j$, $1\leq j\leq q$ ). Since $(c, \tilde c)\in D(c, \tilde c)$, by Proposition \ref{lem:invertsol} either there exists 
$1\leq i\leq p$ such that $a_{i1}\in \ker \Omega_{11}\subset S_1$, where $S_1$ is defined as in (\ref{eq:s}) for $(\psi_1=\varphi^t|_{V_1}, \nu_1=\omega^t|_{V_1}, \beta_1=\beta|_{V_1})$, or there exists $1\leq j\leq q$ such that $a_{j1}\in (\ker \Omega_{11})^c\subset S_1^c$. We argue for  the first possibility, the argument for the second one being similar. 

Thus we have $La_{i1}=\alpha_i a_{i1}$, with $a_{i1}\in S\cap V_1$. Since $\spa\{a_{i1}, \ldots, a_{in_i}\}\subset V_1$, we have
\be\label{eq:v1}0=\<\Omega_{11}a_{i1}, a_{i\ell}\>=\<\Omega a_{i1}, a_{i\ell}\>,\,\,\,1\leq \ell \leq n_i.\ee
%for all $1\leq \ell \leq n_i$. 

 On the other hand,  from $\Omega+\Omega^t={\cal G}^t{\cal G}$ and (\ref{eq:cond3}) that $\Omega_sa_{i1}=0$ and $L^t\Omega_aa_{i1}=\alpha_i\Omega_aa_{i1}$. Thus, for all real eigenvalues $\alpha_j$, $1\leq j\leq p$, with $\alpha_j\neq \alpha_i$, and all complex eigenvalues $\gamma_j$, $1\leq j\leq q$, we have 
%\begin{eqnarray}\label{eq:i1}
%\alpha_i\<\Omega_aa_{i1}, a_{i\ell}\>&=&\<\Omega_aa_{i1}, La_{i\ell}-a_{i,\ell-1}\>\nonumber\\
%&=&\<L^t\Omega_aa_{i1}, a_{i\ell}\>-\<\Omega_aa_{i1}, a_{i,\ell-1}\>\nonumber\\
%&=&\alpha_i\<\Omega_aa_{i1}, a_{i\ell}\>-\<\Omega_aa_{i1}, a_{i,\ell-1}\>,
%\end{eqnarray}
\begin{eqnarray}\label{eq:i2}\alpha_j\<\Omega_aa_{i1}, a_{jt}\>&=&\<\Omega_aa_{i1}, La_{jt}-a_{j,t-1}\>\nonumber\\
&=&\<L^t\Omega_aa_{i1}, a_{jt}\>-\<\Omega_aa_{i1}, a_{j,t-1}\>\nonumber\\
&=&\alpha_i\<\Omega_aa_{i1}, a_{jt}\>-\<\Omega_aa_{i1}, a_{j,t-1}\>,\,\,\,a_{j0}=0,
\end{eqnarray}
\begin{eqnarray}\label{eq:i3}\bar\gamma_j\<\Omega_aa_{i1}, w_{jk}\>&=&\<\Omega_aa_{i1}, Lw_{jk}-w_{j,k-1}\>\nonumber\\
&=&\<L^t\Omega_aa_{i1}, w_{jk}\>-\<\Omega_aa_{i1}, w_{j,k-1}\>\nonumber\\
&=&\alpha_i\<\Omega_aa_{i1}, w_{jk}\>-\<\Omega_aa_{i1}, w_{j,k-1}\>, \,\,\, w_{j0}=0.
\end{eqnarray}
It follows from (\ref{eq:v1}), (\ref{eq:i2}),   (\ref{eq:i3}) and $\Omega_s a_{i1}=0$  that $\Omega a_{i1}=0$, a contradiction.

%Moreover, $(\varphi^t,\omega^t,\beta)$ is $L$-admissible. So if either $\tilde c\geq 0$ or $\tilde c<0$ and $c\in (-\infty,\tilde c)\cup (0,\infty)$, it follows from Theorem \ref{p:ladm} that $E_\alpha\cap S=\{0\}$ (respectively, $E_\alpha\cap S^c=\{0\}$) if $\alpha$ is a real eigenvalue of $L$ with odd algebraic multiplicity (respectively, complex eigenvalue). It means that $E_{\alpha} \cap S_j=\{0\}$ (respectively, $E_{\alpha} \cap S_j^c=\{0\}$) if real (respectively, complex) number $\alpha$ is eigenvalue of $L_j$, where 
%$$
%S_j=\left\{\begin{array}{l} \ker \omega^t_j\cap \ker \beta_j, \ \mbox{if} \ c=\tilde c=0,\\
%            \ker \varphi_j^t\cap \ker \omega^t_j\cap \ker \beta_j, \ \mbox{if} \ (c,\tilde c)\neq (0,0).\end{array}\right.\\
%$$
%and $S_j^c$ the complexified of $S_j$. Thus, again by Theorem \ref{p:ladm} and Affirmation 1, we have that $(\varphi_j^t,\omega_j^t,\beta_j^t)$ is $L_j$-admissible, namely, $\Omega_{jj}$ is invertible. Therefore, by Theorem \ref{c:viiinf} and Definition \ref{df:ribconst}, $f_j$ is an $L_j$-transform of $f$. For $\tilde c<0$ and $c\in[\tilde c,0]$, $f_j$ is an $L_j$-transform of $f$ since $\Omega_{jj}$ is invertible.

%Now, by (\ref{eq:barphi}) we have
%$$D_j\bar{\Phi}_i(X)=\Phi_i(X)-\Phi_j(X)\Omega_{jj}^{-1}\Omega_{ji}.$$

We now show that $(\bar{\va}_i,\bar{\beta}_i,\bar{\Omega}_{ii},L_i)$ defines an $L_i$-transform of $f_j$ for $1\leq i\neq j\leq 2$. We must prove that
$$\bar{\Theta}_i(X):=\bar{\Phi}_i(X)L_i+A^j(X)\bar{\beta}_i+(c-\tilde c)X\bar{\va}_i^t=0=\bar{\Omega}_{ii}L_i+L_i^t\bar{\Omega}_{ii}^t-
(\bar{\beta}_i^t\bar{\beta}_i-(c-\tilde c)\bar{\va}_i\bar{\va}_i^t).$$
%$\bar{\Theta}_i(X):=\bar{\Phi}_i(X)L_i+A^j(X)\bar{\beta}_i+(c-\tilde c)X\bar{\va}_i^t=0$
%and
%$\bar{\Omega}_{ii}L_i+L_i^t\bar{\Omega}_{ii}^t-
%(\bar{\beta}_i^t\bar{\beta}_i-(c-\tilde c)\bar{\va}_i\bar{\va}_i^t)=0.$
Using (\ref{eq:sffs}),  (\ref{eq:barphi}), (\ref{eq:dec1}) and (\ref{eq:dec2}) we obtain
$$\begin{array}{l}
D_j\bar{\Theta}_i(X)=
%D_j\bar{\Phi}_i(X)L_i+D_jA^j(X)\bar{\beta}_i+(c-\tilde c)D_jX\bar{\va}_i^t\vspace{1ex}\\
%\hspace*{9.8ex}=
\Phi_i(X)L_i-\Phi_j(X)\Omega_{jj}^{-1}\Omega_{ji}L_i
+D_jA^j(X)\P_j(\beta_i-\beta_j(\Omega_{jj}^{-1})^t\Omega_{ij}^t)
\vspace{1ex}\\
\hspace*{11.5ex}+(c-\tilde c)D_j(X)(\va_i-\Omega_{ij}\Omega_{jj}^{-1}\va_j)^t\vspace{1ex}\\
\hspace*{9.7ex}=\Phi_i(X)L_i-\Phi_j(X)\Omega_{jj}^{-1}\Omega_{ji}L_i
+(A(X)+\Phi_j(X)\Omega_{jj}^{-1}\beta_j^t)(\beta_i-\beta_j(\Omega_{jj}^{-1})^t\Omega_{ij}^t)\vspace{1ex}\\
\hspace*{11.5ex}
%+\Phi_j(X)\Omega_{jj}^{-1}\beta_j^t(\beta_i-\beta_j(\Omega_{jj}^{-1})^t\Omega_{ij}^t)
%\vspace{1ex}\\\hspace*{11.5ex}
+(c-\tilde c)(X-\Phi_j(X)\Omega_{jj}^{-1}\va_j)(\va_i^t-\va_j^t(\Omega_{jj}^{-1})^t\Omega_{ij}^t)\vspace{1ex}\\
\hspace*{9.7ex}=\Phi_i(X)L_i+A(X)\beta_i+(c-\tilde c)X\va_i^t
-\Phi_j(X)\Omega_{jj}^{-1}(\Omega_{ji}L_i-\beta_j^t\beta_i+(c-\tilde c)\va_j\va_i^t)\vspace{1ex}\\
\hspace*{11.5ex}
+(\Phi_j(X)\Omega_{jj}^{-1}(-\beta_j^t\beta_j+(c-\tilde c)\va_j\va_j^t)-A(X)\beta_j-(c-\tilde c)X\va_j^t)
(\Omega_{jj}^{-1})^t\Omega_{ij}^t=0.
%\vspace{1ex}\\
%\hspace*{9.7ex}=\Phi_j(X)\Omega_{jj}^{-1}L_j^t\Omega_{ij}^t
%-\Phi_j(X)\Omega_{jj}^{-1}L_j^t\Omega_{jj}^t(\Omega_{jj}^{-1})^t\Omega_{ij}^t=0.
\end{array}$$
On the other hand,
$$\begin{array}{l}
\bar{\Omega}_{ii}L_i+L_i^t\bar{\Omega}_{ii}^t-
(\bar{\beta}_i^t\bar{\beta}_i-(c-\tilde c)\bar{\va}_i\bar{\va}_i^t)\vspace{1ex}\\
\hspace*{7.5ex}=(\Omega_{ii}-\Omega_{ij}\Omega_{jj}^{-1}\Omega_{ji})L_i
+L_i^t(\Omega_{ii}^t-\Omega_{ji}^t(\Omega_{jj}^{-1})^t\Omega_{ij}^t)
\vspace{1ex}\\
\hspace*{7.5ex} -(\beta_i^t-\Omega_{ij}\Omega^{-1}_{jj}\beta_j^t)
(\beta_i-\beta_j(\Omega^{-1}_{jj})^t\Omega_{ij}^t)
+(c-\tilde c)(\va_i-\Omega_{ij}\Omega_{jj}^{-1}\va_j)(\va_i^t-\va_j^t(\Omega_{jj}^{-1})^t\Omega_{ij}^t)
\vspace{1ex}\\
\hspace*{7.5ex}=\Omega_{ii}L_i+L_i^t\Omega_{ii}^t-\beta_i^t\beta_i+(c-\tilde c)\va_i\va_i^t
-\Omega_{ij}\Omega_{jj}^{-1}(\Omega_{ji}L_i-\beta_j^t\beta_i+(c-\tilde c)\va_j\va_i^t)\vspace{1ex}\\
\hspace*{7.5ex}-(L_i^t\Omega_{ji}^t-\beta_i^t\beta_j+(c-\tilde c)\va_i\va_j^t)
(\Omega_{jj}^{-1})^t\Omega_{ij}^t-\Omega_{ij}\Omega_{jj}^{-1}(\beta_j^t\beta_j-(c-\tilde c)\va_j\va_j^t)
(\Omega_{jj}^{-1})^t\Omega_{ij}^t\vspace{1ex}\\
\hspace*{7.5ex}=\Omega_{ij}\Omega_{jj}^{-1}L_j^t\Omega_{ij}^t+
\Omega_{ij}L_j(\Omega_{jj}^{-1})^t\Omega_{ij}^t-\Omega_{ij}\Omega_{jj}^{-1}(\Omega_{jj}L_j+L_j^t\Omega_{jj}^t)
(\Omega_{jj}^{-1})^t\Omega_{ij}^t=0.\vspace{1ex}\qed
\end{array}$$

\begin{remarks}\po\label{re:ldecomp}{\em $(i)$ In Theorem \ref{thm:permconst}, if  $L$ decomposes as $L=L_1\oplus L_2$ with respect to any direct sum  decomposition $V=V_1\oplus V_2$, then one can always assume the decomposition to be orthogonal after a suitable change of the inner product on $V$. Namely, define $W_i=V_j^\perp$, $1\leq i\neq j\leq 2$, and let $B\in V^*\otimes V$ be an invertible endomorphism such that $BW_i=V_i$, $1\leq i\leq 2$. Let 
$\<\cdot , \cdot\>^\sim$ be the inner product on $V$ given by 
 %\be\label{eq:innerp}
 $\<\cdot, \cdot\>^\sim = \<\cdot, B\cdot\>$. Then
 $W_1$ and $W_2$ are orthogonal with respect to $\<\cdot , \cdot\>^\sim$, because $W_1$ and $BW_2=V_2$ are orthogonal with respect to $\<\cdot , \cdot\>$.
% $$\<W_1 , W_2\>^\sim=\<W_1, BW_2\>=\<W_1, V_2\>=0$$
 Moreover,  $\hat L=B^{-1}LB$ decomposes as $\hat L=\hat L_1\oplus \hat L_2$ with respect to the  decomposition $V=W_1\oplus W_2$, for
 $$\hat LW_i=B^{-1}LBW_i=B^{-1}LV_i\subset B^{-1}V_i=W_i, \;\;1\leq i\leq 2.$$
 $(ii)$  In particular, let   $\tilde f={\cal R}_{\varphi,\beta,\Omega,L}(f)$
 be an $L$-transform of $f$ such that $L$ is   diagonalizable, with $v_1, \ldots, v_k$ as a basis of eigenvectors. If $e_1, \ldots, e_k\in V$ are such that 
 $\<e_i, v_j\>=\delta_i^j$, $1\leq i,j\leq k$, and we define $B\in V^*\otimes V$ by $B e_i=v_i$, $1\leq i\leq k$,  then  $e_1, \ldots, e_k$ is a basis
 of eigenvectors of $\hat L=B^{-1}LB$ that is orthonormal with respect to $\<\cdot , \cdot\>^\sim= \<\cdot, B\cdot\>$. Thus, if ${\cal R}_{\varphi,\beta,\Omega,L}(f)$ is an $L$-transform of $f$ such that  $L$ is diagonalizable, we may always assume that $L$ is symmetric after an appropriate  change of the inner product on $V$. It follows from  Theorem \ref{thm:permconst} that an $L$-transformation can be obtained by an iteration of scalar $L_i$-transformations, $1\leq i\leq k=\dim V$, if and only if $L$ is diagonalizable.}
 \end{remarks}

%\vspace{3ex}
%\noindent{\bf\large \S 3.2 The Bianchi $L$-cube.} \vspace{3ex}
\subsection{The Bianchi $L$-cube}

Given $L_1\neq L_2\in \mathbb R$, we say that a Bianchi quadrilateral  $\{f,f_1,f_2,f_{12}\}$  is a \textit{Bianchi $(L_1,L_2)$-quadrilateral} if the metric induced by $f$ has constant  curvature $c$, $f_i$ is an $L_i$-transform of $f$,  $1\leq i\leq 2$, and $f_{12}$ is an $L_2$-transform of $f_1$ and an $L_1$-transform of~$f_2$.

We will need the following fact proved in  \cite{dt4}, which can also be derived from Theorem~\ref{thm:permconst}.
%(see Theorem 18 in \cite{dt4})

%\begin{proposition}\label{prop:perm}
%Let $f\colon M^n(c)\rightarrow \mathbb{Q}^{n+p}(\tilde c)$ be an isometric immersion. Assume that either $\tilde c\geq 0$ or that  $\tilde c<0$ and $c\in (-\infty,\tilde c)\cup (0,\infty)$. Given $L_1\neq L_2\in \R$, let 
%$$
%f_r={\cal R}_{\varphi_r,\beta_r}(f)\colon M^n(c) \rightarrow \mathbb{Q}^{n+p}(\tilde c)
%$$
%be an $L_r$-transform of $f$,  $1\leq r\leq 2$. If $[A_{\beta_1}, A_{\beta_2}]=0$, then there exists exactly one  
%isometric immersion $\tilde f \colon \tilde M^n(c) \rightarrow \mathbb{Q}^{n+p}(\tilde c)$ such that $\{f,f_1,f_2,\tilde f\}$ is a Bianchi $(L_1,L_2)$-quadrilateral. 
%\end{proposition}

\begin{proposition}\label{prop:perm}\po
Let $f\colon M^n(c)\rightarrow \mathbb{Q}^{n+p}(\tilde c)$, $(c, \tilde c)\in D(c, \tilde c)$, be an isometric immersion.
% Assume that either $\tilde c\geq 0$ or that  $\tilde c<0$ and $c\in (-\infty,\tilde c)\cup (0,\infty)$. 
If $f_i={\cal R}_{\varphi_i,\beta_i}(f)\colon M_i^n(c) \rightarrow \mathbb{Q}^{n+p}(\tilde c)$
is an $L_i$-transform of $f$,  $1\leq i\leq 2$, with  $L_1\neq L_2$ and $[A_{\beta_1}, A_{\beta_2}]=0$, then there exists exactly one  
isometric immersion $\tilde f \colon \tilde M^n(c) \rightarrow \mathbb{Q}^{n+p}(\tilde c)$ such that $\{f,f_1,f_2,\tilde f\}$ is a Bianchi $(L_1,L_2)$-quadrilateral. 
\end{proposition}

%For each $r\in \{1,...,k\}$, consider the set of multi-indices
%$$
%\Lambda_r=\{\alpha_r=\{i_1,...,i_r\}\subset \{1,...,k\}\,:\, \alpha_r\,\mbox{with $r$ distint elements}\}.
%$$
Given $L_1,...,L_k\in \mathbb R$, with $L_i\neq L_j$ for all $1\leq i\neq j\leq k$, we say that a Bianchi cube
$({\cal{C}}_0, . . . ,{\cal{C}}_k)$ is a \textit{Bianchi $(L_1,...,L_k)$-cube} if, for all $1\leq s\leq k-1$,
\begin{itemize}
	\item[(i)] Each $f_{\alpha_{s+1}}\in {\cal C}_{s+1}$ with  $\alpha_{s+1}=\alpha_{s}\cup \{i_j\}$  is an $L_{i_j}$-transform of $f_{\alpha_{s}}\in {\cal C}_s$.
		\item[(ii)]   $\{f_{\alpha_{s-1}},f_{\alpha_{s-1}\cup \{i_l\}},f_{\alpha_{s-1}\cup \{i_j\}},f_{\alpha_{s+1}}\}$ is a Bianchi $(L_{i_l},L_{i_j})$-quadrilateral when $\alpha_{s+1}=\alpha_{s-1}\cup \{i_l,i_j\}$.
%	\item[(ii)] If $\alpha_{s+1}=\alpha_{s-1}\cup \{i_l,i_j\}$ then  $\{f_{\alpha_{s-1}},f_{\alpha_{s-1}\cup \{i_l\}},f_{\alpha_{s-1}\cup \{i_j\}},f_{\alpha_{s+1}}\}$ is a Bianchi $(L_{i_l},L_{i_j})$-quadrilateral.
\end{itemize}

\begin{theorem}\label{lcubo}\po Let $f\colon M^n(c)\rightarrow \mathbb{Q}^{n+p}(\tilde c)$, $(c, \tilde c)\in D(c, \tilde c)$, be an isometric immersion.  If 
%Assume either that $\tilde c\geq 0$ or that  $\tilde c<0$ and $c\in (-\infty,\tilde c)\cup (0,\infty)$. 
$
f_i={\cal R}_{\varphi_i,\beta_i}(f)\colon M_i^n(c) \rightarrow  \mathbb{Q}^{n+p}(\tilde c)
$, $1\leq i\leq k$, are  $L_i$-transforms of $f$ such that $L_i\neq L_j$ and $[A_{\beta_i},A_{\beta_j}]=0$ for all $1\leq  i\neq j\leq k$,  then there exists a  Bianchi $(L_1,...,L_k)$-cube $({\cal C}_0, . . . , {\cal C}_k)$ such that
${\cal C}_0 = \{f\}$ and ${\cal C}_1 = \{f_1, . . . , f_k\}$, which is unique if  no $f_i$ belongs to the associated family determined by $\{f_j, f_l\}$ for all  $1\leq i\neq j\neq l\neq i\leq k$.
\end{theorem}
%\begin{theorem}\label{lcubo} Let $f\colon M^n(c)\rightarrow \mathbb{Q}^{n+p}(\tilde c)$, $(c, \tilde c)\in D(c, \tilde c)$, be an isometric immersion. 
%Assume either that $\tilde c\geq 0$ or that  $\tilde c<0$ and $c\in (-\infty,\tilde c)\cup (0,\infty)$. 
%For each $1\leq i\leq k$, let
%$
%f_i={\cal R}_{\varphi_i,\beta_i}(f)\colon M^n(c) \rightarrow  \mathbb{Q}^{n+p}(\tilde c)
%$
%be an $L_i$-transform of $f$, with $L_i\neq L_j$ if $i\neq j$. 
%If $[A_{\beta_i},A_{\beta_j}]=0, \ i\neq j,$ and no $f_i$ belongs to the associated family determined by $\{f_j, f_l\}$ for all  $1\leq i\neq j\neq l\neq i\leq k$, then there is a unique Bianchi $(L_1,...,L_k)-$cube $({\cal C}_0, . . . , {\cal C}_k)$ such that
%${\cal C}_0 = \{f\}$ and ${\cal C}_1 = \{f_1, . . . , f_k\}$.
%\end{theorem}
\proof We first prove existence. Set $F=i\circ f$
% where $i\colon \Q^{n+p}(\tilde c)\to \R^{n+p+1}_{\epsilon_0}$ is the inclusion, 
and  ${\cal G}_i=F_*\nabla \varphi_i+i_*\beta_i+\tilde c\varphi_iF$. Since $f_i$ is an $L_i$-transform of $f$, for each $1\leq i\leq k$ the pair $(\varphi_i, \beta_i)$ satisfies
\be\label{vibi}
\alpha(X, \nabla \varphi_i)+\nabla_X^{\perp}\beta_i=0
\ee
for all $X\in \mathfrak{X}(M)$ and  the tensor 
$\Phi_i=\hess \varphi_i-A_{\beta_i}-\tilde c\varphi_iI$ satisfies
\be\label{eq:liphii}
L_i\Phi_i+A_{\beta_i}+(c-\tilde c)\varphi_i I=0.
\ee
 Moreover, $L_i|{\cal G}_i|^2=|\beta_i|^2-(c-\tilde c)\varphi^2_i$.
Define $\varphi\colon M^n(c)\to \R^k$,  $\beta\in \Gamma((\R^k)^*\otimes N_fM)$ and $\Omega\in \Gamma((\R^k)^*\otimes \R^k)$ by  
$$\varphi=(\varphi_1, \ldots, \varphi_k),\,\,\,\,\,\beta=\sum_{i=1}^k e^i\otimes \beta_i \,\,\,\mbox{and}\,\,\,\,\Omega=\sum \Omega_{ij} e^j\otimes e_i,$$
where $e_1, \ldots, e_k$ is the canonical basis of $\R^k$,  $e^1, \ldots, e^k$ is its dual basis, 
\be\label{eq:omegaij}\Omega_{ij}=\frac{1}{L_j-L_i}\left(\rho_{ij}-L_i\left<{\cal G}_i,{\cal G}_j\right>\right),\,\,\,i\neq j,\,\,\,\mbox{and}\,\,\,\Omega_{ii}=\frac{1}{2}|{\cal G}_i|^2,\,\,\,1\leq i\leq k,
\ee
where $\rho_{ij}=\<\beta_i, \beta_j\>-(c-\tilde c)\varphi_i\varphi_j$.

We claim that $(\varphi, \beta, \Omega)$ defines an $L$-transformation of $f$, where $L\in (\R^k)^*\otimes \R^k$ is given by
$L=\sum_{i=1}^kL_ie^i\otimes e_i$.  First, that $(\varphi, \beta)$ satisfies (\ref{eq:alphao}) follows from (\ref{vibi}).
The comutativity relations (\ref{eq:phicom}) follow from the assumption that  $[A_{\beta_i},A_{\beta_j}]=0, \ i\neq j,$
because $\Phi=\Phi(\varphi, \beta)$ 
%$$\Phi_vX=(\nabla_X \omega^t)v-A^f_{\beta(v)}X+\tilde c\varphi^t(v)X,$$
%with $\omega=d\varphi$, 
satisfies (\ref{eq:cconst}) in view of  (\ref{eq:liphii}).

Equation (\ref{eq:o1b}) is equivalent to 
\be \label{eq:o1d}
X(\Omega_{ij})=\<\nabla \varphi_i, \Phi_jX\>
\ee
for all $X\in \mathfrak{X}(M)$ and $1\leq i,j\leq k$. Now, 
\begin{eqnarray}\label{eq:xrhoij}
X(\rho_{ij})&=&\<\nabla^\perp_X\beta_i, \beta_j\>+\<\beta_i, \nabla^\perp_X\beta_j\> -(c-\tilde c)\<\nabla \varphi_i, X\>\varphi_j-(c-\tilde c)\<\nabla \varphi_j, X\>\varphi_i\nonumber\\
&=&-\<(A^f_{\beta_j}+(c-\tilde c)\varphi_jI)X, \nabla \varphi_i\>-\<(A^f_{\beta_i}+(c-\tilde c)\varphi_iI)X, \nabla \varphi_j\>
\end{eqnarray}
and
\begin{eqnarray}\label{eq:xgigj}X(\<{\cal G}_i,{\cal G}_j\>)&=&\<{{\cal G}_i}_*X,{\cal G}_j\>+\<{\cal G}_i,{{\cal G}_j}_*X\>\nonumber\\
%&=&\<F_*\Phi_iX, {\cal G}_j\>+\<{\cal G}_i, F_*\Phi_jX\>\nonumber\\
&=&\<\Phi_iX, \nabla \varphi_j\>+\<\nabla \varphi_i, \Phi_jX\>.
\end{eqnarray}
Since 
$$X(\Omega_{ij})=\frac{1}{(L_j-L_i)}(X(\rho_{ij})-L_iX(\<{\cal G}_i,{\cal G}_j\>)),$$
 (\ref{eq:o1d}) follows from (\ref{eq:liphii}), (\ref{eq:xrhoij}) and (\ref{eq:xgigj}). Finally,   equation  (\ref{eq:omegaij}) implies  (\ref{eq:o2b}) and  (\ref{eq:cond3}).

  Since  $\Omega_{ii}$ is nowhere vanishing for $1\leq i\leq k$, if $c\neq \tilde c$ (respectively,  $c=\tilde c$) then   $\varphi_i$ and  $\beta_i$  (resp., $\beta_i$) do not vanish simultaneously (resp., does not vanish) at any point.  In particular,  no $e_i$, $1\leq i\leq k$, belongs to the subspace $S$ defined by (\ref{eq:s}).   It follows from Proposition \ref{lem:invertsol} that  $(\varphi^t,\omega^t,\beta)$ is $L$-admissible, that is,     $\Omega$ is invertible. Thus $\tilde f=\ral_{\varphi,\beta,\Omega}(f)$ is an $L$-transform of $f$.

 For any multi-index $\alpha_r=\{i_1,...,i_r\}\in \Lambda_r$, define
\begin{equation}\label{eq:varbetaalp}
\varphi^{\alpha_r}=\sum_{j=1}^r\varphi_{i_j}e_{i_j}, \ \beta^{\alpha_r}=\sum_{j=1}^r e^{i_j}\otimes \beta_{i_j}, \, \Omega^{\alpha_r}=\sum_{\ell, j=1}^r\Omega_{i_\ell i_j}e^{i_j}\otimes e_{i_\ell}, \, 
L^{\alpha_r}=\sum_{\ell, j=1}^rL_{i_\ell i_j}e^{i_\ell}\otimes e_{i_j}.
\end{equation}
By Theorem \ref{thm:permconst}, for each $\alpha_r\in \Lambda_r$ the triple $(\varphi^{\alpha_r}, \beta^{\alpha_r}, \Omega^{\alpha_r})$ defines an 
$L^{\alpha_r}$-transform of~$f$. Let  ${\cal C}_r$ be the family of ${k \choose r}$ elements formed by the $L^{\alpha_r}$-vectorial Ribaucour transforms $f_{\alpha_r}=\ral_{\varphi^{\alpha_{r}},\beta^{\alpha_{r}},\Omega^{\alpha_{r}},L^{\alpha_{r}}}(f)$ of $f$.
It also follows from Theorem \ref{thm:permconst} that conditions $(i)$ and $(ii)$ in the definition of a Bianchi $(L_1,...,L_k)$-cube are satisfied by the $(k+1)$-tuple $({\cal{C}}_0, . . . ,{\cal{C}}_k)$.

For the uniqueness, first notice that, by Proposition \ref{prop:perm},  for each pair $\{i,j\}\in \{1,...,k\}$ with $i\neq j$ there is a unique $f_{ij}$ such that  $\{f, f_i,f_j,f_{ij}\}$ is a Bianchi $(L_i, L_j)$-quadrilateral. By Theorem \ref{teo:cuboor},
there exists a unique Bianchi $k$-cube $({\cal C}_0,...,{\cal C}_k)$ such that ${\cal C}_0=\{f\}, \ {\cal C}_1=\{f_1,...,f_k\}$ and ${\cal C}_2=\{f_{ij}\}_{1\leq i\neq j\leq k}$. 
\qed

\begin{remark}\label{re:assfamily}\po\emph{$(i)$ It is worth pointing out that, once the 
$L_i$-transforms $f_i={\cal R}_{\varphi_i,\beta_i}(f)$, $1\leq i\leq k$, of  $f\colon M^n(c)\rightarrow \mathbb{Q}^{n+p}(\tilde c)$  are known (meaning that  all the pairs $(\varphi_i,\beta_i)$, $1\leq i\leq k$, are known),  and hence the first two families ${\cal C}_0=\{f\}$ and ${\cal C}_1=\{f_1,...,f_k\}$ in the Bianchi $(L_1,...,L_k)$-cube $({\cal C}_0, . . . , {\cal C}_k)$ are given, then all $2^k-(k+1)$ elements in the remaining families  ${\cal C}_2,\ldots, {\cal C}_k$ are determined by explicit algebraic formulae, with no integrations required. Namely,  if $\Omega\in \Gamma((\R^k)^*\otimes \R^k)$ is defined by $\Omega=\sum \Omega_{ij} e^j\otimes e_i$, where the $\Omega_{ij}$, $1\leq i,j\leq k$,  are  given  by the algebraic expressions (\ref{eq:omegaij}), which only depend on the $(\varphi_i,\beta_i)$, $1\leq i\leq k$, then all elements of ${\cal C}_r$, $2\leq r\leq k$, are given  explicitly  by 
$$ f_{\alpha_r}=\ral_{\varphi^{\alpha_{r}},\beta^{\alpha_{r}},\Omega^{\alpha_{r}},L^{\alpha_{r}}}(f)=f-{\cal \Fes}^{\alpha_r}(\Omega^{\alpha_{r}})^{-1}\varphi^{\alpha_{r}}
$$
%\begin{eqnarray*} f_{\alpha_r}&=&\ral_{\varphi^{\alpha_{r}},\beta^{\alpha_{r}},\Omega^{\alpha_{r}},L^{\alpha_{r}}}(f)\nonumber\\
%&=&f-{\cal \Fes}^{\alpha_r}(\Omega^{\alpha_{r}})^{-1}\varphi^{\alpha_{r}}
%\end{eqnarray*}
where $\alpha_r\in \{i_1,...,i_r\}\in \Lambda_r$, $(\varphi^{\alpha_{r}},\beta^{\alpha_{r}},\Omega^{\alpha_{r}})$ is given by (\ref{eq:varbetaalp}) and 
${\cal \Fes}^{\alpha_r}= f_*(d\varphi^{\alpha_{r}})^t+\beta^{\alpha_{r}}$.\vspace{1ex}\\
$(ii)$ In Theorem \ref{lcubo}, if $f_i$ belongs to the associated family determined by $\{f_j, f_l\}$,  $1\leq i\neq j\neq l\neq i\leq k$, that is, there exist $a_j, a_l\in \R$ such that $\varphi_i=a_j\varphi_j+a_l\varphi_l$ and 
$\beta_i=a_j\beta_j+a_l\beta_l$, then from $\Phi_i=a_j\Phi_j+a_l\Phi_l$ and (\ref{eq:liphii}) one obtains that
$A_\xi=\lambda I$
 for
 \be\label{vecxi}\xi=a_j(C_i-C_j)\beta_j+a_l(C_i-C_l)\beta_l\,\,\,\mbox{and}\,\,\,\lambda=(\tilde c-c)(a_j(C_i-C_j)\varphi_j+a_l(C_i-C_l)\varphi_l),\ee
%$\xi=k_j(C_i-C_j)\beta_j+k_l(C_i-C_l)\beta_l$ and $\lambda=(\tilde c-c)(k_j(C_i-C_j)\varphi_j+k_l(C_i-C_l)\varphi_l$, 
with 
$C_i=-L_i^{-1}$ for $1\leq i\leq k$. It is easily seen that this can not happen if $f$ satisfies the assumptions of either Proposition 
\ref{thm:hiebetaij} or \ref{teo:curvdif}, depending on whether $c=\tilde c$ or $c\neq \tilde c$, respectively, if the codimension attains its minimum possible values $p=n$ and $p=n-1$, respectively. Therefore, in these cases the last assumption in Theorem  \ref{lcubo} can be dropped.}
\end{remark}

%\section{The $P$-transformation for flat Lagrangian submanifolds}
\section{The $P$-transformation}

In this section we obtain  reductions of the $L$-transformation that preserve the classes  of $n$-dimensional flat Lagrangian submanifolds of  $\mathbb{C}^n=\R^{2n}$ and $n$-dimensional  submanifolds with constant sectional curvature $c$ of  $\Sf_\epsilon^{2n+1}(c)$ that are horizontal with respect to the Hopf fibration $\pi\colon \Sf_\epsilon^{2n+1}(c)\to \tilde M^n(4c)$. %Here $\Sf_\epsilon^{2n+1}(c)$ stands for either the standard Euclidean sphere or the anti-de-Sitter space time of dimension  $(2n+1)$ and constant sectional curvature $c$, corresponding to $\epsilon=1$ or $\epsilon =-1$, respectively, and $\tilde M^n(4c)$  denotes either the complex projective space $\mathbb C\mathbb P^n(4c)$ or complex hyperbolic space $\mathbb C\mathbb H^n(4c)$ of complex dimension n and constant holomorphic sectional curvature 4c, corresponding to $c>0$ or $c<0$, respectively.   
We also obtain a further reduction  that preserves the class  of $n$-dimensional flat Lagrangian submanifolds of  $\R^{2n}$ that are contained in $\Sf^{2n-1}$. 

\subsection{The $P$-transformation for flat Lagrangian submanifolds}

An isometric immersion $f\colon M^n\rightarrow \tilde M^n$ of an $n$-dimensional Riemannian manifold into a Kaehler manifold of complex dimensional $n$ is said to be \emph{Lagrangian} if the almost
complex structure of $\tilde M^m$ carries each tangent space of $M^n$ onto its corresponding
normal space. 
%If in addition $n=m$ then $f$ is said to be \emph{Lagrangian}. 

For a Lagrangian isometric immersion $f\colon M^n \to \mathbb C^n$, comparing  normal and tangential  components of $\tilde \nabla_X Jf_*Y= J\tilde \nabla_Xf_*Y$ yields, respectively,
\begin{equation}\label{eq:dtlflat16}
\nabla^\bot_X Jf_*Y= Jf_*\nabla_XY
\end{equation}
and
\begin{equation}\label{eq:dtlflat17}
-f_*A_{Jf_*Y}X=J\alpha(X,Y)
\end{equation}   
for all $X,Y\in \mathfrak{X}(M)$. It follows from (\ref{eq:dtlflat16}) that 
$$
R^\bot(X,Y)Jf_*Z=Jf_*R(X,Y)Z
$$
for all $X,Y, Z\in \mathfrak{X}(M)$. In particular,  $M^n$ is flat if and only if $f$ has flat normal bundle.\vspace{1ex}

Lagrangian isometric immersions $f\colon M^n(0) \to \mathbb C^n$ with $\nu_f\equiv 0$ have been characterized as follows in \cite{dt3} among isometric immersions $f\colon M^n(0) \to \mathbb C^n$ with $\nu_f\equiv 0\equiv R^\perp$.

\begin{theorem}\po \label{teo:mtlflat10}
If $f\colon M^n(0)\rightarrow \mathbb R^{2n}\cong \mathbb C^n$ is an isometric immersion with  $\nu_f\equiv 0\equiv R^\perp$  and  $(v,h)$ is its associated  solution of (\ref{eq:i}), then $f$ is Lagrangian if and only if $h=h^t$.
\end{theorem}

%\begin{theorem}\po \label{teo:mtlflat10}
%Let $f\colon M^n(0)\rightarrow \mathbb R^{2n}\cong \mathbb C^n$ be an isometric immersion with flat normal bundle and $\nu_f\equiv 0$, and let $(v,h)$ the solution of (\ref{eq:i}) associated to $f$. Then $f$ is Lagrangian if and only if $h=h^t$.
%\end{theorem}

\begin{corollary}\po \label{cor:cocompsl}
Let  $f\colon M^n(0)\rightarrow \mathbb R^{2n}$ be a Lagrangian isometric immersion with $\nu_f\equiv~0$. Then there exist  locally principal coordinates $u_1,..., u_n$ on $M^n(0)$ and positive smooth  functions  $v_1,...,v_n$ such that $ds^2=\sum_j v_j^2du_j^2$ and $\alpha\left(\partial_i,\partial_j\right)=\delta_{ij}J\partial_i$,
%$$
%ds^2=\sum_j v_j^2du_j^2, \ \ \ \ \mbox{and}\ \ \ \alpha\left(\partial_i,\partial_j\right)=\delta_{ij}J\partial_i,
%$$
where $v = (v_1,..., v_n)$ and $h = (h_{ij})$, $h_{ij}=v_i^{-1}\partial_i v_j$,  satisfy the  system of PDE's 
\begin{equation}\label{eq:corlagsist}
\left\{
\begin{array}{lc}
i) \ \partial_j v_i = h_{ji}v_j,\;\;1\leq i\neq j\leq n, & ii) \ (\sum_{\ell=1}^n \partial_\ell) h_{ij}=0,\vspace{0.1cm}\\
iii) \ \partial_\ell h_{ij} = h_{i\ell}h_{j\ell}, \;\;1\leq \ell \neq i\neq j\neq \ell\leq n.
\end{array}
\right.
\end{equation}
%(resp., i), ii), iii) of (\ref{eq:corlagsist}) and iv) $\partial_iv_i=-\sum_{j\neq i}h_{ij}v_j$), with $i\neq j\neq k\neq i.$
Conversely, if $(v, h)$ is a solution of  (\ref{eq:corlagsist}) on an open simply connected subset $U \subset \R^n$ such that  $v_i \neq 0$ everywhere 
%(resp., with initial conditions at some point $x_0\in U$ chosen so that $\sum_iv_i^2(x_0)=1$). Let $(f,Y_1,...,Y_n)$, with $f,Y_i:U\rightarrow \mathbb C^{n}$, be a solution of 
%\begin{equation}\label{eq:ivl}\left\{
%\begin{array}{lc}
%i) \ \partial_i( f) = v_iY_i & ii) \ \partial_j( Y_i)=h_{ij}Y_j, \ i\neq j,\\
%iii) \ \partial_i( Y_i) = -\sum_{k\neq i}h_{ki}Y_k + iY_i,&\\
%\end{array}
%\right.
%\end{equation}
%with initial conditions $(Y_1(u_0),...,Y_n(u_0))$ at some point  $u_0\in U$ chosen so that 
%\begin{eqnarray*}\nonumber
%&\left<Y_i(u_0),Y_j(u_0)\right>=\left<iY_i(u_0),Y_j(u_0)\right>=0, \ i\neq j, \ \left<Y_i(u_0),Y_i(u_0)\right> = 1.\nonumber\\
%\end{eqnarray*}
then there exists a Lagrangian immersion $f:U\rightarrow \mathbb C^n$  with flat induced metric $ds^2=\sum_iv_i^2du_i$ and $\nu_f\equiv 0$.
% (resp., $f(U)\subset \mathbb S^{2n-1}$).
\end{corollary}

%We now introduce a special class of vectorial Ribaucour transforms of Lagrangian isometric immersions $f\colon M^n(0)\to \R^{2n}$, which will be shown to be a subclass of the class of $L$-transforms of $f$.

\begin{definition}\label{df:Ptransformada}\po \emph{Let  $\tilde f={\cal
R}_{\va,\beta,\Omega}(f)\colon\,\tilde{M}^n\to \R^{2n}$ be a vectorial Ribaucour transform of a Lagrangian isometric immersion $f\colon M^n(0)\to \R^{2n}$
determined by $(\va,\beta,\Omega)$, with $\va\in \Gamma(V)$, $\beta\in \Gamma(V^*\otimes N_fM)$ and $\Omega\in \Gamma(V^*\otimes V)$. 
If there exists $P\in V^*\otimes V$ satisfying $\sigma(P)\cap (-\sigma(P))=\emptyset$, where $\sigma(P)$ denotes the set of (complex) eigenvalues of $P$,  such that
\be\label{eq:bJ}
\beta=Jf_*\omega^tP
\ee
where $\omega=d\varphi$, and 
\begin{equation}\label{eq:lyaplag}
\Omega^tP+P^t\Omega^t+T\rho=0
\end{equation}
where $T=-P^t-(P^t)^{-1}$ and $\rho=P^t\omega\omega^tP$, then $\tilde f$  is called a  {\em $P$-vectorial Ribaucour transform\/} of $f$, or simply a  {\em $P$-transform\/} of $f$. 
We write  $\tilde f={\cal
R}_{\va,\beta, \Omega, P}(f)$. }
\end{definition}

\begin{remark}\po\label{re:ptransform}{\em As in the case of the $L$-transformation, if $\<\cdot , \cdot\>^\sim$ is another inner product on $V$, with 
 %\be\label{eq:innerp}
 $ \<\cdot, \cdot\>^\sim = \<\cdot, B\cdot\>$
 for some invertible  $B\in V^*\otimes V$, then conditions (\ref{eq:bJ}) and (\ref{eq:lyaplag}) are satisfied by $(\va, \beta, \Omega, P,  \<\cdot, \cdot\>)$ if and only if they are satisfied by $(\hat\va, \hat\beta, \hat\Omega, \hat P,  \<\cdot, \cdot\>^\sim)$, where $\hat\va=\va$, $\hat\beta=\beta B$, $\hat\Omega=\Omega B$ and $\hat P=B^{-1}PB$. Thus, in the definition of the $P$-transformation one may replace  $(\va, \beta, \Omega, P,  \<\cdot, \cdot\>)$  by any quintuple  $(\hat\va, \hat\beta, \hat\Omega, \hat P,  \<\cdot, \cdot\>^\sim)$  related to $(\va, \beta, \Omega, P,  \<\cdot, \cdot\>)$ in this way.
 As before,  the transposes of $\hat P$ and $P$ with respect to $\<\cdot, \cdot\>^\sim$ and $\<\cdot, \cdot\>$ coincide, and if $P$ is   diagonalizable then one can always assume that $P$ is symmetric after an appropriate  change of the inner product on $V$. }
\end{remark}

%Let $f\colon M^n(0)\to \R^{2n}$ be a Lagrangian isometric immersion and let   
%$P$ be an endomorphism of an Euclidean vector space $V$ satisfying $\sigma(P)\cap (-\sigma(P))=\emptyset$, where $\sigma(P)$ denotes the set of eigenvalues of $P$. A vectorial Ribaucour transform $\tilde f ={\cal
%R}_{\va,\beta,\Omega}(f)$ of $f$ is a  {\em $P$-vectorial Ribaucour transform\/} of $f$, or simply a  {\em $P$-transform\/} of $f$, if 
%\be\label{eq:bJ}
%\beta=Jf_*\omega^tP
%\ee
%where $\omega=d\varphi$, and 
%\begin{equation}\label{eq:lyaplag}
%\Omega^tP+P^t\Omega^t+T\rho=0
%\end{equation}
%where $T=-P^t-(P^t)^{-1}$ and $\rho=P^t\omega\omega^tP$. We write  $\tilde f={\cal
%R}_{\va,\beta, \Omega, P}(f)$. }
%If, in addition,  $f(M)\subset \mathbb S^{2n-1}\subset \mathbb R^{2n}$ and $\varphi=-P^t\omega \sum_i\partial_i$, where $\partial_1, \ldots, \partial_n$ are the coordinate vectors with respect to the local coordinates $(u_1,\ldots, u_n)$ given by Corollary \ref{cor:cocompsl}, then we say that the pair $(\omega,P)$ determines a {\em $P^*$-transform\/} of $f$ and we write $\tilde f={\cal R}_{\omega,P}(f).$ 
%\end{definition}

Given endomorphisms $P,C\in V^*\otimes V$ of a Euclidean  vector space $V$, the  equation
\begin{equation}\label{eq:lyapunov}
P^tX+XP=C
\end{equation} 
for $X\in V^*\otimes V$ is known as the \emph{Lyapunov equation}. Thus, if $\tilde f={\cal
R}_{\va,\beta, \Omega,P}(f)$ is a $P$-transform of $f$, then $\Omega$ satisfies a Lyapunov-type equation. The following result of  \cite{hc} provides a sufficient condition on  $P$ in order that (\ref{eq:lyapunov}) admit a unique solution $X$, and gives an explicit expression for $X$ as a polynomial on  $P$ and $C$.
 
\begin{theorem}\label{lema:antiomega1lag}\po
If $P\in  V^*\otimes V$ satisfies $\sigma (P)\cap (-\sigma(P))=\emptyset$, then (\ref{eq:lyapunov}) admits a unique solution for any $C\in V^*\otimes V$, given by 
\begin{equation}\label{eq:solx}
X=(q_{-P}(P^t))^{-1}\sum_{\ell=1}^k a_\ell\sum_{i=0}^{\ell-1}(-1)^i(P^t)^{\ell-1-i}CP^i,
\end{equation}
where
$
q_{-P}(x)=\sum_{\ell=0}^k a_\ell x^\ell, \ a_k=1,
$
is the characteristic polynomial of $-P$. 
%$$
%\eta(P,C)=
%\sum_{k=1}^n a_k((A^t)^kX-(-1)^kXA^k)=
%\sum_{k=1}^n a_k\sum_{i=0}^{k-1}(-1)^i(P^t)^{k-1-i}CP^i.
%$$
%Moreover, this solution is a polynomial of the matrices  $A$ and $C$.
\end{theorem}
%It is well-known that (\ref{eq:lyapunov}) has a unique solution $X$ for any $C\in\End(V)$ if $\sigma(A)\cap (-\sigma(A))=\emptyset$.

In order to prove that the $P$-transformation  preserves the class  of $n$-dimensional flat Lagrangian submanifolds of  $\mathbb{C}^n=\R^{2n}$ we need the special case for $c=0$ of the following result. The case $c\neq 0$  will  be used in our study in the next section of the $P$-transformation for $n$-dimensional  submanifolds with constant sectional curvature $c$ of  $\Sf_\epsilon^{2n+1}(c)$ that are horizontal with respect to the Hopf fibration $\pi\colon \Sf_\epsilon^{2n+1}(c)\to \tilde M^n(4c)$.

\begin{proposition}\label{teo:allagr}\po
Let  $V$ and $W$ be Euclidean vector spaces, let  $P\in V^*\otimes V$ be such that $\sigma(P)\cap (-\sigma(P))=\emptyset$, let $c\in \R$, $\psi\in V^*$ and $\nu\in V^*\otimes  W$ and 
let $X\in V^*\otimes V$ be the unique solution of the Lyapunov equation
\begin{equation}\label{eq:partsimcomplag}
X^tP+P^tX^t=-T\rho
\end{equation}
where $T=-P^t-(P^t)^{-1}$ and $\rho=P^tQP$, with $Q=\nu^t\nu+c\psi^t\psi$. 
Then  $X$ satisfies 
\begin{equation}\label{eq:partsimlag}
X + X^t=Q+\rho,
\end{equation} 
\begin{equation}\label{eq:tcomx}
TX-X^tT^t=0
\end{equation}
and 
$$
XL+L^tX^t=\rho,\,\,\,\mbox{where}\,\,L=(P^2+I)^{-1}P^2.
$$
%Moreover, if $c\geq 0$  then $X$ is invertible if and only if $\ker (P^c-\alpha I)\cap Y^c=\{0\}$ for any eigenvalue $\alpha$ of $P^c$, where 
%\be\label{eq:s}Y= \left\{ \begin{array}{l} \ker \nu,\,\,\,\mbox{if}\,\,\, c=0;\\
%\ker \psi\cap \ker \nu, \,\,\,\mbox{if}\,\,\,\,c\neq 0.
%\end{array}\right.
%\ee
\end{proposition}
\proof 
Since
$$
%\begin{array}{lll}
(Q+\rho)P+P^t(Q+\rho)+2T\rho
%&=&(Q+\rho)P+P^t(Q+\rho)-2P^t\rho-2QP\\
%&=&-QP+\rho P+P^tQ-P^t\rho\\
=(-QP+\rho P)-(-QP+\rho P)^t
%&=&(-QP+P^tQP^2)-(-QP+P^tQP^2)^t
%\end{array}
$$
is a skew-symmetric matrix,  the symmetric part of the  unique solution  $X$   of (\ref{eq:partsimcomplag}) must coincide with  $(Q+\rho)/2$. This gives (\ref{eq:partsimlag}).
% is the (unique) solution of Lyapunov equation 
%$$XP+P^tX=-(T\rho)_s,$$
%namely,  $(Q+\rho)/2$ is symmetric part of the (unique) solution of (\ref{eq:partsimcomplag}), it shows (\ref{eq:partsimlag}).
Equation (\ref{eq:tcomx}) follows easily from  (\ref{eq:partsimcomplag}) and (\ref{eq:partsimlag}).
%$$
%\begin{array}{lll}
%TX-X^tT^t &=&-P^tX-(P^t)^{-1}X+X^tP+X^tP^{-1}\\ 
%&=&-P^tX-(P^t)^{-1}X-P^tX^t- T\rho  + X^tP^{-1}\\ 
%&=&-P^tX-(P^t)^{-1}X-P^tX^t- T\rho  - (P^t)^{-1}X^t-(P^t)^{-1}T\rho P^{-1}\\ 
%&=& - P^t(X+X^t)-(P^t)^{-1}(X+X^t)+P^t\rho+ (P^t)^{-1}\rho+\rho P^{-1}+(P^{t})^{-2}\rho P^{-1}\\ 
%&=&-P^t(Q+\rho)-(P^t)^{-1}(Q+\rho)+P^t\rho+(P^t)^{-1}\rho+\rho P^{-1}+(P^t)^{-2}\rho P^{-1}\\
%&=&-P^tQ-(P^t)^{-1}+P^tQ+(P^t)^{-1}\\
%&=&0.
%\end{array}
%$$
Now observe that $T^t=-P-P^{-1}=-P^{-1}(I+P^2)$. 
Since $\sigma(P)\cap (-\sigma(P))=\emptyset$, it follows that 
% we have that $-1\not \in \sigma(P^2)$, soon 
$T$ is invertible and 
\begin{equation}\label{eq:tlp}
T^tL=-P.
\end{equation}
Therefore
\begin{eqnarray*}
T(X L +L^tX^t)&=&TX L +TL^tX^t=X^tT^tL+TL^tX^t\\
&=& - X^tP-P^tX^t=T\rho. \qed
\end{eqnarray*}

\begin{lemma}\label{lem:lconstl}\po
Any $P$-transform   
%$\tilde f=\ral_{\varphi,P}(f)$  
of a Lagrangian isometric immersion $f\colon M^n(0)\to \mathbb R^{2n}$  is also an $L$-transform of $f$ with $L=(P^2+I)^{-1}P^2.$
\end{lemma}
\proof It follows from (\ref{eq:bJ}) that
$\nabla_X\beta=Jf_*(\nabla_X\omega^t)P$ for all $X \in \mathfrak{X}(M)$, hence (\ref{eq:alphao}) becomes
\begin{equation}\label{eq:omega1}
\alpha(X,\omega^t(v))+Jf_*\left((\nabla_X\omega^t)Pv\right)=0,
\end{equation}
which can also be writtten as 
%\begin{equation}\label{eq:dop}
$$
d\omega^t(X)P=f_*^tJA(X)^t\omega^t.
$$
%\end{equation}
On the other hand,  from (\ref{eq:dtlflat17})  we obtain
%\begin{equation}\label{eq:fja}
$
f_*^tJA(X)^t=-A(X)Jf_*Y
$
%\end{equation}
for all $X \in \mathfrak{X}(M)$. 
Therefore
%using (\ref{eq:fja}) and (\ref{eq:dop}), we obtain that
$$
A(X)\beta=A(X)Jf_*\omega^tP= -f_*^tJA(X)^t\omega^tP= -d\omega^t(X)P^2.
$$
%$$
%\begin{array}{lll}A(X)\beta&=&A(X)Jf_*\omega^tP\vspace{1ex}\\
%&=& -f_*^tJA(X)^t\omega^tP\vspace{1ex}\\
%&=& -d\omega^t(X)P^2.
%\end{array}
%$$
Thus, the tensor 
$$\Phi(X)v=(\nabla_X \omega^t)v-A_{\beta(v)}X=d\omega^t(X)v-A(X)\beta(v)$$
satisfies
\begin{equation}\label{eq:pdef}
\Phi(X)=d\omega^t(X)(P^2+I),
\end{equation} 
and hence 
%\begin{equation}\label{eq:pab}
$\Phi(X)L+A(X)\beta=0$
%\end{equation}
for all $X \in \mathfrak{X}(M)$. Finally, that $\Omega$ satisfies (\ref{eq:cond3}) follows from Proposition \ref{teo:allagr}.\vspace{1ex}\qed

%Now, from (\ref{eq:dop}) and (\ref{eq:pdef}) we obtain
%$$
%\begin{array}{l}
%d(\Omega^tP+P^t\Omega^t+T\rho)(X)=\Phi(X)^t\omega^tP+P^t\Phi(X)^t\omega^t-(P^t)^2d\omega(X)\omega^tP-d\omega(X)\omega^tP\vspace{1ex}\\
%\hspace*{25ex} -(P^t)^2\omega d\omega^t(X)P-\omega d\omega^t(X)P\vspace{1ex}\\
%\hspace*{25ex} = \Phi(X)^t\omega^tP-((P^t)^2+I) d\omega(X)\omega^tP\vspace{1ex}\\
%\hspace*{27ex} + P^t\Phi(X)^t\omega^t-((P^t)^2+I)\omega d\omega^t(X)P\vspace{1ex}\\
%\hspace*{25ex} =0. \qed
%\end{array} 
%$$
%By Definition \ref{df:ribconst}, we need to prove that (\ref{eq:cond11nf}) and (\ref{eq:cond3}) are true. The first is consequence of (\ref{eq:pab}), and the second derive from Proposition \ref{teo:allagr}.\qed

It follows from Lemma \ref{lem:lconstl} and Theorem \ref{thm:csc} that if $\tilde f={\cal R}_{\va, \beta, \Omega, P}(f)\colon \tilde M^n\to \R^{2n}$ is a $P$-transform of a Lagrangian isometric immersion $f\colon M^n(0)\to \R^{2n}$ then $\tilde M^n$ is also flat. We shall prove in Theorem \ref{thm:lagrang} below that $\tilde f$ is also Lagrangian. First we  express equation (\ref{eq:omega1}) in the local principal coordinates given by 
Corollary \ref{cor:cocompsl}.

\begin{proposition}\label{cor:condint}\po
Let $f\colon M^n(0)\to  \mathbb{R}^{2n}$ be a Lagrangian isometric immersion with $\nu_f\equiv 0$, let $(u_1, \ldots, u_n)$ be  principal  coordinates  given by Corollary \ref{cor:cocompsl} on an open subset  $U\subset M^n(0)$,  let $(v,h)$ be the solution of (\ref{eq:corlagsist}) associated to $f$ and let $P$ be an invertible endomorphism of a Euclidean vector space $V$. If $\varphi\colon U\to V$ is such that $\omega=d\varphi$ satisfies (\ref{eq:omega1}), then $(\varphi, \gamma_1, \ldots, \gamma_n)$, with  $\gamma_i=v_i^{-1}\omega(\partial_i)$ for $1\leq i\leq n$, is a solution of
\begin{equation}\label{ral_{0l}}
\left\{
\begin{array}{ll}
i) \ \partial_i\varphi=v_i\gamma_i, &
ii) \ \partial_i\gamma_i=-(P^t)^{-1}\gamma_i - \sum_{j\neq i}h_{ij}\gamma_j,\\ iii) \ \partial_i\gamma_j = h_{ij}\gamma_i \ j\neq i.\\
\end{array}
\right.
\end{equation}
Conversely, if $(\varphi, \gamma_1, \ldots, \gamma_n)$ is a solution of (\ref{ral_{0l}}) then $\omega=d\varphi$ satisfies (\ref{eq:omega1}) and $\omega(\partial_i)=v_i\gamma_i$ for $1\leq i\leq n$. 
\end{proposition} 

\begin{theorem}\label{thm:lagrang}\po
If $f\colon M^n(0)\to \R^{2n}$ is a Lagrangian isometric immersion with $\nu_f\equiv 0$ then any $P$-transform $\tilde f={\cal R}_{\va,\beta, \Omega, P}(f)\colon \tilde M^n(0)\to \R^{2n}$   of $f$ is also Lagrangian.
% Then $\tilde f$ is also flat Lagrangian. Moreover, if $f(M)\subset \mathbb S^{2n-1}$ and $\tilde f$ is an $Pe$-transform of $f$, then $\tilde f(\tilde M)\subset \mathbb S^{2n-1}$.
\end{theorem}

\proof
%By Lemma \ref{lem:lconstl} and Theorem \ref{thm:csc} $\tilde M$ is flat. It follows from 
According to Theorem \ref{teo:mtlflat10}, it suffices to prove that the pair $(\tilde v, \tilde h)$ associated to $\tilde f$ with respect to the local principal coordinates $u_1, \ldots, u_n$ given by Proposition \ref{thm:hiebetaij} satisfies $\tilde h^t=\tilde h$. In view of (\ref{eq:hbetatrans}), this is equivalent to
% that $P$-transforms $\tilde f$ of $f$ is also Lagrangian if and only if 
\begin{equation}\label{eq:lagigual}
\left<L^{-1}\Omega^{-1}\gamma_i,\beta^j\right>=\left<\beta^i,L^{-1}\Omega^{-1}\gamma_j\right>
\end{equation}
where $\gamma_1,...,\gamma_n$,  $\beta^1,...,\beta^n$ are given by  (\ref{eq:betagamma}). By (\ref{eq:tcomx}), we have
\begin{equation}\label{eq:tomega2}
T\Omega=\Omega^tT^t
\end{equation}
where  $T=-P^t-(P^t)^{-1}$. From  (\ref{eq:betagamma}) and  (\ref{eq:bJ}) we obtain
\begin{equation}\label{eq:relac1}
P^t\gamma_i =P^t\omega(X_i)=\beta^tJf_*X_i=\beta^t\xi_i= \beta^i
\end{equation}
where $X_i=v_i^{-1}\partial_i$, $1\leq i\leq n$. We conclude from (\ref{eq:tlp}), (\ref{eq:tomega2}) and (\ref{eq:relac1}) that (\ref{eq:lagigual}) holds.\vspace{1ex}\qed

In order to prove the existence of $P$-transforms of any Lagrangian isometric immersion $f\colon M^n(0)\to  \mathbb{R}^{2n}$ with $\nu_f\equiv 0$, we need the special case $c=0$ of the folowing algebraic result, whose general form will be used in the next section.

\begin{proposition}\po\label{prop:padmi} Let  $V$ and $W$ be Euclidean vector spaces, let $P\in V^*\otimes V$ be such that $\sigma(P)\cap(-\sigma(P))=\emptyset$, let $\psi\in V^*$,  $\nu\in V^*\otimes W$,  $c\in \R$, and let $X=X(\nu,\psi,c)$ be the unique solution of (\ref{eq:partsimcomplag}). 
\begin{itemize}
\item[(i)] If $c\geq 0$, then $X$ is invertible if and only if $\ker(P^c-\alpha I)\cap Y^c=\{0\}$ for any eigenvalue $\alpha$ of $P^c$, where $Y^c$  is the complexified of
%\be\label{eq:YY}
$$Y= \left\{ \begin{array}{l} \ker \nu,\,\,\,\mbox{if}\,\,\, c=0;\\
\ker \psi\cap \ker \nu, \,\,\,\mbox{if}\,\,\,\,c\neq 0.
\end{array}\right.
$$
%\ee
%\item[(ii)] If $c=0$, then the pair $\nu$ is $P$-admissible if and only if $\ker(P^c-\alpha I)\cap (\ker \nu)^c=\{0\}$ for any eigenvalue $\alpha$ of $P^c$.
\item[(ii)] If $c<0$ and $\ker(P^c-\alpha I) \cap (\ker\nu)^c= \{0\}$ (resp., $\ker(P^c-\alpha I) \cap (\ker\psi)^c= \{0\}$) 
for  any eigenvalue $\alpha$ of $P^c$,  then there exists $\epsilon>0$ such that $X(\nu,\psi,c)$ is invertible if $|\psi|<\epsilon$ (resp., $|\nu|<\epsilon$).
%\item[(iii)] If $c<0$, then
%\begin{itemize}
%\item[(a)] if  $\ker(P^c-\alpha I) \cap (ker\nu)^c= \{0\}$ for  any eigenvalue $\alpha$ of $P^c$,  then there exists $\epsilon>0$ such that $(\psi,\nu)$ is P-admissible if $|\psi|<\epsilon$.
%\item[(b)] if $\ker(P^c-\alpha I) \cap (ker\psi)^c= \{0\}$ for  any eigenvalue $\alpha$ of $P^c$,  then there exists $\epsilon>0$ such that $(\psi,\nu)$ is P-admissible if $|\nu|<\epsilon.$
%\end{itemize}
\end{itemize}
\end{proposition}

%\begin{proposition}\po\label{prop:padmi} Let  $V$ and $W$ be Euclidean vector spaces, let $P\in \End(V)$ be such that $\sigma(P)\cap(-\sigma(P))=\emptyset$, and let  $\psi\in V^*$ and $\nu\in \Hom(V,W)$. \vspace{1ex}
%\begin{itemize}
%\item[(i)] If $c>0$, then the pair $(\psi, \nu)$ is $P$-admissible if and only if $\ker(P^c-\alpha I)\cap S^c=\{0\}$ for any eigenvalue $\alpha$ of $P^c$, where $S=\ker \psi\cap \ker \nu$ and $S^c$ is the complexified of $S$.
%\item[(ii)] If $c=0$, then the pair $\nu$ is $P$-admissible if and only if $\ker(P^c-\alpha I)\cap (\ker \nu)^c=\{0\}$ for any eigenvalue $\alpha$ of $P^c$.
%\item[(iii)] If $c<0$ and $\ker(P^c-\alpha I) \cap (ker\nu)^c= \{0\}$ (respectively, $\ker(P^c-\alpha I) \cap (ker\psi)^c= \{0\}$) 
%for  any eigenvalue $\alpha$ of $P^c$,  then there exists $\epsilon>0$ such that $(\psi,\nu)$ is P-admissible if $|\psi|<\epsilon$ (respectively, $|\nu|<\epsilon$).
%\item[(iii)] If $c<0$, then
%\begin{itemize}
%\item[(a)] if  $\ker(P^c-\alpha I) \cap (ker\nu)^c= \{0\}$ for  any eigenvalue $\alpha$ of $P^c$,  then there exists $\epsilon>0$ such that $(\psi,\nu)$ is P-admissible if $|\psi|<\epsilon$.
%\item[(b)] if $\ker(P^c-\alpha I) \cap (ker\psi)^c= \{0\}$ for  any eigenvalue $\alpha$ of $P^c$,  then there exists $\epsilon>0$ such that $(\psi,\nu)$ is P-admissible if $|\nu|<\epsilon.$
%\end{itemize}
%\end{itemize}
%\end{proposition}
\proof  $(i)$ Assume first that $\ker(P^c-\alpha I)\cap Y^c=\{0\}$ for any eigenvalue $\alpha$ of $P^c$.
%To prove the last assertion, assume that $c\geq 0$  and that  $\ker (P^c-\alpha I)\cap Y^c=\{0\}$ for any eigenvalue $\alpha$ of $P^c$. Then, 
If  $u\in \ker X^t$, we obtain from (\ref{eq:partsimlag}) that
$$
0=|\nu u|^2+  c (\psi u)^2+|\nu Pu|^2+ c(\psi Pu)^2.
$$
Since $c\geq 0$, it follows that $u\in Y$. Thus $\ker X^t \subset Y$, and   (\ref{eq:partsimlag}) implies that $\ker X^t=\ker X$. 
Then  (\ref{eq:partsimcomplag}) by $P(\ker X)\subset \ker X$, and since $\ker X\subset Y$,  $\ker X=\{0\}$.
% if $\ker (P^c-\alpha I)\cap Y^c=\{0\}$ for any eigenvalue $\alpha$ of $P^c$. %if  $E_\alpha\cap S=\{0\}$ (resp., $E_\alpha\cap S^c=\{0\}$) for all real (resp., complex) eigenvalue $\alpha$ of $P$. 

Conversely, if $u \in \ker(P^c-\alpha I)\cap Y^c$, then  
$X_s u=0$ and $P^tX_a(u)=-\alpha X_a(u)$ by (\ref{eq:partsimcomplag}) and (\ref{eq:partsimlag}). Since $\sigma(P)\cap (-\sigma(P))=\emptyset$, it follows that  $X_a(u)=0$,  hence $X(u)=0$. Thus $u=0$.

 $(ii)$ First observe that the unique solution $X(\nu,\psi,c)$ of (\ref{eq:partsimcomplag}) given by (\ref{eq:solx}) depends continuously  on $(\psi, \nu,c)$. Assume first that $\ker(P^c-\alpha I) \cap (\ker\nu)^c= \{0\}$. It follows from case $(i)$ for $c=0$ that $X(\nu,0,c)$ is invertible.  Thus, there exists $\epsilon>0$ such that $X(\nu,\psi,c)$ is invertible for all $\psi$ with $|\psi|<\epsilon$.

 Suppose now that $\ker(P^c-\alpha I) \cap (\ker\psi)^c= \{0\}$ for  any eigenvalue $\alpha$ of $P^c$.  From part $(i)$ for $c>0$, 
 we have that 
 %For (b), setting $\nu=0$ it follows analogous of the $c>0$ with $\nu=0$, that 
 $X(0,\psi,c)$ is invertible.
 % since  $E_\alpha\cap ker \psi=\{0\}$ (resp., $E_\alpha \cap (ker \psi)^c= \{0\}$), for all real (resp., complex) eigenvalue $\alpha$ of $P$. 
 Therefore, there exists $\epsilon>0$ such that $X(\nu,\psi,c)$ is invertible for all $\nu$ with $|\nu|<\epsilon$. \vspace{1ex}\qed  

\begin{theorem}\po Let $f\colon M^n(0)\to  \mathbb{R}^{2n}$ be a Lagrangian isometric immersion with $\nu_f\equiv 0$ 
%(respectively, whit $f(M^n(0))\subset \mathbb S^{2n-1}$), 
and let $P$ an endomorphism of a Euclidean vector space $V$ such that $\sigma(P)\cap(-\sigma(P))=\emptyset$.
Fixed  $x_0\in M^n$, let $\varphi_0\in V$ and  $\omega_0\in T^*_{x_0}M\otimes V$ 
%(respectively, $\omega_0\in T^*_{x_0}M\oplus V$) 
be such that $\ker (P^c-\alpha I)\cap \ker (\omega_0^t)^c=\{0\}$ for any eigenvalue $\alpha$ of $P^c$.
%$E_\alpha\cap ker \omega_0^t=\{0\}$ 
%(respectively, $E_\alpha\cap ker (\omega_0^t)^c=\{0\}$) 
%for all real (respectively, complex)  eigenvalues of $P$.
Then there exist an open neighborhood $U$ of $x_0$ and a unique $P$-transform 
$\tilde f={\cal
R}_{\va,\beta, \Omega, P}(f|_U)$  
%(respectively, $Pe$-transform $\tilde f={\cal R}_{\omega,P}(f|_U)$) 
of $f|_U$ such that $\varphi(x_0)=\varphi^0$ and $d\varphi(x_0)=\omega_0$.
% (respectively, $\omega(x_0)=\omega_0$).
\end{theorem}
\proof Let $(u_1, \ldots, u_n)$ be principal  coordinates given by Corollary \ref{cor:cocompsl} on an open simply-connected neighborhood $U\subset M^n(0)$ of $x_0$ and   let $(v,h)$ be the solution of (\ref{eq:corlagsist}) associated to $f$. It is easily checked that the compatibility conditions of system (\ref{ral_{0l}}) are satisfied by virtue of (\ref{eq:corlagsist}). Therefore, if $\gamma^0_i=v_i^{-1}(x_0)\omega_0(\partial_i(x_0))$ for $1\leq i\leq n$, then there exists a unique solution $(\varphi, \gamma_1, \ldots, \gamma_n)$ of (\ref{ral_{0l}}) such that $\varphi(x_0)=\varphi^0$ and $\gamma_i(x_0)=\gamma_i^0$ for all  $1\leq i\leq n$.  By 
Proposition \ref{cor:condint}, $\omega=d\varphi$ satisfies (\ref{eq:omega1}) and $\omega(\partial_i)=v_i\gamma_i$ for $1\leq i\leq n$.

%Let $\varphi: U\rightarrow V$ be a solution of (\ref{eq:omega1}) on an open simply connected subset $U$ of $x_0\in M$, such that $\varphi(x_0)=\varphi_0$ and $d\varphi(x_0)=\omega_0$. 
Now, by Proposition \ref{teo:allagr} the Lyapunov equation 
$
X^tP+P^tX^t=-T\rho_0,
$
where $T=-P^t-(P^t)^{-1}$ and $\rho_0=P^t\omega_0\omega_0^tP,$ has a unique solution $X=\Omega^0$,  and  $2\Omega^0_s={\cal F}^t(x_0){\cal F}(x_0)$. Moreover,    $\Omega^0$ is invertible by  Proposition \ref{prop:padmi}. 
Define $\beta\in \Gamma(V^*\otimes N_fM)$ by
(\ref{eq:bJ}), and let $\Omega$ be the unique solution of (\ref{eq:o1}) such that $\Omega(x_0)=\Omega_0$. Shrinking $U$ if necessary, we may assume that $\Omega$ is invertible on $U$. 
Using  (\ref{eq:pdef}) we obtain
$$
\begin{array}{l}
d(\Omega^tP+P^t\Omega^t+T\rho)(X)=\Phi(X)^t\omega^tP+P^t\Phi(X)^t\omega^t-(P^t)^2d\omega(X)\omega^tP-d\omega(X)\omega^tP\vspace{1ex}\\
\hspace*{27ex} -(P^t)^2\omega d\omega^t(X)P-\omega d\omega^t(X)P=0.
%\vspace{1ex}\\
%\hspace*{25ex} = \Phi(X)^t\omega^tP-((P^t)^2+I) d\omega(X)\omega^tP\vspace{1ex}\\
%\hspace*{27ex} + P^t\Phi(X)^t\omega^t-((P^t)^2+I)\omega d\omega^t(X)P=0.
%\hspace*{25ex} =0.
\end{array}
$$
Since (\ref{eq:lyaplag}) holds at $x_0$,  it holds on $U$, hence $\tilde f={\cal
R}_{\va,\beta,\Omega}(f|_U)$ is a $P$-transform of $f|_U$.\vspace{1ex}
\qed

In the following corollary we summarize the process given by the  $P$-transformation to generate new Lagrangian isometric immersions $f\colon M^n(0)\to \R^{2n}$ starting with a given one  and a vector-valued solution of the linear system of PDE's (\ref{ral_{0l}}).

\begin{corollary}\label{cor:exlag}\po
Let $f\colon M^n(0)\rightarrow \mathbb{R}^{2n}$ be a Lagrangian isometric immersion with $\nu_f\equiv 0$,   let $(v,h)$ be the solution of (\ref{eq:corlagsist}) associated to $f$ on an open subset $U\subset M^n(0)$ endowed with principal coordinates $(u_1, \ldots, u_n)$ given by Corollary \ref{cor:cocompsl},  let $P$ be an endomorphism of a Euclidean vector space $V$ such that $\sigma(P)\cap (-\sigma(P))=\emptyset$, and 
let $(\varphi,\gamma_1,...,\gamma_n)$ be  a $V$-valued solution  of the linear system of PDE's  (\ref{ral_{0l}})
%\begin{equation}\label{Ral_2l}
%\left\{
%\begin{array}{ll}
%i) \ \partial_i\varphi=v_i\gamma_i, &
%ii) \ \partial_i\gamma_i=-(P^t)^{-1}\gamma_i - \sum_{j\neq i}\gamma_jh_{ij},\\ iii) \ \partial_i\gamma_j = h_{ij}\gamma_i \ j\neq i.\\
%\end{array}
%\right.
%\end{equation}
on an open subset $W\subset U$ where 
\be\label{eq:vtilde}
\tilde v_i= v_i+\left<\gamma_i,(P+P^{-1})\Omega^{-1}\varphi\right>
\ee
does not vanish for all $1\leq i\leq n$, with 
$\Omega^t$ given by (\ref{eq:solx}) for  $C=((P^t)^2+I)d\varphi(d\varphi)^tP$.
% and $A=P$.
%$$A=P\,\,\,\,\mbox{and}\,\,\,\,C=({P^t}^2+I)d\varphi(d\varphi)^tP.$$
 Then $\tilde f:W\rightarrow \mathbb R^{2n}=\mathbb{C}^n$ given by
$$
\begin{array}{l}
\tilde f =f-\sum_{j=1}^n\left(\left<\Omega^{-1}\varphi,\gamma_j\right>+\left<P\Omega^{-1}\varphi,\gamma_j\right>i\right)f_*X_j,\\
\end{array}
$$
where $X_j=v_j^{-1}\partial_j$ for $1\leq j\leq n$, defines a new  Lagrangian isometric immersion with flat induced metric.
Moreover, the pair associated to $\tilde f$ is $(\tilde v,\tilde h)$,   where
%\begin{equation}\label{eq:htilpe}
$$
\tilde{h}_{ij}= h_{ij}+\left<\gamma_j,(P+P^{-1})\Omega^{-1}\gamma_i\right>.
$$
%\end{equation}
In particular, $(\tilde v,\tilde h)$ is a new  solution of  (\ref{eq:corlagsist}).
\end{corollary}

%\begin{corollary}
%Let $(v,h)$ be a solution of (\ref{eq:corlagsist}) and let $(\varphi,\gamma)$ be a solution of (\ref{ral_{0l}}). Then a new solution $(\tilde v,\tilde h)$ of (\ref{eq:corlagsist}) is given in Corollary \ref{cor:exlag}.
%\end{corollary}

\subsection{The $P^*$-transformation}

Lagrangian isometric immersions $f\colon M^n(0) \to \mathbb C^n$ satisfying $f(M^n(0))\subset \mathbb S^{2n-1}$ are of special interest in view of the following result of \cite{dt1}.

\begin{theorem}\po
An isometric immersion  $f\colon M^n(0)\to \R^{2n}$ is Lagrangian and satisfies  $f(M^n(0))\subset \mathbb S^{2n-1}$ if and only if it is the lifting by the Hopf projection $\pi:\mathbb S^{2n-1}\rightarrow \mathbb C\mathbb P^{n-1}$ of a Lagrangian isometric immersion $F\colon M^{n-1}(0)\rightarrow \mathbb C \mathbb P^{n-1}$.
\end{theorem}

We will need the following characterization obtained in \cite{dt2} of isometric immersions $f:M^n(0)\rightarrow \mathbb R^{2n}$ with flat normal bundle  and $\nu_f\equiv 0$ that satisfy $f(M^n(0))\subset \mathbb S^{2n-1}$. 

\begin{corollary}\label{cor:m0r}\po
Let $f:M^n(0)\rightarrow \mathbb R^{2n}$ be an isometric immersion with  $\nu_f\equiv 0\equiv R^\perp$, and let $(v,h)$ be its associated pair with respect to the  local principal coordinates $(u_1,...,u_n)$ on $M^n(0)$ given by Proposition \ref{thm:hiebetaij}. Then
$f(M^n(0))\subset \mathbb S^{2n-1}$ if and only if $\sum_{i=1}^nv_i^2=1$.
%\begin{equation}\label{somavzinhos}
%\sum_iv_i^2=1.
%\end{equation} 
\end{corollary}

We now look for the $P$-transformations that preserve the class of Lagrangian isometric immersions $f\colon M^n(0)\to \mathbb R^{2n}$ such that $f(M^n(0))\subset \mathbb S^{2n-1}$.

\begin{lemma}\po \label{le:vap} Let $f\colon M^n(0)\to \R^{2n}$ be a Lagrangian isometric immersion such that $f(M^n(0))\subset \mathbb S^{2n-1}\subset \mathbb R^{2n}$ and   let $(v,h)$ be its associated  solution of (\ref{eq:corlagsist})  on an open subset $U\subset M^n(0)$ endowed with principal coordinates  given by Corollary \ref{cor:cocompsl}. If $\tilde f ={\cal
R}_{\va,\beta, \omega,P}(f)$  is  a $P$-transform of $f$ and $\gamma_i=v_i^{-1}\omega(\partial_i)$,  $1\leq i\leq n$, then  $ \varphi+ \sum_i v_iP^t\gamma_i
$ 
is constant on~$U$.
\end{lemma}
\proof Using (\ref{ral_{0l}}) and $\sum_{i=1}^n v_i^2=1$ we obtain 
%Define $\bar\varphi=+P^t\omega \sum_i\partial_i$.  It is enough to show that $\partial_i\bar\varphi=\omega(\partial_i)$. Using (iv) of Corollary \ref{cor:cocompsl}, (i) and (ii) of (\ref{{cal R}_{0l}}), we have
$$
\begin{array}{l}
\partial_i(\varphi+ \sum_{j=1}^n v_jP^t\gamma_j)=v_i\gamma_i+P^t\sum_{i\neq j}(v_j\partial_i\gamma_j+\gamma_j\partial_iv_j)+v_iP^t\partial_i\gamma_i+\partial_iv_iP^t\gamma_i\vspace{.2cm}\\
\hspace*{18.4ex}=v_i\gamma_i+P^t\sum_{i\neq j}\left(h_{ij}\gamma_iv_j+\gamma_jh_{ij}v_i\right)-P^t(P^t)^{-1}\gamma_iv_i\vspace{1ex}\\\hspace*{20ex}-P^t\sum_{i\neq j}\gamma_jh_{ij}v_i-P^t\gamma_i\sum_{j\neq i}h_{ij}v_j
%\vspace{.2cm}\\
%\hspace*{18.4ex}
=0.\,\,\,\qed
\end{array}
$$

\begin{definition}\label{df:Ptransformadab}\po \emph{Let $f\colon M^n(0)\to \R^{2n}$ be a Lagrangian isometric immersion such that $f(M^n(0))\subset \mathbb S^{2n-1}\subset \mathbb R^{2n}$ and   let $(v,h)$ be its associated  solution of (\ref{eq:corlagsist})  on an open subset $U\subset M^n(0)$ endowed with principal coordinates  given by Corollary \ref{cor:cocompsl}. A $P$-transform  $\tilde f ={\cal
R}_{\va,\beta,\Omega,P}(f)$ of $f|_U$ is said to be a    {\em $P^*$-transform\/}  of $f|_U$ if 
 \be\label{eq:p*}\varphi+ \sum_{i=1}^n v_iP^t\gamma_i=0.\ee}
%  where $\partial_1, \ldots, \partial_n$ are the coordinate vectors with respect to the local coordinates $(u_1,\ldots, u_n)$ given by Corollary \ref{cor:cocompsl}. We write $\tilde f={\cal R}_{\omega,P}(f)$. }
\end{definition}
\begin{remark}\po\label{re:ptransformb}{\em Notice that  (\ref{eq:p*}) is satisfied by $(\va, \beta, \Omega, P,  \<\cdot, \cdot\>)$ if and only if it is satisfied by $(\hat\va, \hat\beta, \hat\Omega, \hat P,  \<\cdot, \cdot\>^\sim)$, where $ \<\cdot, \cdot\>^\sim = \<\cdot, B\cdot\>$ 
 for some invertible  $B\in V^*\otimes V$, $\hat\va=\va$, $\hat\beta=\beta B$, $\hat\Omega=\Omega B$ and $\hat P=B^{-1}PB$, for the transposes of $\hat P$ and $P$ with respect to $\<\cdot, \cdot\>^\sim$ and $\<\cdot, \cdot\>$ coincide.}
\end{remark}  
%\begin{theorem}
%Let $f\colon M^n(0)\to \R^{2n}$ be a Lagrangian isometric immersion. If  $f(M)\subset \mathbb S^{2n-1}$, then any $P^*$-transform  $\tilde f ={\cal
%R}_{\varphi,\beta, \Omega, P}(f)\colon \tilde{M}^n(0)\to \R^{2n}$   of $f$ also satisfies $\tilde f(\tilde M)\subset \mathbb S^{2n-1}$.
%\end{theorem}

\begin{theorem}\po
Let $f\colon M^n(0)\to \R^{2n}$ be a Lagrangian isometric immersion such that $f(M^n(0))\subset \mathbb S^{2n-1}$. If  $\tilde f ={\cal
R}_{\varphi,\beta, \Omega, P}(f)\colon \tilde{M}^n(0)\to \R^{2n}$ is a $P^*$-transform of $f$, then $\tilde f(\tilde M^n(0))\subset \mathbb S^{2n-1}$.
\end{theorem}
\proof Since $\sum_{i=1}^n v_i^2=1$, then $\tilde v$  given by (\ref{eq:vtilde}) also satisfies $\sum_{i=1}^n \tilde{v}_i^2=1$ 
%same holds for the pair $(\tilde v, \tilde h)$  associated to $\tilde f$, where $\tilde v$ is given by (\ref{eq:vtilde}), 
 if and only if
%\begin{equation}\label{eq:somape}
\begin{eqnarray*}
0&=&\sum_{i=1}^n 2v_i\left<\gamma_i,T^t\Omega^{-1}\varphi\right>+\sum_{i=1}^n \left<\gamma_i,T^t\Omega^{-1}\varphi\right>^2\vspace{1ex}\\
&=&-\sum_{i=1}^n 2v_i\left<\omega^t(T^t\Omega^{-1}\varphi),X_i\right>+\left<\omega^t(T^t\Omega^{-1}\varphi),\omega^t(T^t\Omega^{-1}\varphi)\right>\vspace{1ex}\\
&=&\<T^t\Omega^{-1}\varphi, -2\sum_{i=1}^n v_i\gamma_i\>+\<\varphi, (\Omega^{-1})^tT\omega\omega^tT^t\Omega^{-1}\varphi).
\end{eqnarray*}
%\end{equation}
%It follows from $\omega(X_j)=\gamma_j$ that $\omega^t(T^t\Omega^{-1}\varphi)=\sum_j\left<T^t\Omega^{-1}\varphi,\gamma_j\right>X_j$.
That this is indeed the case follows from the fact that, by  (\ref{eq:lyaplag}) and (\ref{eq:p*}), the last expression is equal to
$$
\begin{array}{l}
%-\sum_i2v_i\left<\gamma_i,T^t\Omega^{-1}\varphi\right>+\sum_i\left<\gamma_i,T^t\Omega^{-1}\varphi\right>^2\vspace{1ex}\\
%\hspace*{18ex}=  -\sum_i2v_i\left<\omega^t(T^t\Omega^{-1}\varphi),X_i\right>+\left<\omega^t(T^t\Omega^{-1}\varphi),\omega^t(T^t\Omega^{-1}\varphi)\right>\vspace{0.2cm}\vspace{1ex}\\
%\hspace*{18ex}= \<T^t\Omega^{-1}\varphi, -2\sum_iv_i\gamma_i\>+\<\varphi, (\Omega^{-1})^tT\omega\omega^tT^t\Omega^{-1}\varphi)
%\>\vspace{1ex}\\
%\hspace*{18ex}=
 \<T^t\Omega^{-1}\varphi, 2(P^t)^{-1}\varphi\>+\<\varphi, (\Omega^{-1})^t(-(P^t)^{-1}\Omega^t-\Omega^tP^{-1})T^t\Omega^{-1}\varphi\>\vspace{1ex}\\
\hspace*{18ex}= \<T^t\Omega^{-1}\varphi, 2(P^t)^{-1}\varphi\>-\<\varphi, (\Omega^{-1})^t(P^t)^{-1}T\varphi\>-\<\varphi, P^{-1}T^t\Omega^{-1}\varphi\>, 
\end{array}
$$
which vanishes by (\ref{eq:tomega2}).
\qed

\begin{corollary}\po\label{cor:p*}
Let $f\colon M^n(0)\rightarrow \mathbb{R}^{2n}$ be a Lagrangian isometric immersion satisfying  $f(M^n(0))\subset  \mathbb S^{2n-1}$, let   $(v,h)$ be its associated solution of (\ref{eq:corlagsist}), with $\sum_{i=1}^nv_i^2=1$, on an open subset $U\subset M^n(0)$ endowed with principal coordinates  given by Corollary \ref{cor:cocompsl}, and let $P$ be an endomorphism of a Euclidean vector space $V$ such that $\sigma(P)\cap (-\sigma(P))=\emptyset$. If  $(\gamma_1,...,\gamma_n)$ satisfies equations $(ii)$ and $(iii)$ of (\ref{ral_{0l}}) 
%be a $V$-valued solution of the linear system  of PDE's
%$$
%\left\{
%\begin{array}{l}
%i) \ \partial_i\gamma_i=-(P^t)^{-1}\gamma_i - \sum_{j\neq i}\gamma_jh_{ij},\\ ii) \ \partial_i\gamma_j = h_{ij}\gamma_i, \ j\neq i,\\
%\end{array}
%\right.
%$$
on an open subset $W\subset U$ where 
\begin{equation}\label{eq:vtilpe}
\tilde v_i= v_i-\sum_{j=1}^nv_j\<\Omega^{-1}\gamma_i,((P^t)^2+I)\gamma_j\>
\end{equation}
%\begin{equation}\label{eq:vtilpe}
%\tilde v_i= v_i-\sum_k\left<\gamma_i,(P+P^{-1})\Omega^{-1}P^tv_k\gamma_k\right>
%\end{equation}
does not vanish for  $1\leq i\leq n$, with 
$\Omega^t$ given by (\ref{eq:solx}) for $C=((P^t)^2+I)\sum_{i=1}^n\gamma_i\gamma_i^tP$,
% and $A=P$,
%$$A=P\,\,\,\,\mbox{and}\,\,\,\,C=-({P^t}^2+I)P^t\sum_{i}\gamma_i\gamma_i^tP.$$
 then $\tilde f\colon W\rightarrow \mathbb R^{2n}$ given by
$$
\begin{array}{l}
\tilde f =f+\sum_{\ell,j=1}^n\left(\left<\Omega^{-1}P^t\gamma_\ell,\gamma_j\right>+\left<\Omega^{-1}P^t\gamma_\ell,P^t\gamma_j\right>i\right)v_\ell f_*X_j,\\
\end{array}
$$
where $X_j=v_j^{-1}\partial_j$ for $1\leq j\leq n$, defines a new  Lagrangian isometric immersion with flat induced metric such that $\tilde f(W)\subset \mathbb S^{2n-1}$. Moreover, $(\tilde v,\tilde h)$ is the solution of (\ref{eq:corlagsist}) associated to $\tilde f$, with
 $$
\tilde{h}_{ij}= h_{ij}+\left<\gamma_j,(P+P^{-1})\Omega^{-1}\gamma_i\right>.
$$
\end{corollary}

\subsection{The $P$-transformation for horizontal submanifolds.} 

%\noindent{\bf\large \S 1.2.2 Horizontal submanifolds of constant sectional curvature.} \vspace{3ex}

Let $\mathbb C^{n+1}_\epsilon$ denote the complex  $(n+1)$-space endowed with the pseudo-Euclidean metric
$$
g_\epsilon = \epsilon dz_1d\bar z_1+\sum_{j=2}^{n+1}dz_jd\bar z_j, \ \epsilon = \pm 1,
$$ 
and let
$$
\mathbb{S}_\epsilon^{2n+1}(c)=\{z\in \mathbb C^{n+1}_\epsilon: \ g_\epsilon(z,z)=\frac{1}{c}, \ \epsilon c>0\}
$$
stand for either the  Euclidean sphere or the anti-de-Sitter space time of dimension  $(2n+1)$ and constant sectional curvature $c$, depending on whether $\epsilon=1$ or $\epsilon =-1$, respectively. 
 
The complex numbers act on $\mathbb C_\epsilon^{n+1}$ by
$
\lambda(z_1,...,z_{n+1})\rightarrow (\lambda z_1,...,\lambda z_{n+1}).
$
The quotient space $\tilde M^n(4c)$ of $\mathbb S^{2n+1}(c)$ under the identification induced by this action is the complex projective space $\mathbb C\mathbb P^n(4c)$ or complex hyperbolic space $\mathbb C\mathbb H^n(4c)$ of complex dimension~$n$ and constant holomorphic  curvature $4c$, corresponding to $c>0$ or $c<0$, respectively. If
$\pi\colon\mathbb S_\epsilon^{2n+1}\to \tilde M^n(4c)$
is the quotient map,  $\tilde J$ is the complex structure on $\mathbb C^{n+1}_\epsilon$ defined by multiplication by $i$ and  $\phi$ is its projection onto the tangent bundle of $\mathbb S^{2n+1}(c)$, then the complex structure $J$ on $\tilde M^n(4c)$ is given by $J\circ \pi_*=\pi_*\circ \phi$.

 An isometric immersion $f\colon M\rightarrow \mathbb S_\epsilon^{2n+1}(c)\subset \mathbb C^{n+1}_\epsilon$ of a Riemannian manifold is said to be \emph{horizontal}   if the \emph{structure vector field} $\xi_f=\sqrt{|c|} \tilde{J} f$ is everywhere normal to  $f$.
If $f\colon M\rightarrow \mathbb S_\epsilon^{2n+1}(c)$  is horizontal, then $f$ is \emph{anti-invariant} with respect to $\phi$, that is, $\phi f_*T_xM\subset N_fM(x)$ for all $x\in M$. Moreover, the second fundamental form of $f$ satisfies
\be\label{eq:sffhor}
A_{\xi_f}=0\,\,\,\mbox{and}\,\,\,\,\phi\alpha(X,Y)=-f_*A_{\phi f_*Y}X
\ee
for all $X,Y\in \mathfrak{X}(M)$, and the following relations hold:
\begin{eqnarray}
\nabla^\perp_X\xi_f&=&\sqrt{|c|}\phi f_*X\label{eq:napxif}\vspace{1ex}\\
\phi^2 X&=&-X+\epsilon\<X,\xi_f\>\xi_f\label{eq:phi2}\vspace{1ex}\\
\<\phi X, \phi Y\>&=&\<X,Y\>+\epsilon \<X, \xi_f\>\<Y, \xi_f\>\label{eq:phiimp}\vspace{1ex}\\
\nabla_X^\perp \phi f_*Y&=&\phi f_*\nabla_XY-\epsilon\sqrt{|c|}\<X,Y\>\xi_f\label{eq:conhor1}\vspace{1ex}\\
R^\perp(X,Y)\xi_f&=&0\nonumber
%\label{eq:conhor2}
\vspace{1ex}\\
\<R^\perp(X,Y)\phi f_*Z,\phi f_*W\>&\!=\!&\<R(X,Y)Z,W\>-c(\<X,W\>\<Y,Z\>-\<X,Z\>\<Y,W\>)\nonumber
%\label{eq:conhor3}
\end{eqnarray}
In particular, a horizontal isometric immersion $f\colon M^n\rightarrow \mathbb S_\epsilon^{2n+1}(c)$ has flat normal bundle if and only if $M^n$ has constant sectional curvature $c$.  \vspace{1ex}

The next result relates Lagrangian isometric immersions $f\colon M\rightarrow \tilde M^n(4c)$  to  horizontal isometric immersions $f\colon M\rightarrow \mathbb S_\epsilon^{2n+1}(c)$. We refer to  \cite{re} for a proof. 

\begin{theorem}\label{teo:re}\po
If $f\colon M^n\rightarrow \mathbb S_\epsilon^{2n+1}(c)$ is horizontal then $g=\pi\circ f$ is Lagrangian. Conversely, if $g\colon M\rightarrow \tilde M^n(4c)$ is a Lagrangian isometric immersions and  $(x_0,y_0)\in M\times \mathbb S_\epsilon^{2n+1}(c)$ is such that   $g(x_0)=\pi(y_0)$, then  there exist a Riemannian manifold $\hat M$, an isometric covering map $\tau:\hat M\rightarrow M$, a horizontal isometric immersion $\hat f:\hat M\rightarrow \mathbb S_\epsilon^{2n+1}(c)$ and  $\hat x\in \hat M$ such that $\pi\circ \hat f=g\circ \tau,$ $\tau (\hat x)=x_0$ and $\hat f(\hat x)=y_0.$
\end{theorem}

It will be convenient to have Proposition \ref{thm:hiebetaij} explicitly restated in this particular case. 

\begin{corollary}\label{cor:hiebetaij}\po Let $f\colon M^n(c) \to \mathbb S_\epsilon^{2n+1}(c)$ be an isometric immersion with flat normal bundle, $\nu_f \equiv 0$, and Riemannian first normal bundle when $c<0$. Then there exist  local principal coordinates  $(u_1,...,u_n)$ on $M^n(c)$, a smooth orthonormal normal frame $\xi_1,...,\xi_{n+1}$ and smooth functions $v_1,...,v_n$ and $\rho_1,...,\rho_n$, with $v_i>0$ for $1\leq i\leq n$,  such that 
%\begin{equation}\label{alphasistnov}
$$
ds^2=\sum_i v_i^2du_i^2, \ \ \alpha\left(\partial_i,\partial_j\right)=v_i\delta_{ij}\xi_i,
$$
%\end{equation}
and
\begin{equation}\label{formadcovnov}
\nabla_{\partial_i}X_j= h_{ji}X_i, \ \nabla^\bot_{\partial_i}\xi_j=h_{ij}\xi_i, i\neq j, \ \nabla^\bot_{\partial_i}\xi_{n+1}=\rho_i\xi_i,
\end{equation}
where $X_i = (1/v_i)\partial_i$ and $h_{ij} = (1/v_i
)\partial_i v_j$ for $i\neq j$. Moreover, the triple $(v, h,\rho)$,
where $v = (v_1,..., v_n)$, $h = (h_{ij})$ and $\rho=(\rho_1,...,\rho_n)$, satisfies the system of PDEs
\begin{equation}\label{eq:iii}
\left\{
\begin{array}{l}
i) \ \partial_j v_i = h_{ji}v_j, \ \ \ ii) \ \partial_j h_{i\ell}=h_{ij}h_{j\ell}, \ \ \ iii) \ \partial_j \rho_i=h_{ij}\rho_j\vspace{0.1cm}\\ 
iv)\ \partial_i h_{ij}+\partial_j h_{ji}+\sum_{\ell=1}^n h_{\ell i}h_{\ell j}+cv_iv_j=0,\vspace{0.1cm}\\
v) \ \partial_j h_{ij}+\partial_i h_{ji}+\sum_{\ell=1}^n h_{i\ell }h_{j\ell}+\epsilon\rho_i\rho_j=0, \ \epsilon = \frac{c}{|c|}, \ i\neq j \neq \ell\neq i.\vspace{0.1cm}\\
\end{array}
\right.
\end{equation}

Conversely, if $(v, h,\rho)$ is a solution of (\ref{eq:iii}) on an open simply connected subset $U \subset \R^n$
such that $v_i > 0$ everywhere, then there exists an immersion $f: U\rightarrow  \mathbb S_\epsilon^{2n+1}(c)$ with flat normal bundle, $\nu_f \equiv 0$, Riemannian first normal bundle and induced metric $ds^2 = \sum_iv_i^2 du_i^2$ of constant sectional curvature $c$.
\end{corollary}

 We shall use the following results  proved in \cite{t}.

 \begin{theorem}\label{teo:toj6}\po
An  isometric immersion $f\colon M^n(c)\rightarrow \mathbb S_\epsilon^{2n+1}(c)$ as in Corollary \ref{cor:hiebetaij} is horizontal if and only if its associated triple $(v,h,\rho)$ satisfies 
\begin{equation}\label{eq:toj12}
h_{ij}=h_{ji} \ \mbox{and} \ \rho_i=\sqrt{|c|}v_i.
\end{equation}
\end{theorem}

\begin{corollary}\label{cor:toj7}\po
Let $f\colon M^n(c)\rightarrow \mathbb S_\epsilon^{2n+1}(c)$ be a horizontal isometric immersion with $\nu_f\equiv 0.$ Then there exist  locally principal coordinates $(u_1,...,u_n)$ on $M^n(c)$ such that
$$
ds^2=\sum_iv_i^2du_i^2, \ v_i>0 \ \mbox{and} \ \alpha\left(\partial_i,\partial_j\right)=\delta_{ij}\phi \partial_i,
$$ 
where $v=(v_1,...,v_n)$ and $h=(h_{ij})$, with $h_{ij}=h_{ji}$ for $1\leq  i\neq j\leq n$, satisfy the  system of PDE's
\begin{equation}\label{eq:iii2}
\left\{
\begin{array}{lc}
i) \ \partial_j v_i = h_{ji}v_j, \;\;1\leq  i\neq j\leq n,& ii) \ (\sum_{\ell=1}^n\partial_\ell) h_{ij}+cv_iv_j=0,\vspace{0.1cm}\\
iii) \ \partial_\ell h_{ij} = h_{i\ell}h_{j\ell}, \;\; 1\leq \ell\neq i\neq j\neq \ell\leq n.\vspace{0.1cm}\\
\end{array}
\right.
\end{equation}

Conversely, if $(v,h)$ is a solution of (\ref{eq:iii2}) on an open simply connected subset $U\subset \mathbb R^n$ such that $v_i\neq 0$ for $1\leq i\leq n$ everywhere, 
%. Let $(F,Y_1,...,Y_n)$, with $F,Y_i:U\rightarrow \mathbb C^{n+1}_\epsilon$, be a solution of system of PDE's
%\begin{equation}\label{eq:iv}\left\{
%\begin{array}{lc}
%i) \ \partial_i( F) = v_iY_i & ii) \ \partial_j( Y_i)=h_{ij}Y_j, \ i\neq j,\\
%iii) \ \partial_i( Y_i) = -\sum_{k\neq i}h_{ki}Y_k + iY_i-cv_iF,&\\
%\end{array}
%\right.
%\end{equation}
%with initial conditions $(F(u_0),Y_1(u_0),...,Y_n(u_0))$ at some point $u_0\in U$ chosen so that
%$$
%\begin{array}{ll}
%&\left<Y_i(u_0),Y_j(u_0)\right>=\left<iY_i(u_0),Y_j(u_0)\right>=0, \ i\neq j, \ \left<Y_i(u_0),Y_i(u_0)\right> = 1, \\
%&\left<F(u_0),Y_i(u_0)\right>=\left< iF(u_0),Y_i(u_0)\right>=0 \ \mbox{and} \ \left<F(u_0),F(u_0)\right>=\frac{1}{c}.\\
%\end{array}
%$$
then there exists a 
%Then $F(U)\subset \mathbb S^{2n+1}_\epsilon (c)\subset \mathbb C^{n+1}_\epsilon$ and the 
horizontal isometric immersion $f:U\rightarrow \mathbb S_\epsilon^{2n+1}(c)$ 
%given by $F=i\circ f$, is horizontal and has 
with induced metric $ds^2=\sum v_i^2du_i$ of constant sectional curvature $c$.
\end{corollary}

\begin{definition}\label{df:Ptransformadanf}\po \emph{Let  $\tilde f={\cal
R}_{\va,\beta,\Omega}(f)\colon\,\tilde{M}^n\to \mathbb{S}_\epsilon^{2n+1}(c)$ be a vectorial Ribaucour transform of a  horizontal  isometric immersion 
$f\colon M^n(c)\rightarrow \mathbb{S}_\epsilon^{2n+1}(c)$ determined by $(\va,\beta,\Omega)$, with $\va\in \Gamma(V)$, $\beta\in \Gamma(V^*\otimes N_fM)$ and $\Omega\in \Gamma(V^*\otimes V)$. 
If there exists $P\in V^*\otimes V$ satisfying $\sigma(P)\cap (-\sigma(P))=\emptyset$  such that
\be\label{eq:betahor}
\beta=(\phi f_*\omega^t+\epsilon\sqrt{|c|}\xi_f\varphi^t)P
\ee
where $\omega=d\varphi$, and 
\begin{equation}\label{eq:lyaplagnf}
\Omega^tP+P^t\Omega^t+T\rho=0
\end{equation}
where $T=-P^t-(P^t)^{-1}$ and $\rho=P^t\omega\omega^tP +cP^t\varphi\varphi^t P=\beta^t\beta$, then $\tilde f$  is said to be 
a  {\em Ribaucour $P$-transform\/} of $f$, or simply a  {\em $P$-transform\/} of $f$. We write  $\tilde f={\cal
R}_{\va,\beta,\Omega,P}(f)$.}
\end{definition}

%Remark \ref{re:ptransform} also holds in this case. The analogue of Lemma \ref{lem:lconstl} is the following.

\begin{lemma}\label{lem:lconst}\po
Any $P$-transform  of a horizontal isometric immersion $f\colon\,M^n(c)\to \mathbb S_\epsilon^{2n+1}(c)$ is also an $L$-transform of $f$  with $L=(P^2+I)^{-1}P^2.$
\end{lemma}
\proof It follows from (\ref{eq:conhor1}) and  (\ref{eq:betahor}) that
$$\nabla_X\beta=\phi f_*(\nabla_X\omega^t)P+\epsilon \sqrt{|c|}\nabla_X^\perp\xi_f\varphi^tP$$
for all $X \in \mathfrak{X}(M)$, hence (\ref{eq:alphao}) becomes
\begin{equation}\label{eq:omega1b}
\alpha(X,\omega^t(v))+\phi f_*(\nabla_X\omega^t)Pv+\epsilon \sqrt{|c|}\nabla_X^\perp\xi_f\varphi^tPv=0,
\end{equation}
which by (\ref{eq:napxif}) and (\ref{eq:phiimp}) can also be writtten as 
\begin{equation}\label{eq:dopb}
d\omega^t(X)P=-f_*^t\phi^tA(X)^t\omega^t-cX\varphi^tP.
\end{equation}
On the other hand,  from (\ref{eq:sffhor}), (\ref{eq:phi2}) and (\ref{eq:phiimp})  we obtain
\begin{equation}\label{eq:fja}
f_*^t\phi^tA(X)^t=
%f_*^tJ\alpha(X,Y)=-f_*^tf_*A_{Jf_*Y}X=
A(X)\phi f_*
\end{equation}
for all $X\in \mathfrak{X}(M)$. 
Therefore, by (\ref{eq:sffhor}), (\ref{eq:dopb}) and (\ref{eq:fja}),
%using (\ref{eq:fja}) and (\ref{eq:dop}), we obtain that
\begin{eqnarray*}
A(X)\beta&=&A(X)(\phi f_*\omega^t+\epsilon \sqrt{|c|}\xi_f\varphi^t)P\vspace{1ex}\\
&=&f_*^t\phi^tA(X)^t\omega^tP\vspace{1ex}\\
&=&-d\omega^t(X)P^2-cX\varphi^tP^2.
\end{eqnarray*}
%$$
%\begin{array}{lll}A(X)\beta&=&A(X)Jf_*\omega^tP\vspace{1ex}\\
%&=& -f_*^tJA(X)^t\omega^tP\vspace{1ex}\\
%&=& -d\omega^t(X)P^2.
%\end{array}
%$$
Thus the tensor 
$\Phi(X)=
%v=(\nabla_X \omega^t)v-A_{\beta(v)}X+cX\varphi^t(v)=
d\omega^t(X)-A(X)\beta+cX\varphi^t$
satisfies
\begin{equation}\label{eq:pdef2}
\Phi(X)=(d\omega^t(X)+cX\varphi^t)(P^2+I),
\end{equation} 
and hence 
%\begin{equation}\label{eq:pab}
$
\Phi(X)L+A(X)\beta=0$
%\end{equation}
for all $X \in \mathfrak{X}(M)$. Finally, that $\Omega$ satisfies (\ref{eq:cond3}) follows from Proposition \ref{teo:allagr}.\vspace{1ex}\qed

It follows from Theorem \ref{thm:csc} and Lemma \ref{lem:lconst}  that if $\tilde f={\cal R}_{\va,\beta, \Omega, P}(f)\colon \tilde M^n\to \mathbb{S}_\epsilon^{2n+1}(c)$ is a $P$-transform of a horizontal isometric immersion $f\colon M^n(c)\to \mathbb{S}_\epsilon^{2n+1}(c)$ then $\tilde M^n$  also has constant sectional curvature $c$. We shall prove that $\tilde f$ is also horizontal. For that we first express equation (\ref{eq:omega1b}) in the local principal coordinates given by 
Corollary \ref{cor:toj7}.

\begin{proposition}\label{cor:condintb}\po
Let $f\colon M^n(c)\to \mathbb{S}_\epsilon^{2n+1}(c)$ be a horizontal isometric immersion  with $\nu_f\equiv~0$, 
let $(v,h)$ be its associated solution of (\ref{eq:iii2}) with respect to  principal  coordinates given by Corollary \ref{cor:toj7},
%let $(v,h)$ be the solution of (\ref{eq:iii2}) associated to $f$ 
and let $P$ be an invertible endomorphism of a Euclidean vector space $V$. If $(\varphi, \omega=d\varphi)$ satisfies (\ref{eq:omega1b}), then $(\varphi, \gamma_1, \ldots, \gamma_n)$, with  $\gamma_i=v_i^{-1}\omega(\partial_i)$ for $1\leq i\leq n$, is a solution of
\begin{equation}\label{ral_{h}}
\left\{
\begin{array}{ll}
i) \ \partial_i\varphi=v_i\gamma_i, &
ii) \ \partial_i\gamma_i=-(P^t)^{-1}\gamma_i - \sum_{j\neq i}\gamma_jh_{ij}-cv_i\varphi,\\ 
iii) \ \partial_i\gamma_j = h_{ij}\gamma_i \ j\neq i,\\
\end{array}
\right.
\end{equation}
Conversely, if $(\varphi, \gamma_1, \ldots, \gamma_n)$ is a solution of (\ref{ral_{h}}) then $(\varphi, \omega=d\varphi)$ satisfies (\ref{eq:omega1b}) and $\omega(\partial_i)=v_i\gamma_i$ for $1\leq i\leq n$. 
\end{proposition} 

\begin{theorem}\label{thm:hor}\po
If $f\colon M^n(c)\to \mathbb{S}_\epsilon^{2n+1}(c)$ is a horizontal isometric immersion with $\nu_f\equiv 0$ then any $P$-transform $\tilde f={\cal R}_{\va,\beta, \Omega, P}(f)\colon \tilde M^n(c)\to \mathbb{S}_\epsilon^{2n+1}(c)$   of $f$ is also horizontal.
\end{theorem}
\proof
By Theorem \ref{teo:toj6}, it suffices to prove that the triple $(\tilde v,\tilde h,\tilde \rho)$ associated to $\tilde f$ satisfies 
(\ref{eq:toj12}).  We claim that such triple is given by
\begin{equation}\label{eq:vhr}
\tilde v_i=v_i-\left<\gamma_i,T^t\Omega^{-1}\varphi\right>, \ \tilde h_{ij}=h_{ij}-\left<\gamma_j,T^t\Omega^{-1}\gamma_i\right
> \ \mbox{and} \ \tilde \rho_i=\rho_i-\sqrt{|c|}\left<T\varphi,\Omega^{-1}\gamma_i\right>,
\end{equation}
where  $T=-P^t-(P^t)^{-1}$. In fact, using (\ref{eq:phiimp}) 
and the fact that  $\xi_i=\phi f_*X_i$, $1\leq i\leq n$, we obtain from  (\ref{eq:tlp}) and (\ref{eq:betahor}) that
%\begin{equation}\label{eq:relac1nf}
$$
\beta^i=\beta^t\xi_i=
P^t\omega f_*^t\phi^t \xi_i=P^t\omega f_*^t\phi^t\phi f_*X_i=P^t\omega(X_i)=P^t\gamma_i=-L^tT\gamma_i
$$
%\end{equation}
where $L=(P^2+I)^{-1}P^2$, and 
$\beta^{n+1}=\beta^t(\xi_f)=\sqrt{|c|}P^t\varphi.
$

 Thus equations (\ref{eq:hbetatrans}) for $1\leq r\leq n$ reduce to the first two equations in  (\ref{eq:vhr}). Now, we have 
$
\tilde \xi_f={\cal P}(\xi_f+\sqrt{|c|}\sum_r\left<T\varphi,\O^{-1}\beta_r\right>\xi_r).
$ 
 Using   (\ref{eq:combe0}), (\ref{eq:o1}), (\ref{eq:ncon}), (\ref{eq:omegabeta}),  (\ref{formadcovnov}), (\ref{eq:pdef2})  and (\ref{ral_{h}}), the last of equations (\ref{eq:vhr}) follows from
$$
\tilde \nabla^\bot_{\partial_i}\tilde \xi_f=\tilde \rho_i{\cal P}(\xi_i-\sum_r\left<\beta^i,L^{-1}\O^{-1}\beta^r\right>\xi_r)=\tilde \rho_i\tilde \xi_i.
$$
Since $T\Omega$ is symmetric by Theorem \ref{teo:allagr}, it follows that 
$(\tilde v,\tilde h,\tilde \rho)$  satisfies 
(\ref{eq:toj12}).\vspace{1ex}\qed

We now prove the existence of $P$-transforms of any horizontal isometric immersion $f\colon\,M^n(c)\to \mathbb S^{2n+1}(c)$ with $\nu_f\equiv 0$ satisfying suitable initial conditions.

\begin{proposition}\label{lem:chavenf}\po Let  $f\colon\,M^n(c)\to \mathbb S^{2n+1}(c)$ be a horizontal isometric immersion with $\nu_f\equiv 0$ and let $P$ be an endomorphism of a Euclidean vector space $V$ such that $\sigma(P)\cap(-\sigma(P))=\emptyset$.
Fixed  $x_0\in M^n$ and  $\varphi_0\in V$,   $\omega_0\in T_{x_0}M^*\otimes V$ such that $\ker (P^c-\alpha I)\cap \ker (\varphi_0^t)^c\cap \ker (\omega_0^t)^c=\{0\}$,  there exist an open neighborhood $U$ of $x_0$ and a unique $P$-transform $\tilde f={\cal R}_{\va,\beta, \Omega, P}(f|_U)$ of $f|_U$ such that $\varphi(x_0)=\varphi_0$ and $d\varphi(x_0)=\omega_0$.
\end{proposition}
\proof Let $(u_1, \ldots, u_n)$ be principal  coordinates given by Corollary \ref{cor:toj7} on an open simply-connected neighborhood $U\subset M^n(c)$ of $x_0$ and   let $(v,h)$ be the solution of (\ref{eq:iii2}) associated to $f$. It is easily checked that the compatibility conditions of (\ref{ral_{h}}) are satisfied by virtue of (\ref{eq:iii2}). Therefore, if $\gamma^0_i=v_i^{-1}(x_0)\omega_0(\partial_i(x_0))$ for $1\leq i\leq n$, there exists a unique solution $(\varphi, \gamma_1, \ldots, \gamma_n)$ of (\ref{ral_{h}}) such that $\varphi(x_0)=\varphi_0$ and $\gamma_i(x_0)=\gamma_i^0$ for all  $1\leq i\leq n$.  By 
Proposition \ref{cor:condintb}, $\omega=d\varphi$ satisfies (\ref{eq:omega1b}) and $\omega(\partial_i)=v_i\gamma_i$ for $1\leq i\leq n$.

%Let $\varphi: U\rightarrow V$ be a solution of (\ref{eq:omega1}) on an open simply connected subset $U$ of $x_0\in M$, such that $\varphi(x_0)=\varphi_0$ and $d\varphi(x_0)=\omega_0$. 
By the assumption on $(\varphi_0,\omega_0)$, it follows from Proposition \ref{prop:padmi} that the unique solution $X=\Omega^0$ of the Lyapunov equation $X^tP+P^tX^t+T\rho_0=0$, 
with  $T=-P^t-(P^t)^{-1}$ and $\rho_0=P^t\omega_0\omega_0^tP+cP^t\varphi_0\varphi_0^tP$,  is invertible. Moreover,  $2\Omega^0_s={\cal G}^t(x_0){\cal G}(x_0)$ by  Proposition~\ref{teo:allagr}.  Now define $\beta\in \Gamma(V^*\otimes N_fM)$ by
(\ref{eq:betahor}), and let $\Omega$ be the unique solution of (\ref{eq:o1b}) on $U$ such that $\Omega(x_0)=\Omega^0$. Shrinking $U$ if necessary, we may assume that $\Omega$ is invertible on $U$.
By (\ref{eq:o1c}) and (\ref{eq:pdef2}) we have
$$
\begin{array}{l}
d(\Omega^tP+P^t\Omega^t + T\rho)(X)= \Phi(X)^t\omega^tP+P^t\Phi(X)^t\omega^t-((P^t)^2+I)d\omega(X)\omega^tP\vspace{0.1cm}\\
\hspace{10.ex}-((P^t)^2+I)\omega d\omega^t(X)P-c((P^t)^2+I)\omega(X)\varphi^tP-c((P^t)^2+I)\varphi\omega(X)^tP=0.
\end{array}
$$
%$$
%\begin{array}{l}
%d(\Omega^tP+P^t\Omega^t + T\rho)(X)= \Phi(X)^t\omega^tP+P^t\Phi(X)^t\omega^t-(P^t)^2d\omega(X)\omega^tP-(P^t)^2\omega d\omega^t(X)P\vspace{0.1cm}\\
%\hspace{10.ex}-d\omega(X)\omega^t P - \omega d\omega(X)^tP - c (P^t)^2\omega(X)\varphi^tP -c (P^t)^2\varphi\omega(X)^tP\vspace{0.1cm}\\
%\hspace{10.ex}- c \omega(X)\varphi^tP -c \varphi\omega(X)^tP\vspace{0.1cm}\\
%\hspace{10.ex}= ((P^t)^2+I)d\omega(X)\omega^tP+c((P^t)^2+I)\varphi X^t\omega^tP+P^t((P^t)^2+I)d\omega(X)\omega^t\vspace{0.1cm}\\
%\hspace{10.ex}+cP^t((P^t)^2+I)\varphi X^t\omega^t-(P^t)^2d\omega(X)\omega^tP-(P^t)^2\omega d\omega^t(X)P-d\omega(X)\omega^t P %\vspace{0.1cm}\\
%\hspace{10.ex}- \omega d\omega(X)^tP - c (P^t)^2\omega(X)\varphi^tP -c (P^t)^2\varphi\omega(X)^tP- c \omega(X)\varphi^tP -c %\varphi\omega(X)^tP
%\vspace{0.1cm}\\
%\hspace{10.ex}
%=0.
%\end{array}
%$$
Since (\ref{eq:lyaplagnf}) holds at $x_0$,  it holds  on $U$,  hence $\tilde f={\cal
R}_{\va,\beta,\Omega}(f|_U)$ is a $P$-transform of $f|_U$.\vspace{1ex}
\qed

\begin{remark}\po \emph{ Proposition \ref{lem:chavenf} also holds for horizontal isometric immersions $f\colon\,M^n(c)\to \mathbb S_\epsilon^{2n+1}(c)$ with $c<0$ and $\epsilon=-1$ if $\varphi_0\in V$ and  $\omega_0\in T_{x_0}M^*\otimes V$ are chosen so that the unique solution $X=\Omega^0$ of the Lyapunov equation $X^tP+P^tX^t+T\rho_0=0$, with  $T=-P^t-(P^t)^{-1}$ and $\rho_0=P^t\omega_0\omega_0^tP+cP^t\varphi_0\varphi_0^tP$,  is invertible. That such a pair $(\varphi_0, \omega_0)$ does exist was shown in part
$(ii)$ of Proposition \ref{prop:padmi}.}
\end{remark}

In the following corollary we summarize the process given by the  $P$-transformation to generate a family of new horizontal isometric immersions $f\colon\,M^n(c)\to \mathbb S_\epsilon^{2n+1}(c)$ starting with a given one and a vector-valued solution of a linear system of PDE's.

\begin{corollary}\label{cor:exem}\po
Let $f\colon M^n(c)\rightarrow \mathbb{S}_\epsilon^{2n+1}(c)$ be a horizontal isometric immersion with $\nu_f\equiv 0$, let $(u_1, \ldots, u_n)$ be local principal coordinates given by Corollary \ref{cor:toj7}, let $(v,h)$ the solution of (\ref{eq:iii2}) associated to $f$, let $P$ be an endomorphism of a Euclidean vector space $V$ such that $\sigma(P)\cap (-\sigma(P))=\emptyset$ and let $(\varphi,\gamma_1,...,\gamma_n)$ be a $V$-valued solution of the linear system of PDE's  
\begin{equation}\label{Ral_2}
\left\{
\begin{array}{l}
i) \ \partial_i\varphi=v_i\gamma_i,\ \  (ii) \ \partial_i\gamma_j= h_{ji}\gamma_i, \ 1\leq i\neq j\leq n,\vspace{0.1cm}\\
iii)  \ \partial_i\gamma_i + \sum_{j\neq i}h_{ji} \gamma_j+(P^t)^{-1}\gamma_i+cv_i\varphi = 0,\vspace{0.1cm}\\
\end{array}
\right.
\end{equation}
on an open subset $W$ where 
%\begin{equation}\label{eq:vtilpeh}
$$
\tilde v_i= v_i+\left<\gamma_i,(P+P^{-1})\Omega^{-1}\varphi\right>
$$
%\end{equation}
does not vanish for  $1\leq i\leq n$, with 
$\Omega^t$ given by (\ref{eq:solx}) for $C=(({P^t})^2+I)\left(d\varphi(d\varphi)^t +c\varphi\varphi^t\right)P$.
% and $A=P$.
%$$A=P\,\,\,\mbox{and}\,\,\,C=({P^t}^2+I)\left(d\varphi(d\varphi)^t +c\varphi\varphi^t\right)P.$$
 Then $\tilde f:W\rightarrow \mathbb S_\epsilon^{2n+1}(c)$ given by
$$ %\begin{equation*}\label{eq:partilf}
\begin{array}{l}
\tilde F =i\circ \tilde f=F-\sum_{j=1}^n\left(\left<\Omega^{-1}\varphi,\gamma_j\right>+\left<P\Omega^{-1}\varphi,\gamma_j\right>i\right)F_*X_j\vspace{0.1cm}\\
\hspace{12.ex}-\left(\left<\Omega^{-1}\varphi,\varphi\right>+\left<P\Omega^{-1}\varphi,\varphi\right>i\right)cF,\vspace{0.1cm}\\
\end{array}
$$%\end{equation*}
where $X_j=v_j^{-1}\partial_j$ for $1\leq j\leq n$,  $i:\mathbb S_\epsilon^{2n+1}(c)\rightarrow \mathbb C^{n+1}_\epsilon$ is the umbilical inclusion and $F=i\circ f$, defines a new horizontal isometric immersion whose  induced metric has constant  curvature $c$.
Moreover, the solution   of (\ref{eq:iii2}) associated to $\tilde f$ is $(\tilde v,\tilde h)$, with
%\begin{equation}\label{eq:htilpeh}
$$
\tilde{h}_{ij}= h_{ij}+\left<\gamma_j,(P+P^{-1})\Omega^{-1}\gamma_i\right>.
$$
%\end{equation}
%In particular, $(\tilde v,\tilde h)$ is a new solution  of (\ref{eq:iii2}).
\end{corollary}

%\begin{corollary}
%Let $(v,h)$ be a solution of (\ref{eq:iii2}) and let $(\varphi,\gamma)$ be a solution of (\ref{ral_{h}}). Then a new solution $(\tilde v,\tilde h)$ of (\ref{eq:iii2}) is given by Corollary \ref{cor:exem}.
%\end{corollary}

%\begin{remark}\po \emph{As was done in \cite{t} in the scalar case, one can  produce new explicit examples of $n$-dimensional flat Lagrangian submanifolds of  $\mathbb{C}^n=\R^{2n}$ and $n$-dimensional  submanifolds with constant sectional curvature $c$ of  $\Sf_\epsilon^{2n+1}(c)$ that are horizontal with respect to the Hopf fibration $\pi\colon \Sf_\epsilon^{2n+1}(c)\to \tilde M^n(4c)$ by applying Corollaries \ref{cor:exlag} and \ref{cor:exem}, respectively,  to the (degenerate) 
%examples that correspond to the trivial solution of systems (\ref{eq:corlagsist}) and (\ref{eq:iii2}) given by 
%\begin{equation}\label{eq:t35}
%v_1=b\neq 0, \ v_i=0, \ 2\leq i\leq n, \ \mbox{and} \ h=0.
%\end{equation}
%$v_1=b\neq 0$, $v_i=0$, $2\leq i\leq n$, and $h=0$. In fact, in this case systems (\ref{ral_{0l}}) and (\ref{Ral_2}) become linear systems of  ODE's which can be explicitly solved. We omit the lengthy and straightforward computations.}
%\end{remark}

\section{A decomposition theorem for the $P$-transformation} 

In this section we prove the following  decomposition theorem for the  $P$-transformation, for which Remarks \ref{re:ldecomp} also apply.

 \begin{theorem}\label{thm:permlagran}\po Let $f\colon M^n(0)\to \R^{2n}$ (respectively,  $f\colon M^n(c)\to \mathbb S^{2n+1}(c)$) be a Lagrangian (respectively, horizontal) isometric immersion and let
 $\tilde f={\cal R}_{\va,\beta, \Omega,P}(f)\colon \tilde{M}^n(0)\to \R^{2n}$ (respectively, ${\cal R}_{\va,\beta, \Omega,P}(f)\colon \tilde{M}^n(c)\to \mathbb S^{2n+1}(c)$)
 be a $P$-transform of $f$ such that $P=P_1\oplus P_2$ with respect to an orthogonal decomposition $V=
V_1\oplus V_2$. Define $\varphi_j, \beta_j, \Omega_{ij}$  by (\ref{vabo2nf}), and $\bar \va_i, \bar \beta_i, \bar\Omega_{ii}$  by (\ref{eq:bar2nf}), $1\leq i,j\leq 2$.
%$\varphi_j=\pi_{V_j}\circ \varphi$, $1\leq j\leq 2$, and
%$$
%\bar \va_j=\varphi_j-\Omega_{ji}\Omega_{ii}^{-1}\varphi_i,
%$$
%where $\Omega_{ij}=\pi_{V_i}\circ \Omega|_{V_j}$.  
Then $(\varphi_j,\beta_j, \Omega_{jj},P_j)$ defines a $P_j$-transform $f_j$ of $f$ for $1\leq j\leq 2$, $(\bar \varphi_i,\bar \beta_i, \bar \Omega_{ii}, P_i)$ a $P_i$-transform $f_{ij}$ of $f_j$ for $1\leq i\neq j\leq 2$ and (\ref{eq:permnf2}) holds.
If, in addition, $f(M^n(0))\subset \Sf^{2n-1}$ and $\tilde f$ is a $P^*$-transform of $f$, then $f_j$  is a $P^*_j$-transform  of $f$ 
and $f_{ji}$ is a $P_j^*$-transform of $f_i$ for $1\leq i\neq  j\leq 2$.
%$$
%{\ral}_{\varphi,\beta, \Omega,P}(f)={\ral}_{\bar \varphi_i,\bar\beta_i, \bar\Omega_{ii},P_i}(\ral_{\varphi_j,\beta_j, \Omega_{jj}, P_j}(f)).
%$$
\end{theorem}

\proof We give the proof for horizontal isometric immersions $f\colon M^n(c)\to \mathbb S^{2n+1}(c)$, the case of  Lagrangian isometric immersions $f\colon M^n(0)\to \R^{2n}$  being similar.  
By the assumption that  $\tilde f=\ral_{\varphi,\beta,\Omega,P}(f)$ is a $P$-transform of $f$, we have that $\beta $ and $\Omega$ satisfy (\ref{eq:betahor}) and 
(\ref{eq:lyaplagnf}), respectively. In particular, $\Omega$ is an invertible solution of the Lyapunov equation (\ref{eq:lyaplagnf}), hence  $\ker(P^c-\alpha I) \cap \ker (\varphi^t)^c\cap \ker (\omega^t)^c=\{0\}$ for any  eigenvalue $\alpha$ of $P$ by  Proposition~\ref{prop:padmi}. Therefore $\ker(P_j^c-\alpha I)\cap \ker (\varphi_j^t)^c\cap \ker (\omega_j^t)^c=\{0\}$ for $1\leq j\leq 2$,  and hence also $\Omega_{jj}$ is invertible by Proposition \ref{prop:padmi}. 
 By Theorem \ref{c:viiinf}, $(\varphi_j, \beta_j, \Omega_{jj})$ satisfies the conditions of Definition \ref{ribvetgeral} for $1\leq j\leq 2$. 
 %We now prove that $f_j=\Ral_{\varphi_j,\beta_j, \Omega_{jj}}(f)$ are $P_j$-transforms of $f$ for $1\leq j\leq 2$. 

Since $P=P_1\oplus P_2$, and hence  $T=T_1\oplus T_2$, equations (\ref{eq:betahor}) and 
(\ref{eq:lyaplagnf}) are equivalent to
$
 \beta_j=(\phi f_*\omega^t_j+\epsilon \sqrt{|c|}\xi_f\varphi_j^t)P_j
$
and
\be\label{eq:ijs}
 \Omega_{ij}^tP_i+P_j^t\Omega_{ij}^t+T_j\rho_{ji}=0, \ 1\leq i,j\leq 2,
\ee
where $\rho_{ji}=P_j^t\omega_j\omega^t_iP_i+cP_j^t\varphi_j\varphi_i^tP_i$. In particular,  $f_j=\Ral_{\varphi_j,\beta_j, \Omega_{jj}}(f)$ is a $P_j$-transform of $f$,  $1\leq j\leq 2$. 

Now by (\ref{eq:tomega2}) we have
\begin{equation}\label{eq:toott}
T_i\Omega_{ij}= \Omega_{ji}^tT_j^t, \,\,\,1\leq i,j\leq 2.
%\ \mbox{and} \  T_j\Omega_{jj}= \Omega_{jj}^{t}T_j^t.
\end{equation}
%Thus, it easily get that 
%$$
%\bar \Omega_{ii}^tP_i+P_i^t\bar \Omega_{ii}^t+T_i\bar \beta_i^t\bar \beta_i=0,
%$$
%Set $\bar \beta_i={\cal P}_j(\beta_i-\beta_j(\Omega_{jj}^{-1})^t\Omega_{ij}^t)$, $1\leq i\neq j\leq 2$. 
Using 
%(\ref{eq:tomega2}) 
this  we obtain
$$
\begin{array}{lll}
\bar{\Omega}_{ii}^tP_i+P_i^t\bar{\Omega}_{ii}^t+T_i\bar \beta_i^t\bar \beta_i&=& \left(\Omega_{ii}^t-\Omega_{ji}^t(\Omega_{jj}^{-1})^t\Omega_{ij}^t\right)P_i+P_i^t\left(\Omega_{ii}^t-\Omega_{ji}^t(\Omega_{jj}^{-1})^t\Omega_{ij}^t\right)\vspace{0.1cm}\\&&+T_i\left(\beta^t_i-\Omega_{ij}\Omega_{jj}^{-1}\beta_j^t\right)\left(\beta_i-\beta_j(\Omega_{jj}^t)^{-1}\Omega_{ij}^t\right)\vspace{0.1cm}\\

&=&\left(\Omega_{ii}^tP_i+P_i^t\Omega_{ii}^t+T_i\beta_i^t\beta_i\right)-\Omega_{ji}^t(\Omega_{jj}^{-1})^t\Omega_{ij}^tP_i-P_i^t\Omega_{ji}^t(\Omega_{jj}^{-1})^t\Omega_{ij}^t\vspace{0.1cm}\\

&&-T_i\beta_i^t\beta_j(\Omega_{jj}^{-1})^t\Omega_{ij}^t-T_i\Omega_{ij}\Omega_{jj}^{-1}\beta_j^t\beta_i+T_i\Omega_{ij}\Omega_{jj}^{-1}\beta_j^t\beta_j(\Omega_{jj}^{-1})^t\Omega_{ij}^t\vspace{0.1cm}\\

%
%&=&-\Omega_{ji}^t(\Omega_{jj}^{-1})^t\Omega_{ij}^tP_i-\left(P_i^t\Omega_{ji}^t+T_i\beta_i^t\beta_j\right)(\Omega_{jj}^{-1})^t\Omega_{ij}^t\vspace{0.1cm}\\
%
%&&-\Omega_{ji}^t(\Omega_{jj}^{-1})^tT_j\beta_j^t\beta_i+\Omega_{ji}^t(\Omega_{jj}^{-1})^tT_j\beta_j^t\beta_j(\Omega_{jj}^{-1})^t\Omega_{ij}^t\vspace{0.1cm}\\
%
&=&-\Omega_{ji}^t(\Omega_{jj}^{-1})^t\left(\Omega_{ij}^tP_i+T_j\beta_j^t\beta_i\right)-\left(P_i^t\Omega_{ji}^t+T_i\beta_i^t\beta_j\right)(\Omega_{jj}^{-1})^t\Omega_{ij}^t\vspace{0.1cm}\\

&&+\Omega_{ji}^t(\Omega_{jj}^{-1})^tT_j\beta_j^t\beta_j(\Omega_{jj}^{-1})^t\Omega_{ij}^t\vspace{0.1cm}\\

%&=&\Omega_{ji}^t(\Omega_{jj}^{-1})^tP_j^t\Omega_{ij}^t+\Omega_{ji}^tP_j(\Omega_{jj}^{-1})^t\Omega_{ij}^t+\Omega_{ji}^t(\Omega_{jj}^{-1})^tT_j\beta_j^t\beta_j(\Omega_{jj}^{-1})^t\Omega_{ij}^t\vspace{0.1cm}\\
&=&\Omega_{ji}^t\left((\Omega_{jj}^{-1})^tP_j^t+P_j(\Omega_{jj}^{-1})^t+(\Omega_{jj}^{-1})^tT_j\beta_j^t\beta_j(\Omega_{jj}^{-1})^t\right)\Omega_{ij}^t=0.\vspace{0.1cm}\\
%&=&0.\\
\end{array}
$$
Now let $\{X_{1},...,X_{n}\}$ be an orthonormal principal frame for $f$ and define $\xi_1, \ldots, \xi_n$ and $X_{j,1},...,X_{j,n}$ by $\xi_l=\phi_*f_*X_l$ and 
  $D_jX_{j,l}=X_l$, $1\leq l\leq n$, $1\leq j\leq 2$.  
Denote $\omega_j(X_{l})=\gamma_{j,l}$ and $\beta_j^t(\xi_{l})=\beta_j^l$,  and define 
 $\bar{\gamma}_{i,l}=\bar{\omega}_iX_{j,l}$ and $\bar \beta_i^l=\bar \beta_i^t(\xi_{j,l})$,  where
$\xi_{j,l}={\cal P}_j(\xi_l-\beta_j(\Omega_{jj}^{-1})^t(L_j^{-1})^t\beta_j^l)$.
It follows from (\ref{eq:baromega2}) that
$$
\bar{\gamma}_{i,l}=\omega_iD_jX_{j,l}-\O_{ij}\O_{jj}^{-1}\omega_jD_jX_{j,l}=\gamma_{i,l}-\O_{ij}\O_{jj}^{-1}\gamma_{j,l}.
$$
On the other hand, from (\ref{eq:omegabeta}),  
(\ref{eq:lyaplagnf}) and (\ref{eq:toott}) we obtain
$$
\begin{array}{lll}
\bar \beta_i^l&=& \left(\beta_i^t-\Omega_{ij}\Omega_{jj}^{-1}\beta_j^t\right){\cal P}_j^t{\cal P}_j(\xi_l-\beta_j(\Omega_{jj}^{-1})^t(L_j^{-1})^t\beta_j^l)\vspace{0.1cm}\\
&=&\beta_i^l-\rho_{ij}(\Omega_{jj}^{-1})^t(L_j^{-1})^t\beta_j^l-\Omega_{ij}\Omega_{jj}^{-1}\beta_j^l+\Omega_{ij}\Omega_{jj}^{-1}\rho_{jj}(\Omega_{jj}^t)^{-1}(L_j^{-1})^t\beta_j^l\vspace{0.1cm}\\
&=&P_i^t\gamma_{i,l}+\rho_{ij}(\Omega_{jj}^t)^{-1}T_j\gamma_{j,l}-\Omega_{ij}\left(\Omega_{jj}^{-1}P_j^t+\Omega_{jj}^{-1}\rho_{jj}T_j^t\Omega_{jj}^{-1}\right)\gamma_{j,l}\vspace{0.1cm}\\
&=&P_i^t\gamma_{i,l}+\rho_{ij}T_j^t\Omega_{jj}^{-1}\gamma_{j,l}+\Omega_{ij}P_j\Omega_{jj}^{-1}\gamma_{j,l}
=P_i^t\gamma_{i,l}+\left(\rho_{ij}T_j^t+\Omega_{ij}P_j\right)\Omega_{jj}^{-1}\gamma_{j,l}\vspace{0.1cm}\\
&=&P_i^t\gamma_{i,l}-P_i^t\Omega_{ij}\Omega_{jj}^{-1}\gamma_{j,l}=P_i^t\bar \gamma_{i,l}.
\end{array}
$$
It follows that $(\bar \varphi_i, \bar \beta_i, \bar\Omega_{ii}, P_i)$  defines a  $P_i$-transform of $f_j$ for $1\leq i\neq j\leq 2$.

We now prove the last asssertion. If $\tilde f$  is  a $P^*$-transform of $f$, then  $\varphi+\sum_{l=1}^n v_lP^t\gamma_l=0$.
Thus  $\varphi_j+ \sum_{l=1}^n v_lP_j^t\gamma_{j,l}=0$ for $1\leq j\leq 2$, which already shows that  $f_j$ is a $P_j^*$-transform of $f$. To prove that $f_{ji}$ is a $P_j^*$-transform of $f_i$ for $1\leq j\leq 2$, we must show that
 $$\bar{\varphi}_j+\sum_{l=1}^n v_{i,l}P_j^t\bar{\gamma}_{j,l}=\bar \varphi_j+P_j^t\bar \omega_j\sum_{l=1}^n\partial_l=0.$$
Using (\ref{eq:ijs}) we obtain
\begin{equation}\label{eq:inipioi}
\begin{array}{lll}
P_j^t\bar \omega_j&=&P_j^t(\omega_jD_i-\Omega_{ji}\Omega_{ii}^{-1}\omega_iD_i)
\vspace{0.1cm}\\
&=& P_j^t\omega_jD_i- \Omega_{ji}\Omega_{ii}^{-1}P_i^t\omega_iD_i-\Omega_{ji}\Omega_{ii}^{-1}\rho_{ii}T_i^t\Omega_{ii}^{-1}\omega_iD_i 
+\rho_{ji}T_i^t\Omega_{ii}^{-1}\omega_iD_i. 
\end{array}
\end{equation}
On the other hand, it follows from (\ref{eq:dembase}), (\ref{eq:vtilpe}) and $\omega_i^t(v)=\sum_{l=1}^n\left<\gamma_{i,l},v\right>X_l$ that 
\begin{equation}\label{eq:di1}
\begin{array}{lll}
\sum_{l=1}^n\omega_zD_i\partial_l&=&\sum_{l=1}^n v_{i,l}\gamma_{z,l}\vspace{0.1cm}\\
&=&\sum_{l=1}^n (v_l-\<\gamma_{i,l}, T_i^t\Omega_{ii}^{-1}\varphi_i\>)\gamma_{z,l}\vspace{0.1cm}\\
&=&\sum_{l=1}^n v_l\gamma_{z,l}-\omega_z\omega_i^tT_i^t\Omega_{ii}^{-1}\varphi_i\vspace{0.1cm}\\
&=&\sum_{l=1}^n v_l\gamma_{z,l}-(P_z^t)^{-1}\rho_{zi}P_i^{-1} T_i^t\Omega^{-1}_{ii}\varphi_i,\ \ \ z=i,j.\vspace{0.1cm}\\
\end{array}
\end{equation}
Thus, using  (\ref{eq:inipioi}) and (\ref{eq:di1})  we obtain
$$
\begin{array}{l}
\bar \varphi_j+P_j^t\bar \omega_j\sum_{l=1}^n\partial_l=\varphi_j-\Omega_{ji}\Omega_{ii}^{-1}\varphi_i+P_j^t\sum_{l=1}^n v_l\gamma_{j,l}-
\rho_{ji}P_i^{-1}T_i^t\Omega_{ii}^{-1}\varphi_i\vspace{1ex}\\
\hspace*{18ex}-\Omega_{ji}\Omega_{ii}^{-1}P_i^t\sum_{l=1}^n v_l\gamma_{i,l}+\Omega_{ji}\Omega_{ii}^{-1}\rho_{ii}P_i^{-1}T_i^t\Omega_{ii}^{-1}\varphi_i\vspace{1ex}\\
\hspace*{18ex}-
\Omega_{ji}\Omega_{ii}^{-1}\rho_{ii}T_i^t\Omega_{ii}^{-1}\sum_{l=1}^n v_l\gamma_{i,l} +\Omega_{ji}\Omega_{ii}^{-1}\rho_{ii}T_i^t\Omega_{ii}^{-1}(P_i^t)^{-1}\rho_{ii}P_i^{-1}T_i^t\Omega_{ii}^{-1}\varphi_i
\vspace{1ex}\\
\hspace*{18ex}+\rho_{ji}T_i^t\Omega_{ii}^{-1}\sum_{l=1}^n v_l\gamma_{i,l}-\rho_{ji}T_i^t\Omega_{ii}^{-1}(P_i^t)^{-1}\rho_{ii}P_i^{-1}T_i^t\Omega_{ii}^{-1}\varphi_i
\vspace{1ex}\\
\hspace*{17ex}= (\rho_{ji}-\Omega_{ji}\Omega_{ii}^{-1}\rho_{ii})\left(P_i^{-1}T_i^t\Omega_{ii}^{-1}P_i^t+T_i^t\Omega_{ii}^{-1}\right.\vspace{1ex}\\
\hspace*{18ex}\left.
+T_i^t\Omega_{ii}^{-1}(P_i^t)^{-1}\rho_{ii}P_i^{-1}T_i^t\Omega_{ii}^{-1}P_i^t\right)\sum_{l=1}^n v_l\gamma_{i,l}=0,
\end{array} 
$$
where  we have used that 
$\rho_{ii}=-P_i^t\Omega_{ii}(T_i^t)^{-1}-\Omega_{ii}P_i(T_i^t)^{-1}$, as follows from (\ref{eq:ijs}). \qed

\subsection{The Bianchi $P$-cube.}

Given $P_1,P_2\in \mathbb R$ with $P_1\neq \pm P_2$, we say that a Bianchi quadrilateral  $\{f,f_1,f_2,f_{12}\}$  is a
 \emph{Bianchi $(P_1,P_2)$-quadrilateral} of $n$-dimensional horizontal submanifolds with constant curvature $c$ of $\mathbb S^{2n+1}(c)$ if  $f\colon M^n(c)\to \mathbb S^{2n+1}(c)$ is horizontal, 
 $f_i$ is a $P_i$-transform of $f$ for $1\leq i\leq 2$ and $f_{12}$ is a $P_2$-transform of $f_1$ and a $P_1$-transform of $f_2$.
 Similarly one defines a Bianchi $(P_1,P_2)$-quadrilateral of $n$-dimensional flat Lagrangian submanifolds of $\R^{2n}$, as well as of   $n$-dimensional flat Lagrangian submanifolds of $\R^{2n}$ that are contained in  $\Sf^{2n-1}\subset \R^{2n}$, in the latter case requiring  $f_i$ to be  a $P^*_i$-transform of $f$ for $1\leq i\leq 2$ and $f_{12}$ to be  a $P^*_2$-transform of $f_1$ and a $P^*_1$-transform of $f_2$. Below we state and prove a Bianchi-cube theorem for the 
 first of these classes, analogous results being true for the others.
We will need the following result of  \cite{dt3} (see also  \cite{t}).

\begin{proposition}\label{prop:permpl}\po
If  $f_i={\cal R}_{\varphi_i, \beta_i, \Omega_i, P_i}(f):M_i^n(c) \to \mathbb S^{2n+1}(c)$, $1\leq i\leq 2$, are  $P_i$-transforms of  $f\colon M^n(c)\to \mathbb S^{2n+1}(c)$  with $P_1\neq \pm P_2$ and  $[d\omega_1^t,d\omega_2^t]=0$, then there is a unique isometric immersion $\tilde f\colon \tilde{M}^n(c)\to \mathbb S^{2n+1}(c)$ such that $\{f,f_1,f_2,\tilde f\}$ is a Bianchi $(P_1,P_2)$-quadrilateral.
\end{proposition}

Given $P_1,...,P_k\in \mathbb R$, with $P_i\neq \pm P_j$ for all $1\leq i\neq j\leq k$, we say that a Bianchi cube
$({\cal{C}}_0, . . . ,{\cal{C}}_k)$ is a \textit{Bianchi $(P_1,...,P_k)$-cube} of $n$-dimensional horizontal submanifolds with constant curvature $c$ of $\mathbb S^{2n+1}(c)$ if  $f\colon M^n(c)\to \mathbb S^{2n+1}(c)$  is  horizontal  and, for all $1\leq s\leq k-1$, 
\begin{itemize}
	\item[(i)] Each $f_{\alpha_{s+1}}\in {\cal C}_{s+1}$ with  $\alpha_{s+1}=\alpha_{s}\cup \{i_j\}$  is a $P_{i_j}$-transform of $f_{\alpha_{s}}\in {\cal C}_s$.
	\item[(ii)]  $\{f_{\alpha_{s-1}},f_{\alpha_{s-1}\cup \{i_l\}},f_{\alpha_{s-1}\cup \{i_j\}},f_{\alpha_{s+1}}\}$ is a Bianchi $(P_{i_l},P_{i_j})$-quadrilateral when  $\alpha_{s+1}=\alpha_{s-1}\cup \{i_l,i_j\}$.
\end{itemize}

\begin{theorem}\po
Let $f\colon M^n(c)\to \mathbb S^{2n+1}(c)$ be a horizontal isometric immersion with $\nu_f\equiv~0$ and let  $f_i=\ral_{\varphi_i,\beta_i, \Omega_i, P_i}(f)\colon M^n(c)\to \mathbb S^{2n+1}(c)$, $1\leq i\leq k$, be  $P_i-$transforms of $f$  such that $P_i\neq \pm P_j$ and  $[d\omega_i^t,d\omega_j^t]=0$ for all $1\leq  i\neq j\leq k$. Then there exists a unique  Bianchi $(P_1, \ldots, P_k)$-cube $({\cal C}_0,...,{\cal C}_k)$ such that ${\cal C}_0=\{f\}$ and ${\cal C}_1=\{f_1,...,f_k\}$.
%, which is unique if no $f_i$ belongs to the associated family determined by $\{f_j,f_l\}$ for all $1\leq i\neq j\neq l\neq i\leq k$.
\end{theorem}
%\begin{theorem}
%If $f_i=\ral_{\varphi_i,\beta_i, \Omega_i, P_i}(f)\colon M^n(c)\to \mathbb S^{2n+1}(c)$, $1\leq i\leq k$, are $P_i-$transforms of a horizontal isometric immersion  $f\colon M^n(c)\to \mathbb S^{2n+1}(c)$ with $\nu_f\equiv~0$ such that $P_i\neq \pm P_j$ and  $[d\omega_i^t,d\omega_j^t]=0$ for all $1\leq  i\neq j\leq k$,  then there exists a unique  Bianchi $(P_1, \ldots, P_k)$-cube $({\cal C}_0,...,{\cal C}_k)$ such that ${\cal C}_0=\{f\}$ and ${\cal C}_1=\{f_1,...,f_k\}$.
%, which is unique if no $f_i$ belongs to the associated family determined by $\{f_j,f_l\}$ for all $1\leq i\neq j\neq l\neq i\leq k$.
%\end{theorem}
\proof We first prove existence. Since $f_i$ is a $P_i$-transform of $f$, we have 
\begin{equation}\label{eq:betigeral}
\beta_i=P_i(\phi f_*\nabla \varphi_i+\epsilon \sqrt{|c|}\varphi_i\xi_f), \,\,\,1\leq i\leq k.
\end{equation}
Define $\varphi\colon M^n(c)\to \R^k$,  $\beta\in \Gamma((\R^k)^*\otimes N_fM)$ and $\Omega\in \Gamma((\R^k)^*\otimes \R^k)$ by  
$$\varphi=(\varphi_1, \ldots, \varphi_k),\,\,\,\,\,\beta=\sum_{i=1}^k e^i\otimes \beta_i \,\,\,\mbox{and}\,\,\,\,\Omega=\sum \Omega_{ij} e^j\otimes e_i$$
where $e_1, \ldots, e_k$ is the canonical basis of $\R^k$,  $e^1, \ldots, e^k$ is its dual basis, and
\be\label{eq:omegaijb}
\Omega_{ij}=\frac{(1+P_j^2)}{P_j(P_i+P_j)}\rho_{ij}, \ 1\leq i,j\leq k,
\ee
with $\rho_{ij}=P_iP_j\left<\nabla\varphi_i,\nabla\varphi_j\right>+cP_iP_j\varphi_i\varphi_j$.  We claim that $(\varphi, \beta, \Omega)$ defines a $P$-transformation of $f$, where $P\in (\R^k)^*\otimes \R^k$ is given by
$P=\sum_{i=1}^kP_ie^i\otimes e_i$. First notice that $f_i$ is an $L_i$-transform of $f$ for  $1\leq i\leq k$, with $L_i=\frac{P_i^2}{P_i^2+1}$, and that 
$$
\Omega_{ij}=\frac{-L_i\left<{\cal G}_i,{\cal G}_j\right>+\rho_{ij}}{L_j-L_i}
$$
where ${\cal G}_i=F_*\nabla \varphi_i+ i_*\beta_i+c\varphi_iF$, with $F=i\circ f$.  It follows from the proof of Theorem~\ref{lcubo} that  
 $(\varphi, \beta, \Omega)$ defines an $L$-transformation of $f$, where $L=(P^2+I)^{-1}P^2$.   Equation (\ref{eq:betahor}) follows immediately from (\ref{eq:betigeral}). It remains to  prove (\ref{eq:lyaplagnf}), 
 which is equivalent to  
$$
\Omega_{ij}(P_i+P_j)+T_j\rho_{ij}=0, \ 1\leq  i, j\leq k.
$$
where $T_j=-(P_j^2+1)/P_j$. For $i=j$ this  follows from the fact that  $f_i$ is a $P_i$-transform of $f$, while for $i\neq j$ this is a consequence of  (\ref{eq:omegaijb}).  
 
 Now for any multi-index $\alpha_r=\{i_1,...,i_r\}\in \Lambda_r$, define $(\varphi^{\alpha_r}, \beta^{\alpha_r}, \Omega^{\alpha_r})$ by (\ref{eq:varbetaalp}) and  set $P^{\alpha_r}=\sum_{\ell, j=1}^rP_{i_\ell i_j}e^{i_\ell}\otimes e_{i_j}$.
By Theorem \ref{thm:permlagran}, for each $\alpha_r\in \Lambda_r$ the triple $(\varphi^{\alpha_r}, \beta^{\alpha_r}, \Omega^{\alpha_r})$ defines a 
$P^{\alpha_r}$-transform of $f$. Let  ${\cal C}_r$ be the family of ${k \choose r}$ elements formed by the $P^{\alpha_r}$-transforms $f_{\alpha_r}=\ral_{\varphi^{\alpha_{r}},\beta^{\alpha_{r}},\Omega^{\alpha_{r}},P^{\alpha_{r}}}(f)$ of $f$.
It also follows from Theorem \ref{thm:permlagran} that conditions $(i)$ and $(ii)$ in the definition of a Bianchi $(P_1,...,P_k)$-cube are satisfied by  $({\cal{C}}_0, . . . ,{\cal{C}}_k)$.

Finally, since $({\cal{C}}_0, . . . ,{\cal{C}}_k)$ is also a Bianchi $(L_1,...,L_k)$-cube, the uniqueness assertion follows from that  in Theorem \ref{lcubo}, once we show that the last assumption in that result is unnecessary in this case. In fact, in view of Remark \ref{re:assfamily}, if $f_i$ belongs to the associated family determined by $\{f_j, f_l\}$,  $1\leq i\neq j\neq l\neq i\leq k$, then the vector field $\xi$ given by (\ref{vecxi}) satisfies $A_\xi=0$. Using  (\ref{eq:betahor}) and the first equation in (\ref{eq:sffhor}), this implies that $A_\eta=0$ for
$$\eta=\phi f_*(k_j(C_i-C_j)P_j\nabla\varphi_j+ k_l(C_i-C_l)P_l\nabla\varphi_l)\in \Gamma(\phi f_*TM).$$ But then the vector  field under brackets must be zero by the second equation in (\ref{eq:sffhor}) and the assumption that $\nu_f\equiv 0$, and this is a contradiction.  \qed

\section{Examples}

In this section we illustrate the procedures developed in the preceding sections. We produce examples of $n$-dimensional flat Lagrangian submanifolds of  $\mathbb{C}^n$ by applying the $P$-transformation to the trivial solution 
\be\label{eq:trivialsol}
v_1=1, \;\;\;v_i=0,\;\;2\leq i\leq n,\;\;\;\mbox{and}\;\;\;h_{ij}=0,\;\;1\leq i\neq j\leq n,
\ee
of the system of PDE's (\ref{eq:corlagsist}).

 Note that, given a solution  $(v,h)$ of (\ref{eq:corlagsist}), with $v=(v_1, \ldots, v_n)$ and $h=(h_{ij})$,   on an open  simply connected subset $U\subset \R^n$ where $v_i\neq 0$ everywhere for $1\leq i\leq k$, in order to determine the corresponding Lagrangian isometric immersion $f\colon\,U\to \mathbb{C}^n$  one has to integrate the system of PDE's
\be\label{lsyst}
\left\{  \begin{array}{l}
i)\; {\displaystyle \partial_i f =
v_iX_i},\;\;\;\;\; 
ii)\; {\displaystyle \partial_j X_{i} =
h_{ij}X_{j}},\;\; i\neq j,\;\;\;\;\; \vspace{1.5ex}\\
iii)\; {\displaystyle \partial_i X_{i} = -\sum_{k\neq i}h_{ki}X_k +iX_i,}
\end{array} \right. 
\ee
with initial conditions $(X_1(u_0),\ldots ,X_n(u_0))$ at some point $u_0\in U$ chosen so that
$$
\<X_i(u_0),X_j(u_0)\>=\<iX_i(u_0),X_j(u_0)\>=0,\;\;i\neq j,\;\;\mbox{and}\;\;\;\<X_i(u_0),X_i(u_0)\>=1.$$

For the trivial solution (\ref{eq:trivialsol}), system (\ref{lsyst}) becomes 
$$
\left\{  \begin{array}{l}
i)\; {\displaystyle \partial_1 f =
X_1},\;\;\;
ii)\; {\displaystyle \partial_i f =
0},\;\;\;2\leq i\leq n, \vspace{1.5ex}\\
iii)\; {\displaystyle \partial_j X_{i} =
0},\;\; 1\leq i\neq j\leq n,\;\;\;
iv)\; {\displaystyle \partial_i X_{i} = iX_i},\;\;\;1\leq i\leq n,
\end{array} \right. 
$$ 
whose solution, with initial conditions $(F(0), X_1(0), \ldots, X_n(0))=(-iE_1, E_1,, \ldots, E_n)$, where $E_1, \ldots, E_n$ is the canonical basis of 
$\mathbb{C}^n$ over $\mathbb{C}$, is 
\be\label{eq:degen}
f=-ie^{iu_1}E_1, \;\;\; X_j=e^{iu_j}E_j,
\;\;1\leq j\leq n.
\ee
System (\ref{ral_{0l}}) reduces in this case to
\begin{equation}\label{ral_{0l}b}
\left\{
\begin{array}{ll}
i) \ \partial_1\varphi=\gamma_1, \;\;\; \partial_i\varphi=0, \;\;2\leq i\leq n,\vspace{1ex}\\
ii) \ \partial_i\gamma_i=-(P^t)^{-1}\gamma_i, \;\;\; iii) \ \partial_i\gamma_j = 0,  \ j\neq i.
\end{array}
\right.
\end{equation}
Given $P_1, \ldots, P_k\in \R$, with $P_i\neq \pm P_j$ for $1\leq i\neq j\leq k$, the solution of (\ref{ral_{0l}b}) for the scalar $P_i$-transformation, $1\leq i\leq k$,  with initial conditions $\varphi_i(0)=-P_i$ and $\gamma^i_j(0)=1$, is 
\be\label{eq:vagamma}
\varphi_i=-P_ie^{-\frac{u_1}{P_{i}}}\;\;\;\mbox{and}\;\;\;\gamma^i_j=e^{-\frac{u_j}{P_{i}}},\;\; 1\leq i\leq k, \;\;1\leq j\leq n.
\ee
The corresponding $P_i$-transform of $f$ is 
$$ f_i=-ie^{iu_1}E_1+\frac{2P_i(1+P_i i)e^{-u_1/P_i}}{(1+P_i^2)\sum_{j=1}^ne^{-2u_j/P_i}}\sum_{j=1}^n e^{\frac{(-1+P_i i)u_j}{P_i}} E_j.$$
The remaining  $2^k -(k+1)$ vertices of the Bianchi $(P_1, \ldots, P_k)$-cube $({\cal C}_0,...,{\cal C}_k)$, with  ${\cal C}_0=\{f\}$ and ${\cal C}_1=\{f_1,...,f_k\}$, are given by explicit algebraic formulae as follows. Let
$$
\Omega_{ij}=\frac{P_i(1+P_j^2)}{P_i+P_j}\sum_{\ell=1}^ne^{-\left(\frac{P_i+P_j}{P_iP_j}\right)u_\ell}, \;\; 1\leq i,j\leq k,
%\ee
$$
and 
$$\Fes_i=\sum_{j=1}^n (1+P_i i) \gamma^i_je^{iu_j}E_j, \;\;1\leq i\leq k.$$
Denote by  $e_1, \ldots, e_k$  the canonical basis of $\R^k$, by  $e^1, \ldots, e^k$  its dual basis, and for each $\alpha_r=\{i_1,...,i_r\}\in \Lambda_r$ define 
%$\varphi\colon \R^n\to \R^k$  and $\Omega\in \Gamma((\R^k)^*\otimes \R^k)$ by  
$$
\varphi^{\alpha_r}=\sum_{j=1}^r\varphi_{i_j}e_{i_j}, \  \, \Omega^{\alpha_r}=\sum_{\ell, j=1}^r\Omega_{i_\ell i_j}e^{i_j}\otimes e_{i_\ell} \, \,\mbox{and}
\, \, \Fes^{\alpha_r}=\sum_{j=1}^r e^{i_j}\otimes \Fes_{i_j}.$$
Then all elements of ${\cal C}_r$, $2\leq r\leq k$,  are given by
$$ f_{\alpha_r}=f-{\cal \Fes}^{\alpha_r}(\Omega^{\alpha_{r}})^{-1}\varphi^{\alpha_{r}}, \,\,\,\alpha_r=\{i_1,...,i_r\}\in \Lambda_r.$$
Note that, since the trivial solution $(v, h)$ we started with satisfies $\sum_{i=1}^n v_i^2=1$, and the corresponding solutions $(\varphi_i, \gamma^i_1, \ldots, \gamma^i_n)$  of (\ref{ral_{0l}}) given by (\ref{eq:vagamma}) satisfy (\ref{eq:p*}) for  $1\leq i\leq k$, all elements of the Bianchi $(P_1, \ldots, P_k)$-cube we just constructed  are flat $n$-dimensional Lagrangian submanifolds of $\mathbb{C}^n$ that are contained in $\Sf^{2n-1}$,  and hence are the liftings, by the Hopf projection,
%$\pi:\mathbb S^{2n-1}\rightarrow \mathbb C\mathbb P^{n-1}$, 
of flat $(n-1)$-dimensional  Lagrangian submanifolds of $\mathbb C \mathbb P^{n-1}$.

      The preceding examples are those that can be obtained by applying a sequence of scalar  $P^*$-transformations starting from the degenerate submanifold (\ref{eq:degen}), that is, those which can be obtained from (\ref{eq:degen}) by applying a  $P^*$-transformation  for a diagonalizable~$P$.
      %operator $P$. 
      
      For a general  $P\in (\R^k)^*\otimes \R^k$, the  two last equations of (\ref{ral_{0l}b}) become ordinary differential equations  whose solution, with initial conditions $\gamma_i(0)=v_i^0\in \R^k$, $1\leq i\leq n$, is 
      $$\gamma_i(u_i)=\exp (-u_i(P^t)^{-1})v_i^0,\;\;1\leq i\leq n.$$
Let $\Omega^t\in (\R^k)^*\otimes \R^k$ be given by (\ref{eq:solx}) for $C=((P^t)^2+I)\sum_{i=1}^n\gamma_i\gamma_i^tP$. By Corollary \ref{cor:p*}, the corresponding $P^*$-transform of  (\ref{eq:degen}) is
\be\label{eq:transft}
\begin{array}{l}
\tilde f =-ie^{iu_1}E_1+\sum_{j=1}^n(a_j+ib_j)e^{iu_j}E_j\\
\end{array}
\ee
where
\be\label{eq:ajbj}
a_j=\left<\Omega^{-1}P^t\gamma_1,\gamma_j\right>\;\;\mbox{and}\;\;b_j=\left<\Omega^{-1}P^t\gamma_1,P^t\gamma_j\right>.
\ee
For instance, for the endomorphism  $P\in (\R^2)^*\otimes \R^2$ whose matrix with respect to the canonical basis $e_1, e_2$ of $\R^2$ is, say, 
$\left(
\begin{array}{l}
\,\,\,\, 1\,\,\,\,\,\,\,\,\,1\,\\
-1\,\,\,\; \;\,\,1\,
\end{array}\right),$
one obtains, for $v_j^0=e_1$,  $1\leq j\leq n$, that 
$$\gamma_j(u_j)=e^{-u_j/2}(\cos (u_j/2), \sin (u_j/2)),\;\;\;1\leq j\leq n,$$ 
%$\gamma_j(u_j)=e^{-\frac{u_j}{2}}(\cos \frac{u_j}{2}, \sin \frac{u_j}{2}), \;\;1\leq j\leq n.$can compute that 
$$\Omega=\sum_{\ell=1}^n \frac{e^{-u_\ell}}{4}\Omega_\ell,\;\; \mbox{with} \;\;\Omega_\ell=\left(
\begin{array}{l}
\,\,\,\, 3+\cos u_\ell-2\sin u_\ell\,\,\,\,\,\,\,\,\,1+2\cos u_\ell+\sin u_\ell\,\\
-1+2\cos u_\ell+\sin u_\ell\,\,\,\; \;\,\,3-\cos u_\ell+2\sin u_\ell\,
\end{array}\right),$$
%\newpage 
$$a_j=\frac{4}{5}D^{-1}e^{-U_j}\left((2\cos V_j-4\sin V_j)G +\sum_{\ell=1}^n e^{-u_\ell}(\sin (u_\ell-U_j)-3\cos (u_\ell-U_j))\right)$$
and
$$b_j=\frac{8}{5}D^{-1}e^{-U_j}\left((3\cos V_j-\sin V_j)G-\sum_{\ell=1}^n e^{-u_\ell}(\sin (u_\ell-U_j)+2\cos (u_\ell-U_j))\right),$$
where $U_j=(u_1+u_j)/2$, $V_j=(u_1-u_j)/2$, $G=\sum_{\ell=1}^n e^{-u_\ell}$ and 
%$$G=\sum_{\ell=1}^n e^{-u_\ell}\;\;\mbox{and} D=5(2G^2-\sum_{\ell,j=1}^n e^{-(u_\ell+u_j)}\cos(u_\ell-u_j).$$
$$D=2G^2-\sum_{\ell,j=1}^n e^{-(u_\ell+u_j)}\cos(u_\ell-u_j).$$
Since $P$ is non-diagonalizable, this provides an example of a  $P^*$-transform of (\ref{eq:degen}) that can not be obtained by an iteration of scalar $P^*$-transformations. 

In a similar way one can produce examples of isometric immersions $f\colon M^n(c)\rightarrow \mathbb Q^{n+p}(\tilde c)$ satisfying the assumptions of 
either of Propositions  \ref{thm:hiebetaij} or \ref{teo:curvdif} by applying the $L$-transformation to the trivial solution 
(\ref{eq:trivialsol}) of either of systems (\ref{alphasist}) or (\ref{eq:ii}), according to whether $c=\tilde c$ or $c\neq \tilde c$, respectively. 
One can also obtain examples of horizontal  isometric immersions $f\colon M^n(c)\to \mathbb{S}_\epsilon^{2n+1}(c)$ by applying the $P$-transformation to the trivial solution (\ref{eq:trivialsol}) of system (\ref{eq:iii2}).

{\renewcommand{\baselinestretch}{1}

\hspace*{-20ex}\begin{tabbing} \indent\= Universidade Federal de Mato Grosso
\indent\indent\=  Universidade Federal de S\~{a}o Carlos\\
\>Avenida Valdon Varj\~ao 6390 
\>Via Washington Luiz km 235 \\
\> 78600-000 -- Barra do Gar\c cas -- Brazil
\> 13565-905 -- S\~{a}o Carlos -- Brazil \\
%\> 22460-320 -- Rio de Janeiro -- Brazil  \>
%13565-905 -- S\~{a}o Carlos -- Brazil \\
\> E-mail: dasigu@ufmt.br \> E-mail: tojeiro@dm.ufscar.br
\end{tabbing}}


\begin{thebibliography}{bbbb}
\bibitem{bi}  L. Bianchi,  ``Lezioni di Geometria
Differenziale", Bologna, 1927.



\bibitem{dft} M. Dajczer, L. Florit and R. Tojeiro,
{\it The vectorial Ribaucour transformation for submanifolds and
applications.} Trans. Amer. Math. Soc. {\bf 359} (2007),  4977--4997.

\bibitem{dt1}  M. Dajczer and R. Tojeiro,
{\it Flat totally real submanifolds of $\mathbb C\mathbb P^n$ and the symmetric generalized wave equation.}
Tohoku Math. J. (2) {\bf 47} (1995), no. 1, 117--123.

\bibitem{dt2}  M. Dajczer and R. Tojeiro,
{\it Isometric immersions and the generalized Laplace and elliptic sinh-Gordon equations.}
J. Reine Angew. Math. {\bf 467} (1995), 109--147.

\bibitem{dt3}  M. Dajczer and R. Tojeiro,
{\it The Ribaucour transformation for flat Lagrangian submanifolds.}
J. Geom. Anal. {\bf 10} (2000), no. 2, 269--280.

\bibitem{dt4}  M. Dajczer and R. Tojeiro,
{\it An extension of the classical Ribaucour transformation.}
Proc. London Math. Soc. {\bf 85} (2002), no. 1, 211--232.


\bibitem{dt5} M. Dajczer and R. Tojeiro,
{\it Commuting Codazzi tensors and the Ribaucour transformations
for submanifolds.}  Result. Math. {\bf 44} (2003), no. 3-4, 258--278.


\bibitem{hc} Q. Hu and D. Cheng, {\it The polynomial solution to the Sylvester matrix equation.} Appl. Math. Lett. {\bf 19} (2006), no. 9, 859--864.


%\bibitem[{\bf J}]{j} Je{\v{z}}ek, J. {\it Symmetric matrix polynomial equations.} Kybernetika (Prague) {\bf 22} (1986), no. 1, 19--30. 


\bibitem{K} V. Ku{\v{c}}era, {\it The matrix equation $AX+XB=C^*$.} SIAM. J. Appl. Math. {\bf 26} (1974), 15--25.


\bibitem{lm} Q. P. Liu and M. Manas,
{\it Vectorial Ribaucour transformation for the Lam\'e equations.}
J. Phys. A. {\bf 31} (1998), 193-200.

\bibitem{M} E. Ma, {\it A finite series solution of the matrix equation $AX-XB=C$.} SIAM J. Appl. Math. {\bf 14} (1966), 490--495.

\bibitem{mo} J. D. Moore, {\it Submanifolds of constant
positive curvature I.} Duke Math. J. {\bf 44} (1977), 449--484.

%\bibitem{M} E. Ma. {\it A finite series solution of the matrix equation $AX-XB=C$.} SIAM J. Appl. Math. {\bf 14} (1966), 490--495.

%\bibitem[{\bf M}]{sm} S. Mardare,
%{\it On Pfaff systems with $L^p$ coefficients and their applications in differential geometry.}
%J. Math. Pures Appl. (9) {\bf 84} (2005), no 12, 2659--1692.

%\bibitem[{\bf Re$_1$}]{re1} RECKZIEGEL, H.: Horizontal lifts of \iis into the bundle space of a pseudo-Riemannian submersion. In: {\it Global Diff. Geom. and Global Analysis (1984).\/}  Lecture Notes in Math. 1156, Springer Verlag 1985, pp. 264--279.

\bibitem{re} H. Reckziegel, \emph{A correspondence between horizontal \sus of Sasakian manifolds and totally real \sus of K\"ahlerian manifolds}, Topics in differential geometry, vol. I, II,  North Holland 1988,  pp. 1063--1081.

%\bibitem[{\bf F}]{f} D. Ferus,  {\it A remark on Codazzi
%tensors in constant curvature spaces.} Lect. Notes in Math. 838,
%Berlin 1981, 247.

%\bibitem[{\bf Fe}]{fe} E. Ferapontov,
%{\it Surfaces with flat normal bundle: an explicit construction.}
%Diff. Geom. Appl. {\bf 14} (2001), 15-37.

%\bibitem[{\bf GT}]{gt} E. Ghanza and S. Tsarev, {\it On superposition of the auto B\"acklund
%transformations for $(2+1)$-dimensional integrable systems\/},
%Preprint solv-int/$9606003$, $1996$.

%\bibitem[{\bf HJ}]{hj} U. Hertrich-Jeromin, {\it Introduction to M\"obius differential
%geometry,\/} London Mathematical Society Lecture Note Series, vol.
%{\bf 300}, Cambridge University Press, Cambridge, $2003$.



\bibitem{TZ} O. Taussky and H. Zassenhaus, {\it On the similarity transformation between a matrix and its transpose.} Pacific J. Math. {\bf 9} (1959), 893--896.


\bibitem{tw} C-L. Terng and E. Wang, \emph{Transformation of flat Lagrangian immersions and Egoroff nets.} 
Asian J. Math. {\bf 12} (2008), 99--119.

\bibitem{t} R. Tojeiro,
{\it Lagrangian submanifolds of constant sectional curvature and their Ribaucour transformation.}  Bull. Belgian Math. Soc. - Simon Stevin {\bf 8} (2001),  29--46.




%\bibitem[{\bf Pi}]{pi} {\it U. Pinkall},  Dupin Hypersurfaces, Math.
%Ann. {\bf 270} (1985), 427--440.







\end{thebibliography}
\end{document}